%% file: paper-sinum.tex
\documentclass[review,hidelinks,onefignum,onetabnum]{siamart220329}

\usepackage{comment}
\usepackage[utf8]{inputenc} 
\usepackage[T1]{fontenc}    
\usepackage{lmodern} 
\usepackage[T1]{fontenc} 
\usepackage{hyperref}       
\usepackage{url}            
\usepackage{booktabs}       
\usepackage{lipsum}
\usepackage{amsfonts} 
\usepackage{multirow}
\usepackage{nicefrac} 
\usepackage{mathtools}
\DeclarePairedDelimiter\abs{\lvert}{\rvert}
\DeclarePairedDelimiter\norm{\lVert}{\rVert}
\usepackage{microtype}      
\usepackage{xcolor}
\usepackage{graphicx}
\usepackage{subcaption}
\usepackage{multicol}
\graphicspath{ {./figs/} }
\usepackage{stmaryrd}
\usepackage{bm}
\usepackage[normalem]{ulem}
\usepackage[tickmarkheight=0.1cm]{todonotes}
\usepackage[toc]{appendix}
\usepackage{enumerate}
\usepackage{enumitem}
\usepackage{stackengine}
\usepackage{orcidlink}

\numberwithin{equation}{section}
\nolinenumbers

\newenvironment{nalign}{
	\begin{equation}
		\begin{aligned}
		}{
		\end{aligned}
	\end{equation}
	\ignorespacesafterend
}

\newcounter{example} 
\renewcommand{\theexample}{\arabic{example}} 

\newenvironment{example}[1][]{%
    \refstepcounter{example} 
    \noindent\textbf{Example \theexample.} 
    \ifx&#1&%
    \else \label{#1}%
    \fi
}{\par}

\hypersetup{
	colorlinks=true, 
	linkcolor=blue, 
	urlcolor=red, 
	citecolor=magenta,
	linktoc=all 
}

\newcommand{\jinz}{{j \in \mathbb{Z}}}
\newcommand{\spc}{{\hspace{0.5cm}}}

\newcommand{\jmh}{{j-\frac{1}{2}}}
\newcommand{\jmtbt}{{j-\frac{3}{2}}}
\newcommand{\jph}{{j+\frac{1}{2}}}
\newcommand{\jpo}{{j+1}}
\newcommand{\nph}{{n+\frac{1}{2}}}
\newcommand{\half}{\frac{1}{2}}

\newcommand*\dif{\mathop{}\!\mathrm{d}}

\newcommand{\cblue}[1]{\textcolor{black}{#1}}

\newcommand{\ubar}[1]{\underline{\smash{#1}}}

\DeclareMathOperator{\sign}{sign}

\newenvironment{frcseries}{\fontfamily{frc}\selectfont}{}

\newcommand\barbelow[1]{\stackunder[1.0pt]{$#1$}{\rule{.8ex}{.075ex}}} 

\input{ex_shared}

\ifpdf
\hypersetup{
	pdftitle={Convergence of a second-order central scheme for conservation laws with discontinuous flux},
	pdfauthor={Nikhil Manoj \orcidlink{0009-0003-6812-8633} and
		Sudarshan~Kumar~K. \orcidlink{0009-0005-1186-4167}}
}
\fi

\begin{document}
\maketitle

\begin{abstract}
In this article, we propose  a second-order central scheme of the 
Nessyahu-Tadmor-type  for a class of scalar conservation laws with discontinuous flux and present its convergence analysis. Since solutions to problems with discontinuous flux generally do not belong to the space of bounded variation  (BV), we employ the theory of compensated compactness to establish the convergence of approximate solutions. A major component of our analysis involves deriving the maximum principle and showing the  $\mathrm{W}^{-1,2}_{\mathrm{loc}}$ compactness of a sequence constructed from approximate solutions. The latter is achieved through the derivation of several essential estimates on the approximate solutions. Furthermore, by incorporating a mesh-dependent correction term in the slope limiter,  we show that the numerical solutions generated by the proposed second-order scheme converge to the entropy solution. Finally, we validate our theoretical results by presenting numerical examples.

\end{abstract}


\section{Introduction}
We are interested in the second-order discretization of the Cauchy problem for scalar conservation laws with spatially varying flux:
\begin{nalign}\label{eq:problem}
    u_{t} + f(k(x), u)_{x} &= 0 \quad \mbox{for} \,\, (x,t) \in \mathbb{R}\times \mathbb{R}_{+}, \\
    u(x,0) &= u_{0}(x), \quad \mbox{for} \,\, x \in \mathbb{R},
\end{nalign}
where $t$  and $x$ are the time and space variables, respectively and $u = u(x,t)$ is the unknown quantity. Here, the coefficient $k(x)$ in the flux function $f$ is allowed to be a discontinuous function of the spatial variable $x$. 
\subsection{A brief review of conservation laws with discontinuous flux}
Conservation laws of the form \eqref{eq:problem} with discontinuous coefficient $k(x)$ appear in various physical models, including two-phase flow in porous media arising in oil reservoirs \cite{gimse1992,kaasschieter1999}, sedimentation applications \cite{burger2005,diehl1996} and traffic flow with varying road surface conditions \cite{burger2003,holden1995,mochon1987}, among others. These equations have been the focus of extensive theoretical and numerical studies over the past few decades, for a detailed review see \cite{burger2008a,mishra2017rev} and references therein. 
A classic example often examined in this context is the two-flux case, which emerges by specifying the flux $f$ in \eqref{eq:problem} as follows:
\begin{nalign}\label{eq:twofluxcase}
    f(H(x),u(x,t)) &= H(x)f_{r}(u) + (1-H(x))f_{l}(u) = \begin{cases}
        f_{l}(u) \quad \mbox{for} \, x < 0,\\
        f_{r}(u) \quad \mbox{for} \, x \geq 0,
    \end{cases}
\end{nalign}
where $H$ is the Heaviside function.
Another common example involves the flux function having a multiplicative dependency on a discontinuous coefficient, as given by:
\begin{nalign}\label{eq:multipl_flux}
    f(k(x),u) = k(x)g(u),
\end{nalign} for some appropriate function $g.$
In this work, we focus on equation \eqref{eq:problem} with a general flux function \( f(k(x),u) \) that does not necessarily have to conform to the cases mentioned above.
\par
It is important to note that standard theoretical tools and numerical methods are not applicable to \eqref{eq:problem} when the coefficient $k$ is discontinuous. This situation necessitates the development of novel theoretical approaches and specialized numerical methods. This topic has been extensively explored in the literature, where well-posedness is primarily established through suitable numerical approximations; see, for example, \cite{adimurthi2005,bressan2019,  karlsen2003, karlsen2004, klingenberg1995}. For a comprehensive framework on the well-posedness of \eqref{eq:problem}, we refer to \cite{boris2011}.
A fundamental aspect in this context is the uniqueness of solutions to these problems, where an entropy condition (or admissibility condition) plays a crucial role in ensuring uniqueness by identifying physically relevant solutions. Several attempts have been made in the literature to distill the physically relevant solution. In particular, the authors in \cite{adimurthi2005} considered the two-flux case described in \eqref{eq:twofluxcase} and introduced an infinite family of entropy conditions, referred to as the (A, B) entropy conditions. Each condition within this family ensures uniqueness, giving rise to the concept of (A, B) entropy solutions. This has led to the fact that there is no unique physically relevant solution for scalar conservation laws with discontinuous flux. Instead, the physics of the problems determines the underlying unique solution, which corresponds to a particular (A, B) connection. 
In this work, we focus on an entropy condition given by a particular (A, B) connection for each problem, which is typically expressed through a Kruzkov-type entropy condition or adapted entropy formulation, mainly in the framework of \cite{ burger2005, burger2009, karlsen2004, karlsen2017}.  Such solutions can be identified as the unique vanishing viscosity solution, a concept extensively studied in numerous works \cite{ andreianov2015,bachmann2006, diehl1996,   karlsen2003,  panov2010, seguin2003, towers2000} and remains relevant across various models that incorporate discontinuous flux.
\par
 A wide range of numerical techniques have been developed in the literature for approximating conservation laws with discontinuous flux. While the list is extensive, we mention a few approaches: Godunov-type schemes  \cite{adimurthi2005, adimurthi2007, Andreianov2012,   ghoshal2022b,  karlsen2017}, relaxation schemes \cite{karlsen2004b},  Enguist-Osher schemes \cite{burger2009, towers2000} and upstream mobility schemes~\cite{mishra2010}. Other approaches include DFLU flux \cite{adimurthi2013}, Roe-type schemes \cite{wang2018}, and Monte-Carlo methods for random conservation laws with discontinuous coefficients \cite{badwaik2021}. Furthermore, general monotone (A, B) entropy stable schemes and a local Lax-Friedrichs scheme were analyzed in \cite{adimurthi2014} and  \cite{Xia2020}, respectively. The problem described in \eqref{eq:problem} has also been investigated from various other perspectives.  Recent contributions include regularity results established in \cite{ghoshal2024a}, error estimates for monotone flux functions in \cite{badwaik2020}, flux stability under similar assumptions in  \cite{ruf2022}, compactness estimates in \cite{karlsen2024}, viscosity approximation in \cite{karlsen2025}, and BV regularity of adapted entropy solutions in \cite{ghoshal2024}, among others. Despite these advancements, rigorous analytical results for second- and higher-order methods remain limited.

\subsection{Scope and outline of the paper}
Although first-order numerical methods are reliable and aid in ensuring the well-posedness of the underlying problems, second- and higher-order methods offer substantially improved accuracy, particularly for two and three-dimensional problems. Specifically, designing efficient second-order schemes achieves a desirable balance between the excessive diffusion of first-order schemes and the computational complexity of higher-order (third-order and above) schemes. Regarding second-order numerical schemes for scalar conservation laws with discontinuous flux, there has been some progress in the literature; see \cite{adigowdasudarshan2014,burger2010,sudarshan2014}. In particular, the convergence analysis of a class of second-order schemes to a weak solution was studied in \cite{burger2010}. This analysis was later extended in \cite{adigowdasudarshan2014} to accommodate a broader class of numerical fluxes and to establish convergence to the (A, B)-entropy solution. Nonetheless, these studies rely on a non-local limiter algorithm to ensure the scheme is FTVD (flux total variation diminishing). While effective, the limiter algorithm introduces additional computational tasks compared to conventional second-order schemes. This naturally leads to a question: Is it possible to design a relatively simple scheme, such as one based on MUSCL-type spatial reconstruction, and establish its convergence to the entropy solution?  To the best of our knowledge, this remains an open question, as also noted in \cite{adigowdasudarshan2014, mishra2017rev}. In this article, we aim to address this problem by proposing and {analyzing} a comparatively simple second-order scheme, specifically, a variant of the Nessyahu-Tadmor central scheme (see \cite{nessyahu1990}) that employs the minmod limiter for reconstruction. 
\par 
 It is worth noting that while upwind schemes provide higher-resolution solutions, they are often more restrictive due to the need for solving Riemann problems, either exactly or approximately, at mesh interfaces. In contrast, central schemes, such as the Lax-Friedrichs (LF) scheme, offer a significant advantage: they eliminate the need to solve Riemann problems, thereby reducing computational complexity and making them more attractive for certain applications. For problems of the type described by \eqref{eq:problem},   a staggered  Lax-Friedrichs central scheme was analyzed in \cite{karlsen2004} and more recently in \cite{karlsen2024}. 
 The second-order central scheme of Nessyahu and Tadmor, introduced in \cite{nessyahu1990}, can be viewed as an extension of the first-order LF scheme to the case where the flux function is continuous or, equivalently, where the coefficient function $k(x)$ in  \eqref{eq:problem} is a constant.  This scheme has been extensively studied further in the literature; see \cite{arminjon1999, mehmatoglu2012, nessyahu1992, popov2006b, popov2006a}. We also  note that in \cite{berres2004},  a modification of the NT scheme, namely the Kurganov-Tadmor \cite{kurganov2000} scheme, was adapted to handle the discontinuous flux case arising in the modeling of continuous separation of polydisperse mixtures, and its superior performance over a first-order scheme was computationally illustrated. It is the purpose of this work to extend the Nessyahu-Tadmor (NT) scheme and provide a rigorous convergence analysis for problems  of the type \eqref{eq:problem}, involving a discontinuous coefficient $k(x)$. 

\par
A key highlight of this work is the application of the theory of compensated compactness to establish the convergence of the proposed second-order scheme. The compensated compactness framework serves as a powerful tool for proving convergence of numerical schemes, particularly in cases where bounded variation ($\textrm{BV}$) estimates are unavailable.  It is well known that for conservation laws with discontinuous flux function, solutions may fail to possess bounded total variation; see \cite{adimurthi2011, ghoshal2024}. A widely adopted approach that establishes the convergence of numerical schemes in such cases is the singular mapping technique, wherein the images of the approximate solutions under a suitable monotone map are shown to have diminishing total variation. This technique was introduced by Temple in \cite{temple1982} to prove the convergence of the Glimm scheme, specifically applied to a $2\times2$ resonant system of conservation laws for modeling oil displacement by water and polymer in reservoirs (also see \cite{adimurthi2004}). However, applying the singular mapping technique to second-order schemes is challenging due to the difficulties in obtaining a time-continuity estimate in the absence of monotonicity of the scheme. This obstacle was overcome in \cite{burger2010} and \cite{adigowdasudarshan2014} by imposing a FTVD limiter on the second-order scheme, which makes the scheme flux-TVD and consequently proving the time-continuity. Nevertheless, this limiting algorithm may seem somewhat tailored, lacking the straightforwardness of the standard slope-limiter methods. 
\par
In this article, we investigate the convergence of the proposed second-order central scheme using the compensated compactness theory, specifically within Tartar's framework~\cite{ DiPerna1983,tartar1979}. This approach was previously employed in \cite{karlsen2004} to establish the convergence of the Lax-Friedrichs scheme for scalar conservation laws with discontinuous coefficients; see \cite{ dutta2016, karlsen2004b} for further references. A novel component of our analysis involves deriving certain key estimates on the spatial differences of the approximate solutions.  Furthermore, we derive a bound using the concept of one-sided Lipschitz stability studied in \cite{popov2006b}, which plays a significant role in the derivation of the aforementioned estimates.  These estimates collectively establish the $\mathrm{W}^{-1,2}_{\textrm{loc}}$ compactness required in the compensated compactness framework. Finally, along with the derived maximum principle, the compensated compactness theory ensures the existence of a strongly convergent subsequence. Further, using a classical Lax-Wendroff-type argument, we can show that the limit of this subsequence is a weak solution to the problem \eqref{eq:problem}.
\par
To establish the convergence of the proposed scheme to the entropy solution, following the approach in \cite{ adigowdasudarshan2014, gowda2023, vila1988}, we incorporate a mesh-dependent term into the slope limiter. The core strategy involves writing the time-stepping in the proposed second-order scheme in a predictor-corrector form, where the predictor step employs a first-order accurate Lax-Friedrichs time-stepping, and the correction terms guarantee second-order accuracy. In broad terms, the mesh-dependent term in the modified slope limiter ensures that these correction terms vanish as the mesh size approaches zero. Building on this idea, we show that, as the mesh size tends to zero, our second-order scheme converges to the limit of the first-order Lax-Friedrichs scheme, which was shown to be the entropy solution in \cite{karlsen2004}. {We note that the analysis presented in \cite{adigowdasudarshan2014, gowda2023, vila1988} relies on bounded variation (BV) estimates of the approximate solutions. However, in the case of discontinuous flux, such BV estimates are not necessarily available; see \cite{adimurthi2011, ghoshal2024}. In this context, a key novelty of our approach lies in utilizing a weaker estimate, which we derive, to establish the entropy convergence.
\par
We have organized this article as follows. Section \ref{sec:prelims} provides preliminary details related to the problem \eqref{eq:problem}. In Section \ref{section:soscheme}, we present the proposed second-order central scheme. The compensated compactness theory is outlined in Section \ref{section:cc}. A maximum principle is proven for the proposed scheme in Section \ref{section:maxpri}. Section \ref{section:entropyests} is dedicated to deriving a priori estimates. In Section \ref{section:weaksolconv}, using the compensated compactness framework, we prove the convergence of the proposed scheme along a subsequence to a weak solution.  
Section \ref{section:entropy} establishes the entropy convergence of the scheme, via the introduction of a mesh-size dependent term in the minmod slopes. Numerical results are presented in Section \ref{section:numerical} and conclusions are drawn in Section \ref{sec:conclusion}.

\section{Preliminaries}\label{sec:prelims}
\subsection{Notations}
We use the following notations throughout the paper:
$\mathbb{R}_+:=[0,\infty),$ denotes the set of non-negative real numbers. For $a,b \in \mathbb{R},$ denote $\mathcal{I}(a,b):= (\min(a,b), \max(a,b))$ and $\norm{\cdot} :=\norm{\cdot}_{\mathrm{L}^{\infty} }.$ For sequences $\{v_{j}\}_{\jinz}$ and $\{w_{\jph}\}_{\jinz}$ we denote $\Delta v_{\jph}:= v_{j+1}-v_{j}$ and $\Delta w_{j}:= w_{\jph}-w_{\jmh},$ respectively, for $\jinz$. Also for $a \in \mathbb{R},$ denote  $a_{+}= \max\{a,0\}$ and $a_{-}= \min\{a,0\}.$  For $a \in \mathbb{R},$ the greatest integer function is denoted by $\lfloor a \rfloor$, and the sign function by $\mathrm{sgn}(a).$ Finally, $\norm{\cdot}_{BV}$ denotes the total variation semi-norm.

\subsection{Hypotheses}
We assume throughout the paper that the initial data $u_{0}$ is such that
\begin{nalign}
    u_{0} \in \mathrm{L}^{\infty}(\mathbb{R});\,  u_{0}(x) \in [\barbelow{u},\bar{u}] \,\,  \mbox{for}\,\, \mbox{a.e.} \, x \in \mathbb{R},
\end{nalign} for some $\barbelow{u}, \bar{u} \in \mathbb{R}$ such that $\barbelow{u} \leq \bar{u}.$ 
Further, we impose the following assumptions on the coefficient $k$ and the flux function $f.$
\begin{enumerate}[label=(\textbf{H\arabic*})]
    \item \label{hyp:H1} For some $\barbelow{k}, \bar{k} \in \mathbb{R}$ such that $\barbelow{k} \leq \bar{k}$ the function $k$ satisfies \begin{align}\label{eq:kinBV}k \in (\mathrm{BV} \cap \mathrm{L}^\infty)(\mathbb{R}) \quad \mbox{and}\quad \barbelow{k}\leq k(x) \leq \bar{k} \,\, \mbox{for} \, \, \mbox{a.e.} \,\, x \in \mathbb{R}.\end{align} 
    \item \label{hyp:H2} For each fixed $k \in [\barbelow{k}, \bar{k}],$ the map $f(k,\cdot):u \mapsto f(k,u) \in C^{3}[\barbelow{u}, \bar{u}]$ and is strictly convex. Moreover, there exists $\gamma_{1}, \gamma_{2} \geq 0$ such that $0 < \gamma_{1} \leq  f_{uu}(k,u) \leq \gamma_{2}$ for all $u\in [\barbelow{u}, \bar{u}].$
    \item \label{hyp:H3}For each fixed $u \in [\barbelow{u}, \bar{u}],$ the map $f(\cdot, u): k \mapsto f(k,u) \in C^{2}[\barbelow{k}, \bar{k}]$ with $\partial^{2}_{kk}f\equiv 0.$
    \item \label{hyp:H4} The map $f_u: k \mapsto f_{u}(k,u) \in C^{1}[\barbelow{k}, \bar{k}].$
    \item \label{hyp:H5}The flux $f$ satisfies $f(k_{1},\bar{u}) = f(k_{2}, \bar{u})$ and $f(k_{1},\barbelow{u}) = f(k_{2}, \barbelow{u})$ for all $k_1, k_2 \in [\barbelow{k}, \bar{k}].$ For multiplicative flux $f(k,u) = k g(u),$ this assumption reduces to $g(\barbelow{u})=g(\bar{u})= 0.$ 
    \item \label{hyp:H6} The coefficient $k$ is piecewise $C^{1}$ and is discontinuous only at finitely many points, say $D=\{ x_{1}, x_{2}, \dots, x_{M}\}.$
    \item \label{hyp:H7} Crossing condition: For any jump in the coefficient $k$ with the corresponding left and right limits $k_{m}^{-}$ and $k_{m}^{+}$ respectively,
    \begin{align*}
        f(k_{m}^{+},u_1)-f(k_{m}^{-}, u_1)< 0 < f(k_{m}^{+},u_2)-f(k_{m}^{-}, u_2) \implies u_{1} < u_{2},
    \end{align*} for any states $u_1, u_2 \in [\barbelow{u}, \bar{u}].$
\end{enumerate}

Also, throughout this article we denote by $\kappa \in \mathbb{R_{+}}$ a positive constant such that \begin{nalign}\label{eq:kappa}
    \lambda\norm{f_{u}} \leq \kappa.
\end{nalign} 
\begin{remark}\label{remark:concaveflux}
    The analysis presented in this article is carried out under the assumption that the flux function satisfies the strict convexity condition stated in \ref{hyp:H2}. However, we note that an analogous analysis can be performed for strictly concave fluxes as well.
\end{remark}
\subsection{Weak and entropy solutions}
It is well established that the Cauchy problem \eqref{eq:problem} does not generally admit classical solutions, even when the coefficient $k$ and the initial datum $u_0$ are smooth. Instead, solutions to \eqref{eq:problem} are interpreted in the following weak sense.
\begin{definition}\label{defn:weaksoln}(Weak solution) A function $u \in \mathrm{L}^{\infty}(\mathbb{R} \times \mathbb{R_{+}})$ is said to be a weak solution of \eqref{eq:problem} if it satisfies
\begin{nalign}
    \int_{\mathbb{R}} \int_{\mathbb{R}_{+}} \left(u \phi_{t} + f(k(x), u)\phi_{x}\right) \dif t \dif x + \int_{\mathbb{R}}u_{0} \phi(x,0) \dif x = 0,
\end{nalign} for all test functions $\phi \in \mathcal{D}(\mathbb{R}\times\mathbb{R}_{+}).$ 
\end{definition}

Weak solutions in the sense defined above need not be unique and an entropy condition needs to be specified to choose the relevant weak solution.
\begin{definition}\label{defn:entropysoln}(Entropy solution) A function $u \in \mathrm{L}^{\infty}(\mathbb{R} \times \mathbb{R_{+}})$ is called an entropy solution of \eqref{eq:problem} if for all $c \in \mathbb{R},$  
\begin{nalign}\label{eq:entropy_ineq}
    &\int\int_{\mathbb{R} \times \mathbb{R}_+} (|u - c|\phi_t + \mathrm{sgn}(u - c)(f(k, u) - f(k, c))\phi_x) \, \dif x \, \dif t 
+\int_\mathbb{R} |u_0 - c|\phi(x, 0) \, \dif x 
\\& \spc + \int\int_{(\mathbb{R} \setminus D) \times \mathbb{R}_+ } |f(k(x), c)_x|\phi \, \dif x \dif t \\
&\spc + \sum_{m=1}^M \int_0^\infty |f(k^+_m, c) - f(k^-_m, c)|\phi(x_m, t) \, \dif t \geq 0,
\end{nalign} for all non-negative test functions $\phi \in \mathcal{D}(\mathbb{R}\times\mathbb{R}_{+}),$ 
where $D$  is the set of discontinuities of $k$ as given in   \ref{hyp:H6}.
\end{definition}
\begin{remark}
    The uniqueness of the entropy solution, as defined above, was established in \cite{karlsen2004} under the assumption of the crossing condition \ref{hyp:H7}. However, in the two-flux case \eqref{eq:twofluxcase} with unimodal fluxes and a single flux crossing, this assumption can be omitted (see \cite{burger2009}). We note that our analysis can be extended to the setting of \cite{burger2009} also.    
\end{remark}

\section{Second-order central scheme}\label{section:soscheme}
The spatial domain is discretized using a uniform mesh of size $\Delta x$ into intervals of the form $[x_{\jmh}, x_\jph]$ where $x_{\jph}-x_{\jmh} = \Delta x.$ For a fixed time step $\Delta t,$ the time domain is discretized into points $t^{n}=n\Delta t$ for $n \in \{0,1, \dots \}.$ The ratio $\lambda = \frac{\Delta t}{\Delta x}$ is kept  as a constant throughout. Further, the initial data $u_{0}$ is discretized as 
\begin{align*}
    u_{j}^{0} = \frac{1}{\Delta x}\int_{x_\jmh}^{x_{\jph}} u_{0}(x) \dif x \quad \mbox{for} \,\, j \in \mathbb{Z}.
\end{align*} 
Given the solution $\{u^n_j\}_{\jinz}$ at the time-level $t^{n},$ we compute the  solutions $\{u^{n+1}_\jph\}_{\jinz},$  at the next time level $t^{n+1}$ on the staggered grid following the approach outlined in \cite{nessyahu1990}.
We begin with a piecewise linear reconstruction of the cell averages $\{u_{j}^{n}\}$ utilizing a minmod slope limiter, as a step towards achieving second-order accuracy in space. At the time-level $t^n,$ the solution in each cell $[x_{\jmh},x_{\jph})$ is then reconstructed as follows
\begin{align}\label{eq:reconstruction}
    \tilde{u}_{j}^{n}(x) = u_{j}^{n} + \frac{(x-x_j)}{\Delta x}\sigma_j^{n}, \quad x_{\jph}\le x \le x_{\jph},
\end{align}
where the slopes are given by 
\begin{align}\label{eq:slopes}
    \sigma_{j}^n =\mbox{ minmod}\left(  ( u^{n}_{j+1}-u^{n}_{j}), \half( u^{n}_{j+1}-u^{n}_{j-1}), ( u^{n}_{j}- u^{n}_{j-1}) \right),
\end{align}
the minmod function is  defined by
\begin{align*}
\textrm{minmod}(a_{1}, \cdots , a_{m}) \coloneqq \begin{cases}\mathrm{sgn}(a_{1}) \min\limits_{1 \leq k \leq  m}\{\lvert a_{k}\rvert\} \hspace{0.55cm} \mbox{if} \ \mathrm{sgn}(a_{1}) = \cdots = \mathrm{sgn}(a_{m})\\
0 \hspace{3.45cm}\mbox{otherwise.}
\end{cases}
\end{align*}
Now, a finite volume integration in the domain $[x_{\jmh}, x_{\jph}] \times [t^n, t^{n+1}]$ yields 
\begin{nalign}\label{eq:fvnt}
    \int_{x_{j}}^{x_{j+1}}u(x,t^{n+1})\dif x&= 
    \int_{x_{j}}^{x_{j+1}}u(x,t^{n}) \dif x - \int_{t^n}^{t^{n+1}}f(k_{j+1}, u(x_{j+1},t)) \dif t\\& \spc+\int_{t^n}^{t^{n+1}}f(k_{j}, u(x_{j},t)) \dif t \\
    & \approx 
    \int_{x_{j}}^{x_{\jph}}\tilde{u}^n_{j}(x,t^{n}) \dif x + \int_{x_{\jph}}^{x_{\jpo}}\tilde{u}^n_{\jpo}(x,t^{n}) \dif x \\&\spc - \Delta t \left(f(k_{j+1}, u(x_{j+1},t^{\nph}))-f(k_{j}, u(x_{j},t^{\nph})) \right),
\end{nalign}
using the midpoint quadrature rule in the time integration.
Applying the definition \eqref{eq:reconstruction} in \eqref{eq:fvnt}, we now write the staggered second-order scheme as \begin{align}\label{eq:NTscheme}
    u^{n+1}_{\jph} = \half( u^n_{j}+u^n_{j+1}) - \frac{1}{8}(\sigma_{j+1}^{n} -\sigma_{j}^{n}) - \lambda \left(f(k_{j+1}, u^{n+1/2}_{j+1})-f(k_{j}, u^{n+1/2}_{j})\right), 
\end{align}
where the mid-time step values are computed as
\begin{nalign}\label{eq:mid-time_values}
    u^{n+\half}_j = u_j^n- \frac{\Delta t}{2\Delta x} f_{u}(k_{j}, u_j^{n})\sigma_j^{n}, \quad \jinz.
\end{nalign}
Alternatively, we can also express the scheme \eqref{eq:NTscheme} as 
\begin{align}
     u^{n+1}_{\jph} = \half( u^n_{j}+u^n_{j+1}) -\lambda \left(g(k_{j+1},u_{j+1}^{n})- g(k_{j},u_{j}^{n})\right),
\end{align}
where $g$ is given by
\begin{nalign}\label{eq:eq:g_defn}
    g(k_{j},u^{n}_{j}) := f(k_{j},u_{j}^{n+\frac{1}{2}})+ \frac{1}{8\lambda}\sigma_{j}^{n}.
\end{nalign}
For a fixed mesh size $\Delta x,$ the piecewise constant approximate solution and the discretized discontinuous coefficient $k$ are represented by the pair 
\begin{nalign}\label{eq:pcsoln}
    &(u_{\Delta}(x,t), k_{\Delta}(x,t)) \\&\spc := \begin{cases}
    (u_{j}^{n}, k_j) \,\, \quad \quad\mbox{if} \, n \, \mbox{is even and} \, (x,t) \in [x_\jmh,x_\jph) \times [t^{n}, t^{n+1}),  \\
    (u_{\jph}^{n},k_{\jph}) \, \mbox{if} \, n \, \mbox{is odd and} \,\, (x,t) \in [x_{j},x_{j+1}) \times [t^{n}, t^{n+1}),
    \end{cases}
\end{nalign} where $n\in \mathbb{N}\cup\{0\}$ and 
\begin{nalign}
    k_{j}:= \frac{1}{\Delta x}\int_{x_{\jmh}}^{x_{\jph}}k(x) \dif x \quad \mbox{and} \quad k_{\jph}:= \frac{1}{\Delta x}\int_{x_{j}}^{x_{\jpo}}k(x) \dif x. 
\end{nalign}

\begin{remark}
    When the slopes $\sigma_{j} = 0$ for all $j,$ the scheme reduces to the first-order Lax-Friedrichs scheme 
\begin{nalign}\label{eq:lxfscheme}
         u^{n+1}_{\jph} = \half( u^n_{j}+u^n_{j+1}) -\lambda \left(f(k_{j+1},u_{j+1}^{n})- f(k_{j},u_{j}^{n})\right),
    \end{nalign} considered in \cite{karlsen2004}. In other words, the scheme \eqref{eq:NTscheme} is a second-order extension of the scheme \eqref{eq:lxfscheme}.
\end{remark}
\section{Compensated compactness}\label{section:cc}
The main ingredient of the compensated compactness framework is a theorem given below, the proof of which can be found in \cite{karlsen2004}. In the sequel, we will show that the proposed scheme \eqref{eq:NTscheme} meets the requirements outlined in this theorem. This will ensure the existence of a convergent subsequence of approximate solutions generated from the scheme \eqref{eq:NTscheme}.
\begin{theorem}[Compensated compactness theorem]\label{lemma:compcompactness}
    Assume that the hypotheses \ref{hyp:H1}-\ref{hyp:H5} hold true. Let $\{u^\varepsilon\}_{\varepsilon > 0}$ be a sequence of measurable functions defined on $\mathbb{R} \times \mathbb{R}^+$ that satisfies the following two conditions: 
    \begin{enumerate}
        \item There exist $a, b\in \mathbb{R} $ with $a < b$, both independent of $\varepsilon$, such that
        \[
            a \leq u^\varepsilon(x, t) \leq b \text{ for a.e. } (x, t) \in \mathbb{R} \times \mathbb{R}^+.
        \]
        \item The two sequences 
        \[
            \{ S_1(u^\varepsilon)_t + Q_1(k(x), u^\varepsilon)_x \}_{\varepsilon > 0} \quad \mbox{and} \quad \{ S_2(k(x), u^\varepsilon)_t + Q_2(k(x), u^\varepsilon)_x \}_{\varepsilon > 0},
        \]
        belong to a compact subset of $\mathrm{W}^{-1,2}_{\mathrm{loc}}(\mathbb{R} \times \mathbb{R}^+)$, where
        \begin{nalign}\label{eq:entropypairs}
             S_1(u) &:= u - c, \quad Q_1(k,u) := f(k,u) - f(k,c),\\
             S_2(k,u) &:= f(k,u) - f(k,c) \quad \mbox{and} \quad Q_2(k,u) := \int_c^u (f_u(k,\xi))^2 \, \dif \xi,
        \end{nalign}
        for any $c \in \mathbb{R}$.
    \end{enumerate}
    Then, there exists a subsequence of $\{u^\varepsilon\}_{\varepsilon > 0}$ that converges pointwise a.e. to a function $u \in \mathrm{L}^{\infty}(\mathbb{R} \times \mathbb{R}^+).$
\end{theorem}
To establish the $\mathrm{W}^{-1,2}_{\mathrm{loc}}$ compactness required by Theorem \ref{lemma:compcompactness}, we will make use of the following interpolation result as well.  For a proof of this result, see \cite{ding1985}. 
\begin{lemma}\label{lemma:interpolation}
Let $\Omega \subset \mathbb{R}^d$ be a bounded open set. Let $q$ and $r$ be a pair of constants satisfying $1 < q < 2 < r < \infty$. If $A$ is compact subset of  $\mathrm{W}^{-1,q}_{\mathrm{loc}}(\Omega)$ and $B$ is a bounded subset of $\mathrm{W}^{-1,r}_{\mathrm{loc}}(\Omega),$ then \\
$$A \cap B \, \text{is compact in}
\,\mathrm{W}^{-1,2}_{\mathrm{loc}}(\Omega).$$
\end{lemma}
\section{Maximum principle and \texorpdfstring{$\mathrm{L}^{\infty}$}{L-infinity}-stability}\label{section:maxpri}
This section establishes that the approximate solutions $u_{\Delta}$, as defined in \eqref{eq:pcsoln}, satisfy a global maximum principle, thereby yielding an $\mathrm{L}^{\infty}$-estimate as well.
\begin{theorem}\label{thm:maxpri}
Let the initial datum $u_{0} \in \mathrm{L}^{\infty}(\mathbb{R})$ with $\barbelow{u} \leq u_{0}(x) \leq \bar{u},$  for all $x \in \mathbb{R}.$ Then, under the CFL condition
\begin{nalign}\label{eq:cfl_maxpri}
    \lambda\norm{f_{u}} \leq \kappa \leq \frac{\sqrt{2}-1}{2},
\end{nalign} and hypotheses \ref{hyp:H1}-\ref{hyp:H5}, the approximate solution $u_{\Delta}$ \eqref{eq:pcsoln} obtained from the scheme \eqref{eq:NTscheme} satisfies the global maximum principle 
\begin{nalign}\label{eq:maxpr}
    \ubar{u} \leq u_{\Delta}(x,t) \leq \bar{u},
\end{nalign} for all $(x,t) \in \mathbb{R}\times\mathbb{R}_{+}.$
Consequently, the approximate solutions $u_{\Delta}$ are $\mathrm{L}^{\infty}$-stable, i.e., 
\begin{nalign}\label{Linf_bd}
    \norm{u_{\Delta}} \leq C_{u_{0}}:=\max\{ \abs{\ubar{u}},\abs{\bar{u}}\}.
\end{nalign}
\begin{proof}
We use the principle of mathematical induction to prove this result. By the assumption on the initial datum $u_0$, the result holds for $n=0.$ For $n\geq 0,$ suppose $u_{j}^{n} \in [\barbelow{u}, \bar{u}],$ for all $\jinz.$ We will now prove  that $u_{\jph}^{n+1} \in [\barbelow{u}, \bar{u}]$ for all $\jinz.$
    Adding and subtracting the term $f(k_{\jpo}, \bar{u}),$ and using the hypothesis \ref{hyp:H5} (i.e., $f(k_{\jpo}, \bar{u})= f(k_{j}, \bar{u})$), the difference $f(k_{\jpo}, u_{\jpo}^{\nph})- f(k_{j}, u_{j}^{\nph})$ can be expressed as 
    \begin{nalign}\label{eq:fdif}
        f(k_{\jpo}, u_{\jpo}^{\nph})- f(k_{j}, u_{j}^{\nph}) & = f(k_{\jpo}, u_{\jpo}^{\nph})- f(k_{\jpo}, \bar{u})+ f(k_{j}, \bar{u})-   f(k_{j}, u_{j}^{\nph})\\
        &= f_{u}(k_{\jpo}, \tilde{\zeta}_{\jpo})(u_{\jpo}^{\nph} - \bar{u}) + f_{u}(k_{j}, \tilde{\zeta}_{j})(\bar{u}- u_{j}^{\nph} )\\
        & = f_{u}(k_{\jpo}, \tilde{\zeta}_{\jpo})(u_{\jpo}^{n} - \bar{u}) + f_{u}(k_{j}, \tilde{\zeta}_{j})(\bar{u}- u_{j}^{n}) \\ & \spc -\frac{1}{2}\lambda f_{u}(k_{\jpo}, u_{\jpo}^{n}) f_{u}(k_{\jpo}, \tilde{\zeta}_{\jpo})\sigma_{\jpo}^{n}  \\& \spc +\frac{1}{2}\lambda f_{u}(k_{j}, u_{j}^{n}) f_{u}(k_{j}, \tilde{\zeta}_{j})\sigma_{j}^{n},
    \end{nalign}
where $\tilde{\zeta}_{\jpo} \in \mathcal{I}(u_{\jpo}^{\nph}, \bar{u})$ and $\tilde{\zeta}_{j} \in \mathcal{I}(u_{j}^{\nph}, \bar{u}).$ From the expression \eqref{eq:fdif}, we get the estimate
\begin{nalign}\label{eq:fluxdifest}
    \abs[\big]{ f(k_{\jpo}, u_{\jpo}^{\nph})- f(k_{j}, u_{j}^{\nph})} \leq \norm{f_{u}} \left(2\bar{u} - u_{j}^{n} - u_{j+1}^{n}\right) + \lambda (\norm{f_u})^2 \abs{u_{j+1}^n -u_j^n}.
\end{nalign}
\par 
From the CFL condition  \eqref{eq:cfl_maxpri}, we obtain $\displaystyle \frac{1}{4}+\kappa^{2} \leq \frac{1}{2} - \kappa.$ This inequality, when combined with the estimate \eqref{eq:fluxdifest} applied on the scheme \eqref{eq:NTscheme}, and \eqref{eq:kappa}, yields
\begin{nalign}\label{uub}
     u^{n+1}_{\jph} &\leq \half( u^n_{j}+u^n_{j+1}) + \frac{1}{4}\abs{u_{j+1}^{n}-u_{j}^{n}} + \lambda \norm{f_{u}} \left(2\bar{u} - u_{j}^{n} - u_{j+1}^{n}\right) \\ & \spc +\lambda^{2} \norm{f_u}^2 \abs{u_{j+1}^n -u_j^n}\\
     & \leq \left(\frac{1}{2}-\kappa\right)( u^n_{j}+u^n_{j+1}) +  \left(\frac{1}{4}+\kappa^2\right)\abs{u_{j+1}^{n}-u_{j}^{n}} +2\kappa \bar{u}\\
     & \leq \left(\frac{1}{2}-\kappa\right) 2\max\{u^n_{j},u^n_{j+1}\} + 2\kappa \bar{u}= \left(1-2\kappa\right) \max\{u^n_{j},u^n_{j+1}\} + 2\kappa \bar{u} \leq \bar{u}.
\end{nalign}
By the CFL condition \eqref{eq:cfl_maxpri}, $2\kappa \in [0,1]$ and hence the term $\left(1-2\kappa\right) \max\{u^n_{j},u^n_{j+1}\}\\ + 2\kappa \bar{u}$ is a convex combination of points in $[\barbelow{u}, \bar{u}].$ As a result, the last inequality in the previous equation holds true. Similar arguments as in \eqref{eq:fdif}, by adding and subtracting the term $f(k_{\jpo}, \barbelow{u})$ in $f(k_{\jpo}, u_{\jpo}^{\nph})- f(k_{j}, u_{j}^{\nph}),$ gives another estimate
\begin{align*}
    \abs[\big]{ f(k_{\jpo}, u_{\jpo}^{\nph})- f(k_{j}, u_{j}^{\nph})} \leq \norm{f_{u}} \left(u_{j}^{n} + u_{j+1}^{n} - 2\barbelow{u}\right) + \lambda (\norm{f_u})^2 \abs{u_{j+1}^n -u_j^n},
\end{align*}
which subsequently yields the lower bound
\begin{align}\label{ulb}
     u^{n+1}_{\jph} &\geq \half( u^n_{j}+u^n_{j+1}) - \frac{1}{4}\abs{u_{j+1}^{n}-u_{j}^{n}} - \lambda \norm{f_{u}} \left(u_{j}^{n} + u_{j+1}^{n}-2\barbelow{u}\right)\notag\\ & \spc - \lambda^{2} \norm{f_u}^2 \abs{u_{j+1}^n -u_j^n} \notag \\
     & \geq \left(\frac{1}{2}-\kappa\right)( u^n_{j}+u^n_{j+1}) - \left(\frac{1}{4}+\kappa^2\right)\abs{u_{j+1}^{n}-u_{j}^{n}} +2\kappa \barbelow{u}\\
       & \geq \left(\frac{1}{2}-\kappa\right)\left(( u^n_{j}+u^n_{j+1}) - \abs{u_{j+1}^{n}-u_{j}^{n}}\right) +2\kappa \barbelow{u} \notag \\ & = \left(1-2\kappa\right)\min\{u^n_{j},u^n_{j+1}\} +2\kappa \barbelow{u} \geq \barbelow{u}. \notag
\end{align}
\par
The expressions \eqref{uub} and \eqref{ulb} together yield the maximum principle \eqref{eq:maxpr}. 
\end{proof}
\end{theorem}   
\section{A priori estimates}\label{section:entropyests}In this section, we derive certain essential a priori estimates on the approximate solutions obtained from the second-order scheme \eqref{eq:NTscheme}. These estimates form the foundation for the convergence analysis presented in Section \ref{section:weaksolconv}.
\par
First, we focus on estimating the term \(\sum_{\jinz}(u_{\jpo}^{n} - u_{j}^{n})^3_{+}\), where \(\{u_{j}^{n}\}_{\jinz}\) are the solutions derived from the second-order scheme \eqref{eq:NTscheme}. The details of this estimate are outlined in the following lemma. The proof is lengthy and is provided in Appendix \ref{app:A}.

\begin{lemma}\label{lemma:osle}
Assume $\{u^{n}_{j}\}_{\jinz}$ is such that $\abs{u_j} \leq C_{u_{0}}, \, \forall \jinz,$ where $C_{u_{0}}$ is as in \eqref{Linf_bd}. Let $\{u^{n+1}_{\jph}\}_{\jinz}$ be obtained from $\{u^{n}_{j}\}_{\jinz}$ by applying the time-update formula \eqref{eq:NTscheme}.
     Under the CFL condition \begin{nalign}\label{eq:cfl}
        \lambda \norm{f_{u}} \leq \kappa \leq \min\left\{\frac{\gamma_{1}}{7500\gamma_{2}}, \frac{1}{4000}\right\},
    \end{nalign} the solution $\{u_{\jph}^{n+1}\}_{\jinz}$ satisfies the estimate \begin{nalign}\label{eq:oslsfinal} \sum_{\jinz}\left(\Delta u_{j}^{n+1}\right)_{+}^{2}
    & \leq \sum_{\jinz}(\Delta u_{\jmh}^{n})_{+}^2- \frac{1}{500} \lambda \gamma_{1}\sum_{\jinz} (\Delta u_{\jmh}^{n})_{+}^{3} + \Psi\norm{k}_{BV},
    \end{nalign}
for all $n \geq 0,$ where \begin{align*}
    \Psi&:=  72\lambda^{2}(C_{u_{0}})^{2}\norm{f_{uk}}+ 114\lambda(C_{u_{0}})^{2}\norm{f_{ku}} + \left(708(C_{u_{0}})^{2}+48\lambda\norm{f_{u}}\right)\lambda^{2}\norm{f_{u}}\norm{f_{uk}}\\
    & \spc + \left(48\lambda^{2} C_{u_0}\norm{k}+132\lambda^{2} (C_{u_0})^{2}\gamma_{2}\norm{f_{u}}+64\lambda\norm{f_k}\norm{k}+88 C_{u_0}\right)\lambda\norm{f_{k}}.
\end{align*}

\end{lemma}

Next, we proceed to obtain a cubic estimate on the spatial differences of the approximate solutions.
\begin{lemma}(Cubic estimate)\label{lemma:cubicest} 
Let the initial datum $u_{0} \in (\mathrm{L}^{\infty}\cap \mathrm{BV})(\mathbb{R})$ be such that $\norm{u_{0}} \leq C_{u_{0}}.$  For any fixed $T>0, X>0,$ define $N:=\lfloor T/\Delta t  \rfloor+1$ and $J:=\lfloor X/\Delta x \rfloor +1.$ Then under the CFL condition
\begin{align}\label{eq:cfl_cubicest}
       \lambda \norm{f_{u}} \leq \kappa \leq \min\left\{\frac{\gamma_{1}}{7500\gamma_{2}}, \frac{1}{4000}, \frac{7}{85+16C_{u_0}}, \frac{\gamma_1}{\gamma_2\chi}\right\},
\end{align}
with $\chi:= 228+13C_{u_0}+ 174C_{u_0}\gamma_2+ 12 (C_{u_0})^{2}\gamma_2,$ the approximate solutions generated by the second-order scheme \eqref{eq:NTscheme} satisfy the uniform bound \begin{align}\label{eq:C_cubiclemma}
       \Delta x\sum_{n=0}^{N-1}\sum_{\substack{\abs{j} \leq J \\ j+\frac{n}{2} \in \mathbb{Z}}}\abs{\Delta u_{\jph}^{n}}^3\leq C(X,T),
    \end{align} for a constant $C(X,T)$ independent of $\Delta x.$    
\begin{proof}
We need to define a linear function $\bar{g}_{\jph}$ which interpolates $g(k_{j}, u_{j}^{n})$ and $g(k_{j+1}, u_{j+1}^{n})$ as follows:

\begin{align*}
    \bar{g}_{\jph}(u) := g(k_{j}, u_{j}^{n}) + \frac{\Delta g_{\jph}^{n}}{\Delta u_{\jph}^n}(u-u_{j}^{n}) \,\, \mbox{for} \,\, u\in \mathcal{I}(u_j^{n}, u_{j+1}^{n}),
\end{align*} 
where $\Delta g_{\jph}^n := g(k_{j+1},u_{j+1}^{n})-g(k_{j},u_{j}^{n}).$ 
Further, we consider functions $S,Q:[\barbelow{k}, \bar{k}] \times [\barbelow{u}, \bar{u}] \rightarrow \mathbb{R},$ with the property that $ \partial_{u} Q= \partial_{u}S  \partial_{u}f.$ Now, we define a quantity  $E_{\jph}^{n}$ associated with the pair $(S,Q)$ as    
\begin{nalign}\label{eq:Ejph} E_{\jph}^{n} &:=  S(k_{\jph},u^{n+1}_{\jph}) - \frac{1}{2}\left( S(k_{j},u_{j}^{n})+S(k_{j+1} ,u_{j+1}^{n})\right) \\ & \spc + \lambda\left(Q(k_{j+1}, u_{j+1}^{n})- Q(k_{j},u_{j}^{n})\right).
\end{nalign}
Our objective now is to reformulate $E_{\jph}^n$ in a suitable form, and use that to obtain the desired estimate. We begin with rewriting $E^{n}_{\jph}$ by adding and subtracting suitable terms, as \begin{align}\label{eq:entr_pro_1}     E_{\jph}^{n}&= S(k_{\jph},u^{n+1}_{\jph}) - \frac{1}{2}\left( S(k_{\jph},u_{j}^{n})+S(k_{\jph} ,u_{j+1}^{n})\right) \\ & \spc + \lambda\left(Q(k_{\jph}, u_{j+1}^{n})- Q(k_{\jph},u_{j}^{n})\right) + R_{1}+ R_{2},\notag
\end{align}
where \begin{align}\label{eq:R1R2}
    R_{1} &:= \frac{1}{2}\left( S(k_{\jph},u_{j}^{n})+S(k_{\jph} ,u_{j+1}^{n})\right)- \frac{1}{2}\left( S(k_{j},u_{j}^{n})+S(k_{j+1} ,u_{j+1}^{n})\right) \\&=\mathcal{O}(k_{\jph}-k_{j}) + \mathcal{O}(k_{\jph} - k_{j+1}),\notag\\
    R_{2} &:= \lambda \left[ Q(k_{j+1}, u_{j+1}^{n})- Q(k_{j},u_{j}^{n}) -\left(Q(k_{\jph}, u_{j+1}^{n})- Q(k_{\jph},u_{j}^{n})\right)\right] \notag\\
    & = \mathcal{O}(k_{j}-k_{\jph}) + \mathcal{O}(k_{j+1}- k_{\jph}). \notag
\end{align}
\par
Now, we define parametrized functions for $s\in[0,1]$ as
\begin{nalign}\label{eq:paramet_fns}
u(s)&:= su_{j}^{n} + (1-s)u_{j+1}^{n},\\
u_{\jph}(s) &:= \frac{1}{2}\left(u(s)+ u_{j+1}^{n}\right) - \lambda(g(k_{j+1},u^{n}_{j+1})- \bar{g}(u(s))),\\
u(r,s) &:=  ru(s) + (1-r)u_{j+1}^{n} \quad \mbox{and}\\
u_{\jph}(r,s) &:= \frac{1}{2}\left( u(s) + u(r,s) \right) - \lambda(\bar{g}(u(r,s))- \bar{g}(u(s))).
\end{nalign}
These functions satisfy the properties
\begin{nalign}
u(0) &= u_{\jpo}^n, \quad u(1)= u_{j}^{n}, \quad
u_{\jph}(0)= u_{\jpo}^{n}, \quad u_{\jph}(1)=  u^{n+1}_{\jph},\\
u(0,s)&= u_{\jpo}^{n}, \quad u(1,s):= u(s), \quad    
    u_{\jph}(0,s) = u_{\jph}(s), \quad u_{\jph}(1,s) = u(s).\\
\end{nalign}
Next, using \eqref{eq:paramet_fns}, we write 
\begin{nalign}\label{eq:Sdif_modified}
    &S(k_{\jph},u^{n+1}_{\jph}) - \frac{1}{2}\left(S(k_{\jph},u_{j}^{n})+S(k_{\jph} ,u_{j+1}^{n})\right) \\&= \int_{0}^{1}S_{u}(k_{\jph}, u_{\jph}(s))\left(\frac{1}{2}+\lambda \bar{g}_{\jph}^\prime(u(s))\right) u^\prime(s) \dif s -\frac{1}{2}\int_{0}^{1}S_{u}(k_{\jph}, u(s))u^\prime(s)\dif s.
\end{nalign}
Also, if $\tilde{Q}$ is defined such that $\displaystyle\tilde{Q}_{u} = S_{u}\bar{g}_{\jph}^\prime = S_{u}\frac{\Delta g_{\jph}^n}{\Delta u_{\jph}^n},$  then adding and subtracting appropriate terms, we get
\begin{nalign}\label{eq:Qdif_modified}
&\lambda\left(Q(k_{\jph}, u_{j+1}^{n}- Q(k_{\jph},u_{j}^{n}))\right) \\&= \lambda(\tilde{Q}(k_{\jph},u_{j+1}^{n})- \tilde{Q}(k_{\jph},u_{j}^{n})) \\&\hspace{0.5cm} + \lambda\left[\left(Q(k_{\jph}, u_{j+1}^{n})-\tilde{Q}(k_{\jph},u_{j+1}^{n})\right)- \left(Q(k_{\jph}, u_{j}^{n})-\tilde{Q}(k_{\jph},u_{j}^{n})\right)\right]\\
    & =  -\lambda\int_{0}^{1}S_{u}(k_{\jph}, u(s))\bar{g}_\jph^\prime(u(s))u^\prime(s) \dif s \\
    &\hspace{0.5cm} -\lambda\int_{0}^{1}S_{u}(k_{\jph}, u(s))\left(f_{u}(k_{\jph}, u(s))- \bar{g}_\jph^\prime(u(s))\right)u^\prime(s) \dif s.
\end{nalign}
Now, in view of \eqref{eq:Sdif_modified} and \eqref{eq:Qdif_modified}, and rearranging  the terms, we write 
\begin{nalign}\label{eq:entropyprod}
    E_{\jph}^{n} 
    & = I+J+ R_{1}^n+ R_{2}^n,
\end{nalign}
where \begin{nalign}\label{eq:IandJ_defn}
I &:= \int_{0}^{1}\left(S_{u}(k_{\jph}, u_{\jph}(s))-S_{u}(k_{\jph}, u(s))\right)\left(\frac{1}{2}+\lambda \bar{g}_{\jph}^\prime(u(s))\right) u^\prime(s) \dif s,\\
J &:= -\lambda\int_{0}^{1}S_{u}(k_{\jph}, u(s))\left(f_{u}(k_{\jph}, u(s))- \bar{g}_{\jph}^\prime(u(s))\right)u^\prime(s) \dif s.
\end{nalign}
Using integration by parts, $J$ can be simplified as 
\begin{nalign}\label{eq:J}
    J 
    &= \bar{J} - \lambda\left[S_{u}(k_{\jph}, u(s))\left(f(k_{\jph}, u(s))- \bar{g}_{\jph}(u(s))\right)\right]_{0}^{1},
\end{nalign}
where $\displaystyle \bar{J}:= \lambda\int_{0}^{1}S_{uu}(k_{\jph},u(s))\left(f(k_{\jph}, u(s))- \bar{g}_{\jph}(u(s))\right) u^\prime(s)\dif s.$
Further, \\by adding and subtracting appropriate terms in \eqref{eq:J}, we obtain
\begin{nalign}\label{eq:Jfinal}
     J &= \bar{J} -\lambda\left[S_{u}(k_{j}, u_j^n)(f(k_{j}, u_j^n) - g(k_{j}, u_{j}^n))\right.\\ & \spc \left.-S_{u}(k_{j+1}, u_{j+1}^n)(f(k_{j+1}, u_{j+1}^n) - g(k_{j+1}, u_{j+1}^n))\right] + \tilde{R}_1^{n} +\tilde{R}_2^{n},
\end{nalign}
where 
\begin{align}\label{eq:R1R2tilde}
    \tilde{R}_1^{n} &:= -\lambda\left[S_{u}(k_\jph, u_j^{n})f(k_{\jph}, u_j^{n})-S_{u}(k_j, u_j)f(k_{j}, u_j^{n})\right]\\
    &\hspace{0.5cm} + \lambda\left[S_{u}(k_\jph, u_{j+1}^{n})f(k_{\jph}, u_{j+1}^{n})-S_{u}(k_{j+1}, u_{j+1}^{n})f(k_{j+1}, u_{j+1}^{n})\right]\notag\\
    &= \mathcal{O}(k_{\jph}-k_j) + \mathcal{O}(k_{\jph}-k_{j+1}), \notag\\
    \tilde{R}_2^{n} &:= -\lambda\left[\left(S_{u}(k_j, u_j^{n})-S_{u}(k_\jph, u_j^{n})\right)g(k_{j}, u_j^{n})\right]\notag\\
    &\hspace{0.5cm} + \lambda\left[\left(S_{u}(k_{j+1}, u_{j+1}^{n})-S_{u}(k_\jph, u_{j+1}^{n})\right)g(k_{j+1}, u_{j+1}^{n})\right]\notag\\
    &= \mathcal{O}(k_{\jph}-k_j) + \mathcal{O}(k_{j+1}-k_{\jph}). \notag
\end{align}
At this stage, equating the right hand sides of  \eqref{eq:Ejph} and \eqref{eq:entropyprod} and using \eqref{eq:Jfinal}, we may write 
\begin{nalign}\label{eq:numentropyprod}
    \bar{E}_{\jph}^{n} 
    &= I+\bar{J} + R_{1}^{n}+ R_{2}^{n} + \tilde{R}_{1}^{n}+ \tilde{R}_{2}^{n},
\end{nalign}
where we define
\begin{align}\label{eq:numentropyflux}\bar{E}_{\jph}^{n} 
   &:=  S(k_{\jph},u^{n+1}_{\jph}) - \frac{1}{2}\left( S(k_{j},u_{j}^{n})+S(k_{j+1} ,u_{j+1}^{n})\right) \\ & \spc + \lambda\left(G(k_{j+1}, u_{j+1}^{n})- G(k_{j},u_{j}^{n}))\right), \notag\\
    G(k_{j}, u_{j}) &:= Q(k_j, u_j^n) - S_{u}(k_{j}, u_j^n)\left(f(k_{j}, u_j^n) - g(k_{j}, u_{j}^n)\right).\notag 
\end{align}
Plugging in the derivative 
$\frac{\partial}{\partial r}u_{\jph}(r,s) = -s\Delta u_{\jph}^{n} \left(\frac{1}{2} -\lambda \bar{g}_{\jph}^\prime(u(r,s))\right)
$ in \eqref{eq:IandJ_defn}, the term $I$ can be represented as
\begin{nalign}\label{eq:I_b}
    I = -(\Delta u_{\jph}^{n})^2 \tilde{I}, \quad \mbox{where}
\end{nalign}  $\tilde{I}:=\int_{0}^{1}\int_{0}^{1}s S_{uu}(k_{\jph}, u_{\jph}(r,s))\left(\frac{1}{2}- \lambda \bar{g}_{\jph}^{\prime}(u(r,s))\right)\left(\frac{1}{2}+ \lambda \bar{g}_{\jph}^{\prime}(u(s))\right) \dif s \dif r.$
\par
Now, we set  $\displaystyle S(k,u)= \frac{u^{2}}{2},$ and re-work on \eqref{eq:numentropyprod}, specifically focusing on $I$ and $\bar{J}.$ Noting $S_{uu}=1,$ performing a change of variable $z=u(s)$ and subsequently applying the trapezoidal rule, the term $\bar{J}$ simplifies to 
\begin{nalign}\label{eq:Jbar_aa}
    \bar{J} 
    &= -\lambda\int_{u_{j}}^{u_{j+1}}f(k_{\jph}, z) \dif z + \lambda \int_{u_{j}}^{u_{j+1}} \bar{g}_{\jph}(z)\dif z\\
    &= -\frac{1}{2}\lambda \Delta u_{\jph} \left[f(k_{\jph}, u_{j}^{n})-g(k_j,u_j^{n})+ f(k_{\jph}, u_{j+1}^{n})-g(k_{j+1},u_{j+1}^{n})\right] \\& \spc + \frac{\lambda}{12}(\Delta u_\jph^{n})^{3} f_{uu}(k_{\jph}, \zeta_{1}),
\end{nalign}
for some $\zeta_{1} \in \mathcal{I}(u_{j},u_{j+1}).$
Recalling the definition $g(k_{j},u^{n}_{j}) := f(k_{j},u_{j}^{n+\frac{1}{2}})+ \frac{1}{8\lambda}\sigma_{j}^{n}$ and using Taylor series expansions in the second variable of $g,$ we obtain
\begin{nalign}\label{eq:gtaylor_j}
g(k_j,u_j^{n}) 
& = f(k_{j},u_j^n)-\frac{\lambda}{2}(a_{j}^{n})^2\sigma_{j}^n+ \frac{1}{8}(\lambda a_{j}^{n}
\sigma_{j}^n)^{2}f_{uu}(k_j,\zeta_{2})+ \frac{1}{8\lambda}\sigma_{j}^n,\\
g(k_{j+1},u_{j+1}^{n}) 
& = f(k_{j+1},u_{j+1}^n)-\frac{\lambda}{2}(a_{j+1}^{n})^2\sigma_{j+1}^n \\ & \spc+ \frac{1}{8}(\lambda a_{j+1}^{n}\sigma_{j+1}^n)^{2}f_{uu}(k_{j+1},\zeta_{3}) + \frac{1}{8\lambda}\sigma_{j+1}^n,
\end{nalign} for $\zeta_{2} \in \mathcal{I}(u_{j}^n, u^{n+\frac{1}{2}}_j)$ and $\zeta_{3} \in \mathcal{I}(u_{j+1}^n, u^{n+\frac{1}{2}}_{j+1}),$
where $a_{j}^{n} = f_{u}(k_j, u_j^{n}),$ as in \eqref{eq:ajn_defn}. By adding and subtracting the terms $f(k_{j},u_j^{n})$ and $f(k_{j+1},u_{j+1}^{n}),$ respectively, and using \eqref{eq:gtaylor_j} we obtain
\begin{nalign}\label{eq:fminusg_taylor}
    f(k_{\jph}, u_{j}^{n})-g(k_j,u_j^{n}) 
    &= \frac{\lambda}{2}(a_{j}^{n})^2\sigma_{j}^n- \frac{1}{8}(\lambda a_{j}^{n}\sigma_{j}^n)^{2}f_{uu}(k_j,\zeta_{2}) - \frac{1}{8\lambda}\sigma_{j}^n \\ & \spc+ \mathcal{O}(k_{\jph}-k_j),\\
    f(k_{\jph}, u_{j+1}^{n})-g(k_{j+1},u_{j+1}^{n})
    &= \frac{\lambda}{2}(a_{j+1}^{n})^2\sigma_{j+1}^n- \frac{1}{8}(\lambda a_{j+1}^{n}\sigma_{j+1}^n)^{2}f_{uu}(k_{j+1},\zeta_{3}) \\ & \spc - \frac{1}{8\lambda}\sigma_{j+1}^n + \mathcal{O}(k_{\jph}-k_\jpo),
\end{nalign}
Further, substituting \eqref{eq:fminusg_taylor} in \eqref{eq:Jbar_aa} yields
\begin{align}\label{eq:barJ_refor}
    \bar{J} 
    &= -\frac{1}{2}\lambda \Delta u_{\jph}^{n} \left[ \frac{\lambda}{2}(a_{j}^{n})^2\sigma_{j}^n+ \frac{\lambda}{2}(a_{j+1}^{n})^2\sigma_{j+1}^n - \frac{1}{8\lambda}\sigma_{j}^n- \frac{1}{8\lambda}\sigma_{j+1}^n \right] \\
    & \hspace{0.5cm} + \frac{1}{16}\lambda (\Delta u_{\jph}^{n})^3\left[\left(\lambda a_{j+1}^{n}\left(\frac{\sigma_{j+1}^{n}}{\Delta u_{\jph}^{n}}\right)\right)^2 + \left(\lambda a_{j}^{n}\left(\frac{\sigma_{j}^{n}}{\Delta u_{\jph}^{n}}\right)\right)^2 \right.\notag\\&\spc \spc  \left.+\frac{4}{3}f_{uu}(k_{\jph}, \zeta_{1}) \right] + \mathcal{O}(k_{\jph}-k_j) + \mathcal{O}(k_{\jph}-k_{j+1}), \notag\\
    & =  \frac{1}{16}\Delta u_{\jph}^{n} \mathcal{A}_{1} + \frac{1}{16}\lambda (\Delta u_{\jph}^{n})^3 \mathcal{A}_{2} + \mathcal{O}(k_{\jph}-k_j) + \mathcal{O}(k_{\jph}-k_{j+1}), \notag
\end{align}
where we define
\begin{align}\label{eq:A1A2defn}
    \mathcal{A}_{1} &:=  \sigma_{j}^n + \sigma_{j+1}^n - 4\lambda^2(a_{j}^{n})^2\sigma_{j}^n - 4\lambda^2(a_{j+1}^{n})^2\sigma_{j+1}^n,\\
    \mathcal{A}_{2} &:= \left(\lambda a_{j+1}^{n}\left(\frac{\sigma_{j+1}^{n}}{\Delta u_{\jph}^{n}}\right)\right)^2 + \left(\lambda a_{j}^{n}\left(\frac{\sigma_{j}^{n}}{\Delta u_{\jph}^{n}}\right)\right)^2 +\frac{4}{3}f_{uu}(k_{\jph}, \zeta_{1}). \notag
\end{align}
Now, we concentrate on the terms $\mathcal{A}_{1}$ and $\mathcal{A}_{2}.$ The hypothesis \ref{hyp:H2} along with \eqref{eq:kappa}  give us\begin{nalign}\label{eq:A2bound}
    \mathcal{A}_{2} \geq \frac{4}{3}\gamma_{1} \quad \mbox{and} \quad \abs{\mathcal{A}_{2}} \leq 2\kappa^{2}+\frac{4}{3}\gamma_{2}. 
\end{nalign}
Next, introducing the notations \begin{nalign}\label{eq:ujphdef}
\displaystyle\tilde{k}_\jph := \frac{k_{j}+k_{j+1}}{2}, \quad \displaystyle \tilde{u}^{n}_\jph := \frac{u^n_{j}+u^n_{j+1}}{2},
    \quad \mbox{and} \quad \tilde{a}_{\jph}^{n}:= f_{u}(\tilde{k}_{\jph},\tilde{u}_{\jph}^{n}),
\end{nalign} and using the Taylor expansion of $f_{u}$ about the point $(\tilde{k}_{\jph},\tilde{u}_{\jph}^{n})$, we obtain
\begin{nalign}\label{eq:a_taylor}
    a_{j}^{n} &= \tilde{a}_{\jph}^{n} - \frac{1}{2}(k_{j+1}-k_{j})f_{uk}(c_{4},\zeta_{4}) - \frac{1}{2}(u_{j+1}^{n}-u_{j}^{n})f_{uu}(c_{4},\zeta_{4}),\\
    a_{j+1}^{n} &= \tilde{a}_{\jph}^{n} + \frac{1}{2}(k_{j+1}-k_{j})f_{uk}(c_{5},\zeta_{5}) + \frac{1}{2}(u_{j+1}^{n}-u_{j}^{n})f_{uu}(c_{5},\zeta_{5}).
\end{nalign}
 where $c_{4} \in \mathcal{I}(k_{j}, \tilde{k}_{\jph}),$ $c_{5} \in \mathcal{I}(k_{j+1}, \tilde{k}_{\jph}),$ $\zeta_{4} \in  \mathcal{I}(u_{j}^{n}, \tilde{u}_{\jph}^{n})$ and $\zeta_{5} \in \mathcal{I}(u_{j+1}^{n}, \tilde{u}_{\jph}^{n} ).$ From these, we obtain 
\begin{nalign}\label{eq:aj_sq&ajph_sq}
    (a_{j}^{n})^{2} &= (\tilde{a}_\jph^{n})^2 + \frac{1}{4}(\Delta u_{\jph}^{n})^{2}(f_{uu}(c_4,\zeta_4))^2 - \tilde{a}_{\jph}\Delta u_{\jph}^{n} f_{uu}(c_4,\zeta_4) + \mathcal{O}(\Delta k_{\jph}
),\\
    (a_{j+1}^{n})^{2} &= (\tilde{a}_\jph^{n})^2 + \frac{1}{4}(\Delta u_{\jph}^{n})^{2}(f_{uu}(c_5,\zeta_5))^2 + \tilde{a}_{\jph}\Delta u_{\jph}^{n} f_{uu}(c_5,\zeta_5) + \mathcal{O}(\Delta k_{\jph}
),\\
\end{nalign}
Further, using \eqref{eq:aj_sq&ajph_sq} and the notation $\beta := \lambda \tilde{a}_{\jph}^{n},$ the term $\mathcal{A}_{1}$ in \eqref{eq:A1A2defn} can be reformulated  as
\begin{align}\label{eq:A1_refor}
    \mathcal{A}_{1} 
        &= \left(1-4\beta^{2}\right)(\sigma_{j}^n + \sigma_{j+1}^n) - \lambda (\Delta u_{\jph}^{n})^{2}\mathcal{A}_{3} + \mathcal{O}(\Delta k_{\jph}), 
\end{align} where 
\begin{nalign}\label{eq:A3}
    \mathcal{A}_{3} &:=  4\beta\frac{\sigma_{j+1}^{n}}{\Delta u_{\jph}^{n}}f_{uu}(c_5, \zeta_5)- 4\beta\frac{\sigma_{j}^{n}}{\Delta u_{\jph}^{n}}f_{uu}(c_4, \zeta_4) \\& \spc + \lambda (f_{uu}(c_4,\zeta_{4}))^2\sigma_{j}^n + \lambda (f_{uu}(c_5,\zeta_{5}))^2\sigma_{j+1}^n.
\end{nalign}
Replacing the term $
\mathcal{A}_{1}$ in \eqref{eq:barJ_refor} with its expression  \eqref{eq:A1_refor}, $\bar{J}$ can be rewritten as 
\begin{nalign}\label{eq:Jbar}
    \bar{J} & =  \frac{1}{16}(\Delta u_{\jph}^{n})^2 \left(1-4\beta^{2}\right)\frac{(\sigma_{j}^n + \sigma_{j+1}^n)}{\Delta u_{\jph}^{n}} + \frac{1}{16}\lambda(\Delta u_{\jph}^{n})^3(\mathcal{A}_2-\mathcal{A}_3)\\ &\spc + \mathcal{O}(k_{\jph}-k_j) + \mathcal{O}(k_{\jph}-k_{j+1}) + \mathcal{O}(\Delta k_{\jph}),
\end{nalign}
\par
 Next, we focus on the term $I$ in \eqref{eq:I_b}. The choice $S(k,u)=\frac{u^2}{2}$  simplifies $I$ as  
\begin{nalign}\label{eq:I}
    I 
    & = -\frac{1}{8}(\Delta u_{\jph}^{n})^2\left(1- \left(2\lambda\frac{\Delta g_{\jph}^{n}}{\Delta u_{\jph}^{n}}\right)^{2}\right).
\end{nalign}
Here, using \eqref{eq:gtaylor_j}, the term $\displaystyle\frac{\Delta g_{\jph}^{n}}{\Delta u_{\jph}^{n}}$ in \eqref{eq:I} can be expressed as follows 
\begin{nalign}\label{eq:dgbydu}
    \frac{\Delta g_{\jph}^{n}}{\Delta u_{\jph}^{n}}
    & = L + M + N,
\end{nalign}
where 
\begin{nalign}\label{eq:L}
    L &:=\frac{f(k_{j+1},u_{j+1}^{n})- f(k_{j},u_{j}^{n})}{\Delta u_{\jph}^{n}}, \\  M&:=- \frac{\left(\sigma_{j+1}^n\left(\frac{\lambda}{2}(a_{j+1}^{n})^2-\frac{1}{8\lambda}\right)- \sigma_{j}^n\left(\frac{\lambda}{2}(a_{j}^{n})^2-\frac{1}{8\lambda}\right)\right)}{\Delta u_{\jph}^{n}}, \,\,\, \\
    N&:= \frac{1}{8}\frac{\left((\lambda a_{j+1}^{n}\sigma^n_{j+1})^2f_{uu}(k_{j+1},\zeta_3)- (\lambda a_{j}\sigma^n_{j})^2f_{uu}(k_{j},\zeta_2)\right)}{\Delta u_{\jph}^{n}}.
\end{nalign}
Next, in order to simplify $L,$ we consider the following Taylor series expansions
\begin{nalign}\label{eq:ftaylor}
f(k_j, u^{n}_{j})&= f(\tilde{k}_\jph, \tilde{u}^{n}_{\jph}) - \frac{1}{2}\Delta k_{\jph}f_{k}(\tilde{k}_{\jph}, \tilde{u}_{\jph}^{n}) -\frac{1}{2}\Delta u_{\jph}^{n}f_{u}(\tilde{k}_{\jph}, \tilde{u}_{\jph}^{n})\\
&\hspace{0.5cm}+ \frac{1}{8}(\Delta u_{\jph}^{n})^2 f_{uu}(c_{6},\zeta_{6}) + \mathcal{O}(\Delta k_{\jph}),\\
f(k_{j+1}, u^{n}_{j+1})&= f(\tilde{k}_\jph, \tilde{u}^{n}_{\jph}) + \frac{1}{2}\Delta k_{\jph}f_{k}(\tilde{k}_{\jph}, \tilde{u}_{\jph}^{n}) +\frac{1}{2}\Delta u_{\jph}^{n}f_{u}(\tilde{k}_{\jph}, \tilde{u}_{\jph}^{n})\\
&\hspace{0.5cm}+ \frac{1}{8}(\Delta u_{\jph}^{n})^2 f_{uu}(c_{7},\zeta_{7}) + \mathcal{O}(\Delta k_{\jph}),
\end{nalign}
where  $c_{6} \in \mathcal{I}(k_{j}, \tilde{k}_{\jph}),$ $c_{7} \in \mathcal{I}(k_{j+1}, \tilde{k}_{\jph}),$ $\zeta_{6} \in  \mathcal{I}(u_{j}^{n}, \tilde{u}_{\jph}^{n}),$ $\zeta_{7} \in \mathcal{I}(u_{j+1}^{n}, \tilde{u}_{\jph}^{n} ).$ 
Now, in view of\eqref{eq:ftaylor}, the term $L$ in \eqref{eq:dgbydu} reads as 
\begin{nalign}\label{eq:Lmodified}
L 
&= \tilde{a}_{\jph}^{n}+ \frac{1}{8}\Delta u_{\jph}^{n} \left( f_{uu}(c_{7},\zeta_{7})-  f_{uu}(c_{6},\zeta_{6})\right) + \frac{\mathcal{O}(\Delta k_\jph)}{\Delta u_{\jph}^{n}}.
\end{nalign}
Further, the expressions in \eqref{eq:aj_sq&ajph_sq} allow us to write the term $M$ as
\begin{align}\label{eq:M}
    M 
    &= \frac{1}{\lambda}\left(\frac{1}{8}-\frac{1}{2}\beta^2\right)\frac{\Delta \sigma_{\jph}^n}{\Delta u_{\jph}^{n}} -\frac{\beta}{2}\left(\sigma_{j+1}^{n}f_{uu}(c_{5}, \zeta_{5})+ \sigma_{j}^{n}f_{uu}(c_{4}, \zeta_{4})\right) \\
    &\hspace{0.5cm} -\frac{1}{8}\left(\lambda\left(\frac{\sigma_{j+1}^{n}}{\Delta u_{\jph}^{n}}\right)(\Delta u_{\jph}^{n}f_{uu}(c_5, \zeta_5))^2 - \lambda\left(\frac{\sigma_{j}^{n}}{\Delta u_{\jph}^{n}}\right)(\Delta u_{\jph}^{n}f_{uu}(c_4, \zeta_4))^2\right)\notag \\ & \spc + \frac{\mathcal{O}(\Delta k_{\jph})}{\Delta u_{\jph}^{n}}. \notag
\end{align}
Combining the expressions \eqref{eq:Lmodified} and  \eqref{eq:M}, we may write \eqref{eq:dgbydu} as
\begin{nalign}\label{eq:gdifbyudif_fin}
    \lambda \frac{\Delta g_{\jph}^{n}}{\Delta u_{\jph}^{n}}
    & = \beta + \frac{1}{8}\left(1-4\beta^2\right)\frac{\Delta \sigma_{\jph}^n}{\Delta u_{\jph}^{n}} + \mathcal{A}_4 +  \frac{\mathcal{O}(\Delta k_{\jph})}{\Delta u_{\jph}^{n}},
\end{nalign}
where we define
\begin{nalign}\label{eq:A_4}
    \mathcal{A}_{4} &:= \frac{\lambda}{8}\Delta u_{\jph}^{n} \left( f_{uu}(c_{7},\zeta_{7})-  f_{uu}(c_{6},\zeta_{6})\right)  -\lambda\frac{\beta}{2}\left(\sigma_{j+1}^{n}f_{uu}(c_{5}, \zeta_{5})+ \sigma_{j}^{n}f_{uu}(c_{4}, \zeta_{4})\right)\\
    &\hspace{0.5cm} -\frac{1}{8}\left[\left(\frac{\sigma_{j+1}^{n}}{\Delta u_{\jph}^{n}}\right)(\lambda\Delta u_{\jph}^{n}f_{uu}(c_5, \zeta_5))^2 - \left(\frac{\sigma_{j}^{n}}{\Delta u_{\jph}^{n}}\right)(\lambda\Delta u_{\jph}^{n}f_{uu}(c_4, \zeta_4))^2\right] \\ & \spc + \lambda N.
\end{nalign}
In view of \eqref{eq:gdifbyudif_fin}, the expression \eqref{eq:I} now reads as
\begin{align}\label{eq:I_a}
    I 
    & = -\frac{1}{8} (\Delta u_{\jph}^{n})^2 (1-4\beta^2)\left[1-\frac{1}{16}(1-4\beta^2)\left(\frac{\Delta \sigma_{\jph}^{n}}{\Delta u_{\jph}^{n}}\right)^2 - \beta\left(\frac{\Delta \sigma_{\jph}^{n}}{\Delta u_{\jph}^{n}}\right)\right] \\
    & \hspace{0.5cm}  -\frac{1}{8} (\Delta u_{\jph}^{n})^2 \left[- \mathcal{A}_4\left(8\beta+(1-4\beta^2)\left(\frac{\Delta \sigma_{\jph}^{n}}{\Delta u_{\jph}^{n}}\right)\right)-4(\mathcal{A}_4)^2\right] + \mathcal{O}(\Delta k_{\jph}).\notag
\end{align}
Further, combining \eqref{eq:Jbar} and \eqref{eq:I_a}, we write 
\begin{nalign}\label{eq:IplusJbar}
    I+\bar{J} 
    & = -\frac{1}{8}\eta_{\jph}^{n} (\Delta u_{\jph}^{n})^2 + \frac{1}{16} \theta_{\jph}^{n}\lambda(\Delta u_{\jph}^{n})^3  + \mathcal{O}(\Delta k_{\jph}) \\ &\spc+ \mathcal{O}(k_{\jph}-k_j) + \mathcal{O}(k_{\jph}-k_{j+1}),
\end{nalign}
where 
\begin{align}
    \eta_{\jph}^{n} &:=(1-4\beta^2)\left[1-\frac{1}{16}(1-4\beta^2)\left(\frac{\Delta \sigma_{\jph}^{n}}{\Delta u_{\jph}^{n}}\right)^2\right. \label{eq:eta}\\ & \spc\spc  \left.- (\beta+\mathcal{A}_{4})\left(\frac{\Delta \sigma_{\jph}^n}{\Delta u_{\jph}^{n}}\right)-\frac{(\sigma_{j}^n + \sigma_{j+1}^n)}{2\Delta u_{\jph}^n}\right], \notag \\
    \theta_{\jph}^{n} &:= \mathcal{A}_2-\mathcal{A}_3+ 16 \frac{\beta\mathcal{A}_4}{\lambda\Delta u_{\jph}^n}+8\frac{(\mathcal{A}_4)^2}{\lambda \Delta u_\jph^n}. \label{eq:theta}
\end{align}
Now, our immediate aim is to reformulate the expression \eqref{eq:IplusJbar} in the form 
\begin{align}\label{eq:IplusJbar_mod}
    I+\bar{J} 
    & = -\frac{1}{8}\tilde{\eta}_{\jph}^{n} (\Delta u_{\jph}^{n})^2 + \frac{1}{16} \tilde{\theta}_{\jph}^{n}\lambda(\Delta u_{\jph}^{n})^3 + \mathcal{O}(\Delta k_{\jph})\\ & \spc + \mathcal{O}(k_{\jph}-k_j) + \mathcal{O}(k_{\jph}-k_{j+1}),
\end{align}
with the coefficients $\tilde{\eta}_{\jph}^{n}, \tilde{\theta}_{\jph}^{n}\ge 0.$
To achieve this, we begin with  decomposing the term $\mathcal{A}_3$ as 
\begin{nalign}\label{eq:A3-A3bar}
    \mathcal{A}_3 =\tilde{\mathcal{A}}_3 + \mathcal{O}(\Delta k_{\jph}),
\end{nalign}
where we define
\begin{align*}
    \tilde{\mathcal{A}}_{3} &:=     \left[4\beta\frac{\sigma_{j+1}^{n}}{\Delta u_{\jph}^{n}}f_{uu}(c_5, \zeta_5)- 4\beta\frac{\sigma_{j}^{n}}{\Delta u_{\jph}^{n}}f_{uu}(c_4, \zeta_4)\right]\\
    & \hspace{0.5cm} + \left[\lambda (f_{uu}(c_4,\zeta_{4}))(2\tilde{a}_{\jph}^{n}-2a_{j}^{n})\left(\frac{\sigma_{j}^n}{\Delta u^{n}_{\jph}}\right)\right. \\&\spc \spc   \left.+ \lambda f_{uu}(c_5,\zeta_{5})(2a_{\jpo}^{n}-2\tilde{a}_{\jph}^{n})\left(\frac{\sigma_{j+1}^n}{\Delta u^{n}_\jph}\right) \right].
\end{align*}
Here,  the term $\tilde{\mathcal{A}}_{3}$ can be bounded as 
\begin{nalign}\label{eq:A3barbound}
    \abs{\tilde{\mathcal{A}}_{3}} & \leq 8\gamma_{2}(\abs{\beta}+\kappa). 
\end{nalign}
This is derived using \eqref{eq:kappa} and the hypothesis \ref{hyp:H2}, along with the fact that $0\leq \frac{\sigma_{j+1}^n}{\Delta u^{n}_\jph} \leq 1$ (see \eqref{eq:slopes}).
Using the same arguments, we obtain a bound on the term $\mathcal{A}_{4}$ in \eqref{eq:A_4} as
\begin{align}\label{eq: A4bound1}
    \abs{\mathcal{A}_{4}} \leq C_{0} \lambda \abs{\Delta u_\jph^{n}} \gamma_2,
\end{align}
where $\displaystyle C_{0} := \frac{1}{4}+\kappa+ \frac{\kappa^2}{4}+ \frac{1}{2}\lambda C_{u_0}\gamma_2.$ In addition to this bound, we also need to split $\mathcal{A}_4$ into a suitable form. To do this, we rewrite the term $\lambda N$ in \eqref{eq:A_4} by using the expressions for $f_{uu}(k_j,\zeta_{2})$ and $f_{uu}(k_{\jpo},\zeta_{3})$ taken from the Taylor expansions given in \eqref{eq:gtaylor_j},
\begin{nalign}\label{eq:lambdaN}
    \lambda N = E_{1}+ E_{2} + E_{3},
\end{nalign}
where 
\begin{nalign}\label{eq:E1E2E3}
    E_{1} &:= \lambda \frac{\left(f(k_{j+1}, u_{j+1}^{n+\frac{1}{2}}) - f(k_{j}, u_{j}^{n+\frac{1}{2}})\right)}{\Delta u_{\jph}^{n}}, \\ E_2 &:= - \lambda \frac{\left(f(k_{j+1}, u_{j+1}^{n}) - f(k_{j}, u_{j}^{n})\right)}{\Delta u_{\jph}^{n}}, \,\, 
 E_3 := \lambda^2 \frac{\left(a_{j+1}^2(\sigma_{j+1}^n)^{2}-a_{j}^2(\sigma_{j}^n)^{2}\right)}{2\Delta u_{\jph}^{n}}.
\end{nalign}
 Adding and subtracting $\lambda f(k_{j}, u_{j+1}^{n+\frac{1}{2}})$ to the numerator and  applying the mean value theorem, we may write  $E_1$ as 
\begin{align}\label{eq:E1}    
E_1 
    & = \frac{\mathcal{O}(\Delta k_{\jph})}{\Delta u_{\jph}^{n}} + \lambda f_{u}(k_j, \zeta_8)\left(1 - \frac{\lambda}{2}\left(a_{j+1}^{n}\frac{\sigma_{j+1}^{n}}{\Delta u_{\jph}^{n}}- a_{j}^{n}\frac{\sigma_{j}^{n}}{\Delta u_{\jph}^{n}}\right)\right), 
\end{align}
 for some $c_{8} \in \mathcal{I}(k_j, k_{j+1})$ and $\zeta_{8} \in \mathcal{I}(u_{j}^{\nph}, u_{j+1}^{\nph}),$ where we have used the definition of $u_{\jph}^{\nph}$ from \eqref{eq:mid-time_values}.  Similarly, adding and subtracting $\lambda(f(k_{j}, u_{j+1}^{n})),$ we may write 
$E_{2}$ as
\begin{nalign}\label{eq:E2}
    E_{2} &= \lambda\frac{\mathcal{O}(\Delta k_{\jph})}{\Delta u_{\jph}^{n}} - \lambda f_{u}(k_j, \zeta_9),
\end{nalign} for some $c_{9} \in \mathcal{I}(k_j, k_{j+1})$ and $\zeta_{9} \in \mathcal{I}(u_{j}^{n}, u_{j+1}^{n}).$ 
\par
Now, substituting the expressions for $E_{1}, E_2$ and $E_3$ from \eqref{eq:E1}, \eqref{eq:E2} and  \eqref{eq:E1E2E3}, respectively, and using the relations \eqref{eq:a_taylor} and \eqref{eq:ftaylor} for the terms $f_{uu}(c_{4},\zeta_{4}),$ $f_{uu}(c_{5},\zeta_{5}),$ $f_{uu}(c_{6},\zeta_{6})$ and $f_{uu}(c_{7},\zeta_{7}),$ 
the term $\mathcal{A}_4$  can be reformulated as

\begin{align}\label{eq:A4-A4bar}
    \mathcal{A}_4 &= \tilde{\mathcal{A}}_{4} + \frac{\mathcal{O}(\Delta k_{\jph})
    }{\Delta u_{\jph}^{n}},
\end{align}
where $\tilde{\mathcal{A}}_{4}:= \tilde{\mathcal{A}}_{4}^{1} + \tilde{\mathcal{A}}_{4}^{2}+ \tilde{\mathcal{A}}_{4}^{3}+\tilde{\mathcal{A}}_{4}^{4}+\tilde{\mathcal{A}}_{4}^{5} + \tilde{\mathcal{A}}_{4}^{6},$ and
\begin{align*}
    \tilde{\mathcal{A}}_{4}^{1} &:= \frac{1}{2}\lambda \left[f_{u}(k_{j+1}, \zeta_{11})+ f_{u}(k_{j}, \zeta_{10})-2f_{u}(\tilde{k}_{\jph}, u_{\jph}^{n}))\right], \\
    \tilde{\mathcal{A}}_{4}^{2} &:=-\lambda \beta\left[\frac{\sigma^n_{j+1}}{\Delta u_\jph^{n}}(a_{\jpo}^{n}-\tilde{a}_{\jph}^{n})+ \frac{\sigma^n_{j}}{\Delta u_{\jph}^{n}}(\tilde{a}_{\jph}^{n}-a_{j}^{n})\right],\\
    \tilde{\mathcal{A}}_{4}^{3} &:= -\frac{1}{8}\lambda^{2}\left[\frac{\sigma_{j+1}^n}{\Delta u_\jph^{n}}\left(4(a_{j+1}^{n})^2+ 4(\tilde{a}_{\jph}^{n})^2-8a_{j+1}^{n}\tilde{a}_{\jph}^{n}\right)\right]\\
    &\hspace{0.5cm} + \frac{1}{8}\lambda^2\left[\frac{\sigma_{j}^n}{\Delta u_\jph^{n}}\left(4(a_{j}^{n})^2+ 4(\tilde{a}_{\jph}^{n})^2-8a_{j}^{n}\tilde{a}_{\jph}^{n}\right)\right],\\
    \tilde{\mathcal{A}}_{4}^{4} &:= \lambda f_{u}(k_j, \zeta_8)\left[1 - \frac{\lambda}{2}\left(a_{j+1}^{n}\frac{\sigma_{j+1}^{n}}{\Delta u_{\jph}^{n}}- a_{j}^{n}\frac{\sigma_{j}^{n}}{\Delta u_{\jph}^{n}}\right)\right], \\
    \tilde{\mathcal{A}}_{4}^{5}&:=- \lambda f_{u}(k_j, \zeta_9), \quad   
    \tilde{\mathcal{A}}_{4}^{6} := \lambda^2 \frac{\left(a_{j+1}^2(\sigma_{j+1}^n)^{2}-a_{j}^2(\sigma_{j}^n)^{2}\right)}{2\Delta u_{\jph}^{n}}.
\end{align*}
Using \eqref{eq:kappa}, \eqref{Linf_bd} and the fact that $0\leq \frac{\sigma_{j+1}^n}{\Delta u^{n}_\jph} \leq 1,$ we  easily obtain the bounds 
\begin{align*}
    \abs{\tilde{\mathcal{A}}_{4}^{1}} &\leq 2\kappa, \quad \abs{\tilde{\mathcal{A}}_{4}^{2}} \leq 4\kappa^{2},  \quad \abs{\tilde{\mathcal{A}}_{4}^{3}} \leq 4\kappa^2, \quad \abs{\tilde{\mathcal{A}}_{4}^{4}} \leq \kappa(1+\kappa), \\ \abs{\tilde{\mathcal{A}}_{4}^{5}} &\leq \kappa \quad \mbox{and} \quad \abs{\tilde{\mathcal{A}}_{4}^{6}} \leq 2 C_{u_0}\kappa^2,
\end{align*}
and hence
\begin{nalign}\label{eq:A4_bound2}
\abs{\tilde{\mathcal{A}}_4} \leq C_{1} \kappa,
\end{nalign}
where $C_{1} := 4+(9+2C_{u_{0}})\kappa.$ Upon rewriting \eqref{eq:theta} using the expressions  \eqref{eq:A3-A3bar} and \eqref{eq:A4-A4bar} for $\mathcal{A}_3$ and $\mathcal{A}_4,$ respectively, we arrive at \begin{nalign}
    \theta_{\jph}^{n} &= \tilde{\theta}_{\jph}^{n} + \frac{\mathcal{O}(\Delta k_{\jph})}{(\Delta u_{\jph}^{n})^2},
\end{nalign}
where we define
\begin{nalign}\label{eq:theta_modified}
    \tilde{\theta}_{\jph}^{n}&:=\mathcal{A}_2-\tilde{\mathcal{A}}_3 + 16 \frac{\beta\mathcal{A}_4}{\lambda\Delta u_{\jph}^{n}}+ 8\frac{\mathcal{A}_4\tilde{\mathcal{A}}_4}{\lambda \Delta u_\jph^{n}}.
\end{nalign}
Analogously, in view of \eqref{eq:A4-A4bar}, \eqref{eq:eta} yields
\begin{align*}
    \eta_{\jph}^{n} 
    &= \tilde{\eta}_{\jph}^{n} + \frac{\mathcal{O}(\Delta k_{\jph})}{\Delta u_{\jph}^{n}},
\end{align*} 
where we define
\begin{nalign}\label{eq:eta_modified}
    \tilde{\eta}_{\jph}^{n}&:=(1-4\beta^2)\left[1-\frac{1}{16}(1-4\beta^2)\left(\frac{\Delta \sigma_{\jph}^{n}}{\Delta u_{\jph}^{n}}\right)^2 \right.\\ & \spc \spc  \left.- (\beta+\tilde{\mathcal{A}}_{4})\left(\frac{\Delta \sigma_{\jph}^{n}}{\Delta u_{\jph}^{n}}\right)-\frac{(\sigma_{j}^n + \sigma_{j+1}^n)}{2\Delta u_{\jph}^{n}}\right],
\end{nalign}
Hence, we obtain the desired reformulation \eqref{eq:IplusJbar_mod} of  \eqref{eq:IplusJbar}.
\par
Next, we show that $\tilde{\eta}_{\jph}^{n} \geq 0$ and $\tilde{\theta}_{\jph}^{n} > 0, \, \mbox{for all}\, \jinz.$ 
{First, we consider the term $\tilde{\eta}_{\jph}^{n},$}
 and note that 
$1-4\beta^2 \geq 0,$ by the CFL condition \eqref{eq:cfl_cubicest}. Now, using the bound \eqref{eq:A4_bound2}, a portion of the term $\tilde{\eta}^{n}_{\jph}$ can be bounded as
\begin{align}\label{eq:etageqzero}
    &1-\frac{1}{16}(1-4\beta^2)\left(\frac{\Delta \sigma_{\jph}^{n}}{\Delta u_{\jph}^{n}}\right)^2 - (\beta+\tilde{\mathcal{A}}_{4})\left(\frac{\Delta \sigma_{\jph}^{n}}{\Delta u_{\jph}^{n}}\right)-\frac{(\sigma_{j}^n + \sigma_{j+1}^n)}{2\Delta u_{\jph}^{n}}\\
    &\hspace{0.5cm} \geq \left(\frac{1}{2}-\frac{1}{16}(1-4\beta^2)\abs{\frac{\Delta \sigma_{\jph}^{n}}{\Delta u_{\jph}^{n}}}-\abs{\beta+\tilde{\mathcal{A}}_{4}}\right)\abs[\Big]{\frac{\Delta \sigma_{\jph}^n}{\Delta u_{\jph}^{n}}} +1 \notag\\ & \spc \spc -\left( \abs[\Big]{\frac{\Delta \sigma_{\jph}^{n}}{2\Delta u_{\jph}^{n}}} +\frac{(\sigma_{j}^n + \sigma_{j+1}^n)}{2\Delta u_{\jph}^{n}}\right) \notag\\
    &\hspace{0.5cm}\geq \left(\frac{1}{2}-\frac{1}{16}(1-4\beta^2)-\abs{\beta+C_{1}\kappa}\right) \abs[\Big]{\frac{\Delta \sigma_{\jph}^{n}}{\Delta u_{\jph}^{n}}}\geq 0.\notag
\end{align} Here, we use the bound $\displaystyle (1+C_{1})\kappa \leq \frac{7}{16}, $ that follows from the CFL condition \eqref{eq:cfl_cubicest}. Thus, we obtain $\tilde{\eta}_\jph^{n} \geq 0,$ for all $\jinz.$\par
Now, collecting the estimates \eqref{eq:A2bound}, \eqref{eq:A3barbound}, \eqref{eq: A4bound1}  and \eqref{eq:A4_bound2}, we deduce that 
\begin{nalign}\label{eq:thetageq0}
    \tilde{\theta}_{\jph}^{n}
    & \geq \frac{4}{3}\gamma_{1}  -   8\gamma_{2}(\abs{\beta}+\kappa) - 16 \abs{\beta} \gamma_2 C_{0} - 8\gamma_2 \kappa C_{0} C_{1} \geq  \gamma_1 > 0, \, \, \mbox{for all} \,\, \jinz, 
\end{nalign} as the CFL condition \eqref{eq:cfl_cubicest} implies that $\displaystyle \kappa \leq \frac{\gamma_{1}}{24\gamma_{2}\left(2(1+C_{0})+C_{0}C_{1}\right)}$ for $\Delta t \leq \Delta x.$ We also obtain an upper bound on  $\tilde{\theta}_{\jph}^{n},$ by invoking \eqref{eq:A2bound}, \eqref{eq:A3barbound}, \eqref{eq: A4bound1} and \eqref{eq:A4_bound2}, as follows 
\begin{align}\label{eq:thetabar_bd}
    \tilde{\theta}_{\jph}^{n}\leq C_{2}, \,\, \mbox{where} \,\, C_{2}:= 2\kappa^{2} +\frac{4}{3}\gamma_{2} + 8\kappa\gamma_{2}(2+2C_{0}+C_{0}C_{1}).
\end{align}
Finally, in light of the equation \eqref{eq:IplusJbar_mod}, the expression \eqref{eq:numentropyprod} reduces to
\begin{nalign}\label{eq:entropprod-final}  &S(k_{\jph},u^{n+1}_{\jph}) - \frac{1}{2}\left( S(k_{j},u^n_{j})+S(k_{j+1} ,u^n_{j+1})\right) + \lambda\left(G(k_{j+1}, u^n_{j+1})- G(k_{j},u^n_{j}))\right) \\
    &= -\frac{1}{8}\tilde{\eta}_{\jph}^{n} (\Delta u_{\jph}^{n})^2 + \frac{1}{16} \tilde{\theta}_{\jph}^{n}\lambda(\Delta u_{\jph}^{n})^3 + \Phi_{\jph} ,
\end{nalign} 
where $\Phi_{\jph} = \mathcal{O}(\Delta k_{\jph})+ \mathcal{O}(k_{\jph}-k_j) + \mathcal{O}(k_{j+1}-k_{\jph}).$
Summing \eqref{eq:entropprod-final} over $n= 0,1, \dots, N-1$ and $j+\frac{n}{2} \in \mathbb{Z}$ with $\abs{j} \leq J,$  and with the choice  $\displaystyle S(k,u)= \frac{u^{2}}{2}$ yields 
\begin{align}\label{eq:sum_entropy_prod}
&\frac{1}{2}\sum_{\abs{j}\leq J} (u^N_{j})^2 - \frac{1}{2}\sum_{\abs{j}\leq J} (u^0_{j})^2 + \lambda \sum_{n=0}^{N-1} \left(G(k_{J},u^n_{J})- G(k_{-J},u^n_{-J})\right)\\
&=  -\frac{1}{8}\sum_{n=0}^{N-1}\sum_{\substack{\abs{j} \leq J \\ j+\frac{n}{2} \in \mathbb{Z}}}\tilde{\eta}_{\jph}^{n} (\Delta u_{\jph}^{n})^2 + \frac{1}{16}\sum_{n=0}^{N-1}\sum_{\substack{\abs{j} \leq J \\ j+\frac{n}{2} \in \mathbb{Z}}}\tilde{\theta}_{\jph}^{n}\lambda(\Delta u_{\jph}^{n})^3_{-} \notag\\& \hspace{0.5cm} +\frac{1}{16}\sum_{n=0}^{N-1}\sum_{\substack{\abs{j} \leq J \\ j+\frac{n}{2} \in \mathbb{Z}}}\tilde{\theta}_{\jph}^{n}\lambda(\Delta u_{\jph}^{n})^3_{+} + \sum_{n=0}^{N-1}\sum_{\substack{\abs{j} \leq J \\ j+\frac{n}{2} \in \mathbb{Z}}} \Phi_{\jph}. \notag
\end{align}
Rearranging \eqref{eq:sum_entropy_prod} and taking into account $\tilde{\eta}_{\jph}^{n} \geq 0$ (see \eqref{eq:etageqzero}), it follows that 
\begin{align}\label{eq:cubicest_prefinal}
     &-\frac{1}{16}\sum_{n=0}^{N-1}\sum_{\substack{\abs{j} \leq J \\ j+\frac{n}{2} \in \mathbb{Z}}}\tilde{\theta}_{\jph}^{n}\lambda(\Delta u_{\jph}^n)^3_{-}\\ &\spc \leq \frac{1}{2}\sum_{\abs{j}\leq J} (u^0_{j})^{2} -\frac{1}{2} \sum_{\abs{j}\leq J} (u^N_{j})^{2} - \lambda \sum_{n=0}^{N-1} \left(G(k_{J},u^n_{J})- G(k_{-J},u^n_{-J})\right)\notag\\
     &\spc \hspace{0.5cm} +\frac{1}{16}\sum_{n=0}^{N-1}\sum_{\substack{\abs{j} \leq J \\ j+\frac{n}{2} \in \mathbb{Z}}}\tilde{\theta}_{\jph}^{n}\lambda(\Delta u_{\jph}^n)^3_{+}  + C_{3}N \norm{k}_{BV} . \notag
\end{align} for some constant $C_{3} \geq 0.$
\par
Now,  summing \eqref{eq:oslsfinal}  of Lemma \ref{lemma:osle} over $n= 0, \dots, N-1,$ we see that
\begin{align}\label{eq:oslssummed}
        \frac{1}{500} \lambda \gamma_{1}\sum_{n=0}^{N-1}\sum_{\substack{\abs{j} \leq J \\ j+\frac{n}{2} \in \mathbb{Z}}} (\Delta u_{\jmh}^{n})_{+}^{3} 
    &\leq  \frac{1}{500} \lambda \gamma_{1}\sum_{n=0}^{N-1}\sum_{j+\frac{n}{2} \in \mathbb{Z}} (\Delta u_{\jmh}^{n})_{+}^{3} 
   \\&\leq \sum_{n=0}^{N-1}\sum_{j+\frac{n}{2} \in \mathbb{Z}}\left((\Delta u^{n}_{\jmh})_{+}^2-\left(\Delta u^{n+1}_{j}\right)_{+}^{2} \right) + \sum_{n=0}^{N-1}\Psi\norm{k}_{BV} \notag\\
    & \leq \sum_{\jinz}(\Delta u^{0}_{\jmh})_{+}^2- \sum_{j \in \mathbb{Z}}\left(\Delta u^{N}_{j}\right)_{+}^{2} + \sum_{n=0}^{N-1}\Psi\norm{k}_{BV}\notag \\
    &\leq \sum_{j \in \mathbb{Z}}(\Delta u^{0}_{\jmh})_{+}^2 + \sum_{n=0}^{N-1}\Psi\norm{k}_{BV}, \notag\\
    & \leq 2\norm{u_0}\norm{u_{0}}_{BV}+ \Psi\norm{k}_{BV}N. \notag
    \end{align}
The estimate \eqref{eq:oslssummed} derived above, together with the bound $\gamma_1 <\tilde{\theta}_{\jph}^{n}\leq C_2$ given in  \eqref{eq:thetageq0} and \eqref{eq:thetabar_bd}, applied to \eqref{eq:cubicest_prefinal}, yields
\begin{align}\label{eq:negcubessum}
     &-\frac{1}{16}\gamma_{1}\lambda\sum_{n=0}^{N-1}\sum_{\substack{\abs{j} \leq J \\ j+\frac{n}{2} \in \mathbb{Z}}}(\Delta u_{\jph}^{n})^3_{-} \\&\spc \leq \frac{1}{2}\sum_{\abs{j}\leq J} (u^0_{j})^{2} -\frac{1}{2} \sum_{\abs{j}\leq J} (u^N_{j})^{2} - \lambda \sum_{n=0}^{N-1} \left(G(k_{J},u^n_{J})- G(k_{-J},u^n_{-J})\right)\notag\\
     &\spc \hspace{0.5cm}+\frac{125}{4}\frac{C_2}{\gamma_{1}}\left(2\norm{u_0}\norm{u_{0}}_{BV}\right)+\left(\frac{125}{4}\frac{C_2}{\gamma_{1}} \Psi N +  C_{3}N\right)\norm{k}_{BV}\notag\\
     & \spc \leq \frac{1}{2}\norm{u_0}^{2}(2J+1) + 2\lambda \norm{G}N+\frac{125}{4}\frac{C_2}{\gamma_{1}}\left(2\norm{u_0}\norm{u_{0}}_{BV}\right)\notag\\ & \spc \spc +\left(\frac{125}{4}\frac{C_2}{\gamma_{1}} \Psi N +  C_{3}N\right)\norm{k}_{BV}. \notag
\end{align}
Observe that 
\begin{align}\label{eq:cubicest_pref}
    \Delta x\sum_{n=0}^{N-1}\sum_{\substack{\abs{j} \leq J \\ j+\frac{n}{2} \in \mathbb{Z}}}\abs{\Delta u_{\jph}^{n}}^3 = \Delta x \sum_{n=0}^{N-1}\sum_{\substack{\abs{j} \leq J \\ j+\frac{n}{2} \in \mathbb{Z}}} \left((\Delta u_{\jmh}^{n})_{+}^{3}- (\Delta u_{\jmh}^{n})_{-}^{3}  \right).
\end{align}
Now, using the estimates \eqref{eq:oslssummed} and \eqref{eq:negcubessum} in the above equation, we obtain the desired estimate \eqref{eq:C_cubiclemma}
 with \begin{align*}
    C(X,T) &:= \frac{C_4}{\lambda\gamma_1} 1000 \left(1+\frac{C_2}{\gamma_1}\right) \norm{u_0}\norm{u_{0}}_{BV}+
      \frac{32}{\lambda\gamma_1}\norm{u_0}^2 X  \\
     &\spc \hspace{0.5cm} + \frac{32}{\lambda\gamma_1} \norm{G}T +\frac{1}{\lambda^2\gamma_1}T \left(500\Psi+500\frac{C_2}{\gamma_{1}} \Psi  +  16C_{3}\right)\norm{k}_{BV},
\end{align*} when $\Delta x \leq C_4,$  for a constant $C_4 > 0.$ This concludes the proof.
\end{proof}
\end{lemma}
Now, we derive a quadratic estimate on the spatial differences of the approximate solutions.
\begin{lemma}(Quadratic estimate)\label{lemma:quadest}
Let the initial datum $u_{0} \in (\mathrm{L}^{\infty}\cap \mathrm{BV})(\mathbb{R})$ be such that  $\norm{u_{0}} \leq C_{u_{0}}.$  Consider the cell-average approximate solutions $\{u^{n}_{j}\}$ generated by the scheme \eqref{eq:NTscheme} and define
\begin{align}\label{eq:nujph}
    \nu_{\jph}^n&:= \frac{1}{8}  (1-4(\beta^{n}_{\jph})^2)\left[1-\frac{1}{16}(1-4(\beta^{n}_{\jph})^2)\left(\frac{\Delta \sigma_{\jph}^n}{\Delta u_{\jph}^n}\right)^2 \right.\\& \left. \spc - \beta^{n}_{\jph}\left(\frac{\Delta \sigma_{\jph}^n}{\Delta u_{\jph}^n}\right) - \frac{\sigma_{j+1}^{n}+  \sigma_{j}^{n}}{2\Delta u_{\jph}^{n}}\right]f_{uu}(k_{\jph}, \tilde{u}_{\jph}^{n}), \notag
\end{align}
with $\beta_{\jph}^{n}:= \lambda \tilde{a}_{\jph}^n,$ where $ \tilde{a}_{\jph}^n$ and $\tilde{u}_{\jph}^{n}$ are as in \eqref{eq:ujphdef}. If the CFL condition \eqref{eq:cfl_cubicest} holds, then we have:\\
(i) $\nu_{\jph}^n \geq 0$ for $\jinz,$ and \\ (ii) For any fixed $T,X,N$ and $J$ as in Lemma \ref{lemma:cubicest}, the following estimate holds:\begin{align}\label{eq:quadestimate}
        \Delta x \sum_{n =0}^{N-1} \sum_{j \in \mathbb{Z}}\nu_{\jph}^{n} (\Delta u_{\jph}^n)^2 \leq K(X,T),
    \end{align}
where $K(X,T)$ is a constant independent of $\Delta x.$
\begin{proof}
Recall the expression \eqref{eq:numentropyprod}: 
    \begin{align}\label{eq:Ejph_quad_estimate}
    \bar{E}_{\jph}^{n} 
         & = I+\bar{J} + R_{1}^n+ R_{2}^n + \tilde{R}_{1}^n+ \tilde{R}_{2}^n,
    \end{align}
where $I$ and $\bar{J}$ is as in \eqref{eq:I_b} and \eqref{eq:J}, respectively.    
Choosing $S(k,u) = f(k,u) - f(k,c)$ in \eqref{eq:I_b} and using the fact that $\bar{g}_{\jph}^\prime= \frac{\Delta g_{\jph}^n}{\Delta u_{\jph}^n},$ we write
\begin{align*}
     I 
    &= -(\Delta u_{\jph}^n)^2\left(\frac{1}{4}- \left(\lambda\frac{\Delta g_{\jph}^n}{\Delta u_{\jph}^n}\right)^{2}\right)\int_{0}^{1}\int_{0}^{1}s f_{uu}(k_{\jph}, u_{\jph}(r,s))\dif s \dif r.
\end{align*} Now, owing to hypothesis \ref{hyp:H2}, we expand $f_{uu}$ in the second variable about the point $(k_{\jph}, \tilde{u}_{\jph}^n)$ using a Taylor series, and write the term $I$ as 
\begin{align}\label{eq:I_quadest}
    I 
    & =-\frac{1}{8}(\Delta u_{\jph}^{n})^2\left(1- \left(2\lambda\frac{\Delta g_{\jph}^{n}}{\Delta u_{\jph}^{n}}\right)^{2}\right) f_{uu}(k_{\jph}, \tilde{u}_{\jph}^{n}) + \mathcal{O}((\Delta u_{\jph})^{3}) \\ & \spc + \mathcal{O}(\Delta k_{\jph}), \notag
\end{align}
where $\tilde{u}_{\jph}^{n}$ is defined as in \eqref{eq:ujphdef}. The above expression uses the fact that
\begin{align*}
   u_{\jph}(r,s) - \tilde{u}_{\jph}^{n} &= \frac{1}{2}\Delta u_{\jph}^{n}\left(1-s(1+r)\right) - 
    \lambda \frac{\Delta g_{\jph}^n}{\Delta u_{\jph}^n}s(1-r) \Delta u_{\jph}^{n} \\&= \mathcal{O}(\Delta u_{\jph}^n) + \mathcal{O}(\Delta k_{\jph}).
\end{align*} 
Further, the choice $S(k,u) = f(k,u) - f(k,c)$ and the trapezoidal rule of integration along with  \eqref{eq:fminusg_taylor} simplify $\bar{J}$ in \eqref{eq:J} to
\begin{align*}
   \bar{J}
   &=  \lambda\int_{0}^{1}f_{uu}(k_{\jph},u(s))\left(f(k_{\jph}, u(s))- \bar{g}_{\jph}(u(s)\right) u^\prime(s)\dif s\\
   &= -\frac{\lambda}{2} \Delta u_{\jph}^{n}\left[ f_{uu}(k_{\jph},u_{\jpo}^{n})\left(f(k_{\jph},u_{\jpo}^{n})-g(k_{\jpo},u_{\jpo}^{n})\right) \right.\\
   &\spc + \left. f_{uu}(k_{\jph},u_{j}^{n})\left(f(k_{\jph},u_{j}^{n})-g(k_{j},u_{j}^{n})\right) \right] + \mathcal{O}((\Delta u_{\jph}^{n})^{3})\\
   & =  -\frac{\lambda}{2} \Delta u_{\jph}^{n}\left[ f_{uu}(k_{\jph},u_{\jpo}^{n})\left(\frac{\lambda}{2}(a_{j+1}^n)^2\sigma_{j+1}^n- \frac{1}{8}(\lambda a_{j+1}^n\sigma_{j+1}^n)^{2}f_{uu}(k_{j+1},\zeta_{3}) \right. \right. \\ & \spc \spc \left. \left. - \frac{1}{8\lambda}\sigma_{j+1}^n\right) \right.\\
   &\spc  + \left. f_{uu}(k_{\jph},u_{j}^{n})\left( \frac{\lambda}{2}(a_{j}^n)^2\sigma_{j}^n- \frac{1}{8}(\lambda a_{j}^n\sigma_{j}^n)^{2}f_{uu}(k_j,\zeta_{2}) - \frac{1}{8\lambda}\sigma_{j}^n \right) \right] \\ & \spc  + \mathcal{O}((\Delta u_{\jph}^n)^{3}) +  \mathcal{O}(\Delta k_{\jph})+ \mathcal{O}(k_{\jph}-k_j) + \mathcal{O}(k_{\jph}-k_{j+1}).
\end{align*}
Again, noting that  $\sigma_{j}^n, \sigma_{\jpo}^n$ are $\mathcal{O}(\Delta u_{\jph}^n),$ replacing  $(a_{j}^n)^2$ and $(a_{j+1}^n)^2$  from their expression in \eqref{eq:aj_sq&ajph_sq}, and subsequently expanding $f_{uu}$ in a Taylor series about the point $(k_\jph, \tilde{u}_\jph^n),$ the term $\bar{J}$ reduces to
\begin{align}\label{eq:Jinquadest}
     \bar{J} &=\frac{1}{16}\Delta u_{\jph}^{n}(1-4(\beta_\jph^n)^{2})\left[ f_{uu}(k_{\jph},u_{\jpo}^{n})\sigma_{j+1}^{n}+  f_{uu}(k_{\jph},u_{j}^{n})\sigma_{j}^{n} \right] \\ & \spc + \mathcal{O}((\Delta u_{\jph}^n)^{3}) +  \mathcal{O}(\Delta k_{\jph}) + \mathcal{O}(k_{\jph}-k_j) + \mathcal{O}(k_{\jph}-k_{j+1})\notag\\
     &= \frac{1}{16}(\Delta u_{\jph}^{n})^{2}(1-4(\beta_\jph^n)^{2})f_{uu}(k_{\jph},\tilde{u}_{\jph}^{n})\frac{(\sigma_{j+1}^{n}+\sigma_{j}^{n})}{\Delta u_{\jph}^n} + \mathcal{O}((\Delta u_{\jph}^n)^{3})\notag\\ & \spc +  \mathcal{O}(\Delta k_{\jph})+ \mathcal{O}(k_{\jph}-k_j) + \mathcal{O}(k_{\jph}-k_{j+1}).\notag
\end{align}
Now, we note that the term $$\displaystyle-\frac{1}{8}(\Delta u_{\jph}^{n})^2\left(1- \left(2\lambda\frac{\Delta g_{\jph}^{n}}{\Delta u_{\jph}^{n}}\right)^{2}\right),$$
which appears in \eqref{eq:I_quadest}, is nothing but the value of $I$ obtained with the choice $S(k,u) = \frac{u^2}{2},$ see \eqref{eq:I}. Therefore, we replace this term in \eqref{eq:I_quadest} with the right-hand side of \eqref{eq:I_a}. The resulting $I$ of \eqref{eq:I_quadest} is added with the term $\bar{J}$ in \eqref{eq:Jinquadest}, to yield
\begin{nalign}\label{eq:IplusJbar_quad}
    I+\bar{J} 
    & = -\frac{1}{8}(1-4(\beta_{\jph}^n)^2)\left[1-\frac{1}{16}(1-4(\beta_{\jph}^n)^2)\left(\frac{\Delta \sigma_{\jph}^{n}}{\Delta u_{\jph}^{n}}\right)^2\right. \\ & \spc\spc  \left.- (\beta_{\jph}^n+\mathcal{A}_{4})\left(\frac{\Delta \sigma_{\jph}^n}{\Delta u_{\jph}^{n}}\right)-\frac{(\sigma_{j}^n + \sigma_{j+1}^n)}{2\Delta u_{\jph}^n}\right]f_{uu}(k_{\jph}, \tilde{u}_{\jph}^{n})(\Delta u_{\jph}^{n})^2 \\& \spc + \mathcal{O}((\Delta u_{\jph}^{n})^3)  + \mathcal{O}(\Delta k_{\jph}) + \mathcal{O}(k_{\jph}-k_j) + \mathcal{O}(k_{\jph}-k_{j+1}),
\end{nalign} where $\mathcal{A}_4$ is as in \eqref{eq:A_4}. Since $\mathcal{A}_4=\mathcal{O}(\Delta u_{\jph}^n)$ (see  \eqref{eq: A4bound1}), the term above turns into 
\begin{nalign}\label{eq:IpJbarJ_fin} 
    I+ \bar{J}
    &= -\nu_{\jph}^n(\Delta u_{\jph}^n)^2 + \mathcal{O}((\Delta u_{\jph}^n)^{3}) +  \mathcal{O}(\Delta k_{\jph}),
\end{nalign}
where $\nu_{\jph}^n$ is given by \eqref{eq:nujph}.
Now, to show that $\nu_{\jph}^n \geq 0,$ first we first observe that $1-4(\beta_{\jph}^n)^2 \geq 0,$ by the CFL condition \eqref{eq:cfl_cubicest} and $f_{uu}(k_{\jph}, \tilde{u}_{\jph}^{n}) \geq 0,$ by hypothesis \ref{hyp:H2}. This, together with an argument similar to \eqref{eq:etageqzero} now implies that $\nu_{\jph}^n \geq 0.$ Finally, observing that $R_1^n+ R_2^n+ \tilde{R}_{1}^n+\tilde{R}_{2}^n=\mathcal{O}(k_{\jph}-k_j)+ \mathcal{O}(k_{\jph}-k_\jpo)$ (see \eqref{eq:R1R2} and \eqref{eq:R1R2tilde}) and calling the expression \eqref{eq:IpJbarJ_fin} in \eqref{eq:numentropyprod}, we write
\begin{nalign}\label{eq:entrprod_quadest}
    \bar{E}_{\jph}^{n}&= -\nu_{\jph}^n(\Delta u_{\jph}^n)^2 + \mathcal{O}((\Delta u_{\jph}^n)^{3}) +  \mathcal{O}(\Delta k_{\jph})\\& \spc + \mathcal{O}(k_{\jph}-k_j)+ \mathcal{O}(k_{\jph}-k_\jpo).
\end{nalign}
Recalling the notation (see \eqref{eq:numentropyflux}),
\begin{align*}
  \bar{E}_{\jph}^{n}&:=S(k_{\jph},u^{n+1}_{\jph}) - \frac{1}{2}\left( S(k_{j},u_{j}^n)+S(k_{j+1} ,u_{j+1}^n)\right) \\& \spc + \lambda\left(G(k_{j+1}, u_{j+1}^n)- G(k_{j},u_{j}^n)\right),  
\end{align*}
summing \eqref{eq:entrprod_quadest} over $n= 0,1, \dots, N-1,$ $j+\frac{n}{2} \in \mathbb{Z}$ with $\abs{j} \leq J,$ and multiplying by $\Delta x,$ we obtain
\begin{align}\label{eq:quadest_prefinal}
    &\Delta x\sum_{n=0}^{N-1}\sum_{\substack{\abs{j} \leq J \\ j+\frac{n}{2} \in \mathbb{Z}}}\nu_{\jph}^n(\Delta u_{\jph}^n)^2 
    \\ &\leq \Delta x\sum_{\substack{\abs{j} \leq J \\ j+\frac{n}{2} \in \mathbb{Z}}}\frac{1}{2}\left( S(k_{j},u_{j}^{0})+S(k_{j+1} ,u_{j+1}^{0})\right)-\Delta x\sum_{\substack{\abs{j} \leq J \\ j+\frac{n}{2} \in \mathbb{Z}}}S(k_{\jph},u^{N}_{\jph})\notag \\ & \spc + \Delta t\sum_{n=0}^{N-1}G(k_{-J},u_{-J}^{n})  - \Delta t\sum_{n=0}^{N-1}G(k_{J}, u_{J}^{n})
   \notag \\ & \spc + K_{1}\Delta x\sum_{n=0}^{N-1}\sum_{\substack{\abs{j} \leq J \\ j+\frac{n}{2} \in \mathbb{Z}}}\abs{\Delta u_{\jph}^{n}}^{3} + K_{2}N\Delta x\norm{k}_{BV}, \notag
\end{align}
for some constants $K_{1},K_{2} \geq 0.$ Finally, since $\norm{S} \leq 2\norm{f}$ for the choice  $S(k,u) = f(k,u) - f(k,c),$  using the cubic estimate \eqref{eq:C_cubiclemma} from Lemma \ref{lemma:cubicest}  and the boundedness of $G$ in \eqref{eq:quadest_prefinal}, we obtain the desired estimate \eqref{lemma:quadest}, with $$K(X,T) := 8X\norm{f}+ 2 \norm{G}T + K_{1}C(X,T) + \frac{1}{\lambda}K_{2}\norm{k}_{BV}T,$$ 
and $C(X,T)$ is as in  \eqref{eq:C_cubiclemma}, thus completing the proof.
\end{proof}
\end{lemma}
\section{Convergence of the second-order scheme}\label{section:weaksolconv}
As a first step in proving convergence to a weak solution, we provide the $\mathrm{W}^{-1,2}_{\mathrm{loc}}$ compactness in the following lemma, which is necessary for applying the compensated compactness result in Theorem \ref{lemma:compcompactness}.

\begin{lemma}[$\mathrm{W}^{-1,2}_{\mathrm{loc}} \mbox{compactness}$]\label{lemma:W-12compactness}
Let $u_{0} \in (\mathrm{L}^{\infty}\cap \mathrm{BV})(\mathbb{R}).$ Under the CFL condition \eqref{eq:cfl_cubicest}, for the approximate solutions $u_{\Delta}$  in \eqref{eq:pcsoln} and the pairs $(S_{i},Q_{i})$ considered in \eqref{eq:entropypairs}, the sequence of distributions
\[
\left\{ S_{i}(k(x), u_{\Delta})_t + Q_{i}(k(x), u_{\Delta})_x \right\}_{\Delta > 0},
\]
is contained in a compact subset of \( \mathrm{W}^{-1,2}_{\mathrm{loc}}(\mathbb{R} \times \mathbb{R}^{+}) \), for $i=1,2.$

\begin{proof}
        
Suppressing the index of $(S_{i}, Q_{i}),$ we denote $$\mathcal{L}^{\Delta
}:=-(S(u^\Delta)_t + Q(k(x), u^\Delta)_x)$$ which is defined by the action 
\begin{align}\label{eq:Ldelta}
    \left\langle \mathcal{L}^{\Delta}, \phi \right\rangle &:= \int_{\mathbb{R_{+}}}\int_{\mathbb{R}} \left( S(k(x),u_{\Delta})\phi_t + Q(k(x),u_{\Delta})\phi_x \right) \dif x \dif t,
\end{align}
for $\phi \in W^{1,2}_{0,\mathrm{loc}}(\mathbb{R}\times \mathbb{R}_{+}).$
Using $k_{\Delta}$ as defined in \eqref{eq:pcsoln}, we add and subtract $S(k_{\Delta}(x), u_{\Delta})$ and $Q(k_{\Delta}(x), u_{\Delta})$ in the integrand of \eqref{eq:Ldelta}, which in turn allows us to write $\mathcal{L}^{\Delta} = \hat{\mathcal{L}}^{\Delta} + \tilde{\mathcal{L}}^{\Delta},$ where  
\begin{align*}
    \langle \hat{\mathcal{L}}^{\Delta}, \phi \rangle &:= \int_{\mathbb{R}^+} \int_{\mathbb{R}} \left( S(k(x), u_{\Delta}) - S(k_{\Delta}(x), u_{\Delta}) \right) \phi_t \, \dif x \, \dif t \nonumber \\
    &\spc + \int_{\mathbb{R}^+} \int_{\mathbb{R}} \left( Q(k(x), u_{\Delta}) - Q(k_{\Delta}(x), u_{\Delta}) \right) \phi_x \, \dif x \, \dif t, \\
    \langle \tilde{\mathcal{L}}^{\Delta}, \phi \rangle &:= \int_{\mathbb{R}^+} \int_{\mathbb{R}} \left( S(k_{\Delta}(x), u_{\Delta}) \phi_t + Q(k_{\Delta}(x), u_{\Delta}) \phi_x \right) \, \dif x \, \dif t.
\end{align*}
Now, we consider a bounded open subset $\Omega$ of $\mathbb{R}\times \mathbb{R}_{+}$ and let $X > 0,$ $T > 0$ be such that $\Omega \subseteq [-X, X] \times [0, T].$ Further, choose smallest integers $J,N \in \mathbb{N}$ such that $J\Delta x > X+\Delta x,$ and $N\Delta t > T.$ Additionally, define  $S_{j}^{n}:=S(k_{\Delta}(x_{j},t^{n}), u_{\Delta}(x_{j},t^{n}))= S(k_{j}^{n}, u_{j}^{n}),$ $Q_{j}^{n}:=Q(k_{\Delta}(x_{j},t^{n}), u_{\Delta}(x_{j},t^{n}))= Q(k_{j}^{n}, u_{j}^{n}).$  

Let $q \in (1,2)$ and  $\displaystyle p=\frac{q}{q-1}.$ Now, applying H\"{o}lder's inequality to the term $\langle \hat{\mathcal{L}}^{\Delta}, \phi \rangle$ with $\phi \in W_{0}^{1,p}(\Omega),$ we obtain
\begin{align*}
    \lim_{\Delta \rightarrow 0} \abs{\langle \hat{\mathcal{L}}^{\Delta}, \phi \rangle} \leq \lim_{\Delta \rightarrow 0}\left(\norm{\partial_{k}S}+\norm{\partial_{k}Q}\right)\norm{k-k_{\Delta}}_{L^{q}(\Omega)}\norm{\phi}_{W_{0}^{1,p}(\Omega)} = 0,
\end{align*}
from which we can conclude that
\begin{nalign}\label{eq:Lhatcompact}
    \{\hat{\mathcal{L}}^{\Delta}\} \,  \mbox{is compact in} \, \mathrm{W}^{-1,q}(\Omega) \,\,  \mbox{for any} \, q \in (1,2).
\end{nalign}
\par
Next, we focus on showing that  $\{\tilde{\mathcal{L}}^{\Delta}\}_{\Delta > 0} \,  \mbox{is compact in}\,  \mathrm{W}^{-1,q}(\Omega) \,  \mbox{for some} \, q \in (1,2).
$ We begin with applying summation by parts to expand the term $\left\langle \tilde{\mathcal{L}}^{\Delta}, \phi \right\rangle $ as
\begin{align*}
\left\langle \tilde{\mathcal{L}}^{\Delta}, \phi \right\rangle &=\sum_{n=0}^{N-1}\sum_{\substack{\abs{j} \leq J \\ j+\frac{n}{2} \in \mathbb{Z}}}\int_{t^{n}}^{t^{n+1}} \int_{x_{j-\frac{1}{2}}}^{x_{j+\frac{1}{2}}}  \left( S(k_{\Delta}(x), u_{\Delta}) \phi_t + Q(k_{\Delta}(x), u_{\Delta}) \phi_x \right) \, \dif x \, \dif t.\\
&= \sum_{n=0}^{N-1}\sum_{\substack{\abs{j} \leq J \\ j+\frac{n}{2} \in \mathbb{Z}}}\Bigg[ \int_{x_{j-\frac{1}{2}}}^{x_{j+\frac{1}{2}}} 
S_j^n \left( \phi(x, t_{n+1}) - \phi(x, t_n) \right) \dif x \\
&\spc + \int_{t_n}^{t_{n+1}} Q_j^n \left( \phi(x_{j+\frac{1}{2}}, t) - \phi(x_{j-\frac{1}{2}}, t) \right) \dif t \Bigg] 
\\
&= \sum_{\substack{\abs{j} \leq J \\ j+\frac{n}{2} \in \mathbb{Z}}} \int_{x_{j-\frac{1}{2}}}^{x_\jph} S_j^{N-1}  \phi(x, N\Delta t) \dif x- \sum_{\substack{\abs{j} \leq J \\ j+\frac{n}{2} \in \mathbb{Z}}} \int_{x_{j-\frac{1}{2}}}^{x_\jph} S_j^{0}  \phi(x, 0) \dif x\\
&\spc- \sum_{n=1}^{N-1}\sum_{\substack{\abs{j} \leq J \\ j+\frac{n}{2} \in \mathbb{Z}}} \Bigg[ \int_{x_{j-\frac{1}{2}}}^{x_j} \left( S_j^n - S_{j-\frac{1}{2}}^{n-1} \right) \phi(x, t_n) \dif x \\ & \spc \spc 
+ \int_{x_j}^{x_{j+\frac{1}{2}}} \left( S_j^n - S_{j+\frac{1}{2}}^{n-1} \right) \phi(x, t_n) \dif x \Bigg] \\
&
\spc - \sum_{n=0}^{N-1}\sum_{\substack{\abs{j} \leq J \\ j+\frac{n}{2} \in \mathbb{Z}}} \int_{t^{n}}^{t_{n+1}} \left( Q_{j+1}^{n} - Q_{j}^{n} \right) \phi(x_\jph, t) \dif t.
\end{align*}
Further,  denoting $\phi_j^{n}:= \phi(x_{j}, t^{n})$ and adding and subtracting suitable terms, we write 
\begin{align}\label{eq:Lphi}
\left\langle \tilde{\mathcal{L}}^{\Delta}, \phi \right\rangle 
&=:< \tilde{\mathcal{L}}^{\Delta}_{0}, \phi> + < \tilde{\mathcal{L}}^{\Delta}_{1}, \phi> + < \tilde{\mathcal{L}}^{\Delta}_{2}, \phi> + < \tilde{\mathcal{L}}^{\Delta}_{3}, \phi> + < \tilde{\mathcal{L}}^{\Delta}_{4}, \phi>,
\end{align}
where 
\begin{align}
    < \tilde{\mathcal{L}}^{\Delta}_{0}, \phi> &:= \sum_{\substack{\abs{j} \leq J \\ j+\frac{n}{2} \in \mathbb{Z}}} \int_{x_{j-\frac{1}{2}}}^{x_\jph} S_j^{N-1}  \phi(x, N\Delta t) \dif x- \sum_{\substack{\abs{j} \leq J \\ j+\frac{n}{2} \in \mathbb{Z}}} \int_{x_{j-\frac{1}{2}}}^{x_\jph} S_j^{0}  \phi(x, 0) \dif x,\\
    < \tilde{\mathcal{L}}^{\Delta}_{1}, \phi>&:=\sum_{n=1}^{N-1}\sum_{\substack{\abs{j} \leq J \\ j+\frac{n}{2} \in \mathbb{Z}}} \int_{x_{j-\frac{1}{2}}}^{x_j} \left( S_j^n - S_{j-\frac{1}{2}}^{n-1} \right) \left( \phi_j^n - \phi(x, t_n) \right) \dif x, \notag   \\
    < \tilde{\mathcal{L}}^{\Delta}_{2}, \phi> &:= \sum_{n=1}^{N-1}\sum_{\substack{\abs{j} \leq J \\ j+\frac{n}{2} \in \mathbb{Z}}} \int_{x_j}^{x_{j+\frac{1}{2}}} \left( S_j^n - S_{j+\frac{1}{2}}^{n-1} \right) \left( \phi_j^n - \phi(x, t_n) \right) \dif x, \notag\\
    < \tilde{\mathcal{L}}^{\Delta}_{3}, \phi> &:= \sum_{n=1}^{N}\sum_{\substack{\abs{j} \leq J \\ j+\frac{n-1}{2} \in \mathbb{Z}}} \int_{t_{n-1}}^{t_{n}} \left( Q_{j+\frac{1}{2}}^{n-1} - Q_{j-\frac{1}{2}}^{n-1} \right) \left( \phi_j^n - \phi(x_j, t) \right) \dif t, \notag\\
     < \tilde{\mathcal{L}}^{\Delta}_{4}, \phi>&:=  -\Delta x \sum_{n=1}^{N-1}\sum_{\substack{\abs{j} \leq J \\ j+\frac{n}{2} \in \mathbb{Z}}} \phi_j^n \left( S_j^n - \frac{1}{2} \left( S_{j+\frac{1}{2}}^{n-1} + S_{j-\frac{1}{2}}^{n-1} \right) + \lambda \left( Q_{j+\frac{1}{2}}^{n-1} - Q_{j-\frac{1}{2}}^{n-1} \right) \right). \notag
\end{align}

Hereafter, we set $S(k,u)=S_2(k,u) = f(k,u) - f(k,c), \,\, Q(k,u)=Q_2(k,u) = \int_c^u (f_u(k,\xi))^2 \, \dif \xi.$ The proof for the case $S_1(u) = u - c, \,\, Q_1(k,u) = f(k,u) - f(k,c)$ follows analogously and is omitted.
 Now, using a change of index, we write the term $< \tilde{\mathcal{L}}^{\Delta}_{4}, \phi>$ in \eqref{eq:Lphi} as
\begin{nalign}\label{eq:L4}
    < \tilde{\mathcal{L}}^{\Delta}_{4}, \phi> &=\Delta x \sum_{n=0}^{N-2}\sum_{\substack{\abs{j} \leq J \\ j+\frac{n+1}{2} \in \mathbb{Z}}} \phi_{\jph}^{n+1} \left[S_{\jph}^{n+1} - \frac{1}{2} \left( S_{j+1}^{n} + S_{j}^{n} \right) + \lambda \left( Q_{j+1}^{n} - Q_{j}^{n} \right) \right].
\end{nalign}
Proceeding as in the proof of Lemma \ref{lemma:cubicest} (see \eqref{eq:Ejph}-\eqref{eq:I_b}) and using the same notations,  we can write
\begin{align}\label{eq: entropypro_compact}
    &S_{\jph}^{n+1} - \frac{1}{2} \left( S_{j+1}^{n} + S_{j}^{n} \right) + \lambda \left( Q_{j+1}^{n} - Q_{j}^{n} \right)\\
   & = I+ \bar{J}-\lambda\Bigl(S_{u}(k_{j}, u_j^n)(f(k_{j}, u_j^n) - g(k_{j}, u_{j}^n)) \notag\\ & \spc -S_{u}(k_{j+1}, u_{j+1}^n)(f(k_{j+1}, u_{j+1}^n) - g(k_{j+1}, u_{j+1}^n))\Bigr)+ \tilde{R}_1 +\tilde{R}_2+ R_{1}+ R_{2}. \notag
\end{align}
Now, invoking \eqref{eq: entropypro_compact} in \eqref{eq:L4}, we write 
\begin{nalign}\label{eq:L4_b}
    < \tilde{\mathcal{L}}_{4}^{\Delta}, \phi> &= < \tilde{\mathcal{L}}_{4,1}^{\Delta}, \phi> + < \tilde{\mathcal{L}}_{4,2}^{\Delta}, \phi> + < \tilde{\mathcal{L}}^{\Delta}_{4,3}, \phi>,
\end{nalign}
where we define
\begin{align}\label{eq:L4splitting}
 < \tilde{\mathcal{L}}^{\Delta}_{4,1}, \phi>&:=    \sum_{n=0}^{N-2}\sum_{\substack{\abs{j} \leq J \\ j+\frac{n+1}{2} \in \mathbb{Z}}} \Delta x \phi_{\jph}^{n+1}\left(I+ \bar{J}\right), \\
<\tilde{\mathcal{L}}^{\Delta}_{4,2}, \phi>&:=  -\sum_{n=0}^{N-2}\sum_{\substack{\abs{j} \leq J \\ j+\frac{n+1}{2} \in \mathbb{Z}}} \Delta x \phi_{\jph}^{n+1} \lambda\Bigl(S_{u}(k_{j}, u_j^n)(f(k_{j}, u_j^n) - g(k_{j}, u_{j}^n))\notag\\ & \spc \spc -S_{u}(k_{j+1}, u_{j+1}^n)(f(k_{j+1}, u_{j+1}^n) - g(k_{j+1}, u_{j+1}^n))\Bigr),\notag\\  
    < \tilde{\mathcal{L}}^{\Delta}_{4,3}, \phi> \notag&:=  \sum_{n=0}^{N-2}\sum_{\substack{\abs{j} \leq J \\ j+\frac{n+1}{2} \in \mathbb{Z}}} \Delta x \phi_{\jph}^{n+1} \left(\tilde{R}_1 +\tilde{R}_2+ R_{1}+ R_{2}\right). \notag
\end{align}
At this stage, by summarizing \eqref{eq:Lphi} and \eqref{eq:L4splitting}, we write 
\begin{align}
\left\langle \tilde{\mathcal{L}}^{\Delta}, \phi \right\rangle 
&=:< \tilde{\mathcal{L}}^{\Delta}_{0}, \phi> + < \tilde{\mathcal{L}}^{\Delta}_{1}, \phi> + < \tilde{\mathcal{L}}^{\Delta}_{2}, \phi> + < \tilde{\mathcal{L}}^{\Delta}_{3}, \phi> \\ & \spc +< \tilde{\mathcal{L}}_{4,1}^{\Delta}, \phi> + < \tilde{\mathcal{L}}_{4,2}^{\Delta}, \phi> + < \tilde{\mathcal{L}}^{\Delta}_{4,3}, \phi>, \notag
\end{align} and  proceed with the proof of compactness of $\tilde{\mathcal{L}}^{\Delta}$ in two steps. \\
\textbf{Step 1} (Compactness of $ \tilde{\mathcal{L}}^{\Delta}_{4,1}, \tilde{\mathcal{L}}^{\Delta}_{4,3}
\, \,\mbox{and}\, \,  \tilde{\mathcal{L}}^{\Delta}_{0} $):
Using \eqref{eq:IpJbarJ_fin} from Lemma \ref{lemma:quadest}, we write
\begin{nalign}\label{eq:L41}
    &<\tilde{\mathcal{L}}^{\Delta}_{4,1}, \phi> \\
    & \spc = \sum_{n=0}^{N-2}\sum_{\substack{\abs{j} \leq J \\ j+\frac{n+1}{2} \in \mathbb{Z}}} \Delta x \phi_{\jph}^{n+1}\left(-\nu_{\jph}^n(\Delta u_{\jph}^n)^2 + \mathcal{O}((\Delta u_{\jph}^n)^{3}) +  \mathcal{O}(\Delta k_{\jph})\right).
\end{nalign}
Again, the cubic estimate \eqref{eq:C_cubiclemma} from Lemma \ref{lemma:cubicest} and the quadratic estimate \eqref{eq:quadestimate} from Lemma \ref{lemma:quadest}, together with the hypothesis \ref{hyp:H1} yield
\begin{align}\label{eq:L41bd}
    \abs{<\tilde{\mathcal{L}}^{\Delta}_{4,1}, \phi>} &\leq \Delta x \norm{\phi}\sum_{n=0}^{N-2}\sum_{\substack{\abs{j} \leq J \\ j+\frac{n+1}{2} \in \mathbb{Z}}} \nu_{\jph}^n(\Delta u_{\jph}^n)^2 \\ & \spc + \mathcal{K}_{1} \Delta x \norm{\phi}\sum_{n=0}^{N-2}\sum_{\substack{\abs{j} \leq J \\ j+\frac{n+1}{2} \in \mathbb{Z}}} \abs{\Delta u_{\jph}^n}^{3} +\mathcal{K}_{2}\Delta x \norm{\phi}\sum_{n=0}^{N-2}\sum_{\substack{\abs{j} \leq J \\ j+\frac{n+1}{2} \in \mathbb{Z}}}  \abs{\Delta k_{\jph}}\notag\\ & \leq \mathcal{K}_{3}\norm{\phi}, \quad \mbox{for} \quad \phi \in C_{0}(\Omega),\notag
\end{align} and some constants $\mathcal{K}_1, \mathcal{K}_2 >0$ and $\displaystyle K_{3} := K(X,T) + \mathcal{K}_{1}C(X,T) + \frac{1}{\lambda}\mathcal{K}_{2}\norm{k}_{BV}T.$ Here,  $C_{0}(\Omega)$ denotes the space of continuous functions on $\Omega$ that vanish at the boundary.
\par 
In the next step, we estimate $ < \tilde{\mathcal{L}}^{\Delta}_{4,3}, \phi>$ 
for $\phi \in C_0(\Omega)$ by recalling \eqref{eq:R1R2}, \eqref{eq:R1R2tilde} and using hypothesis \ref{hyp:H1}, as
\begin{align}\label{eq:L43bd}
    \abs{<\tilde{\mathcal{L}}^{\Delta}_{4,3}, \phi>} \leq  \norm{\phi}\mathcal{K}_{4} \sum_{n=0}^{N-2} \Delta x \norm{k}_{\mathrm{BV}} \leq \mathcal{K}_{5}\norm{\phi},
\end{align}
where $\displaystyle \mathcal{K}_{5}:= \norm{k}_{BV}\frac{\mathcal{K}_{4}}{\lambda} T$  and $\mathcal{K}_{4} >0$ is some constant. Combining the bounds \eqref{eq:L41bd} and \eqref{eq:L43bd}, we obtain the following estimate

\begin{align}\label{eq:L41L43bd}
    \norm{\tilde{\mathcal{L}}^{\Delta}_{4,1}}_{\mathcal{M}(\Omega)}, \norm{\tilde{\mathcal{L}}^{\Delta}_{4,3}}_{\mathcal{M}(\Omega)} \leq \mathcal{K}_{6}= \max\{\mathcal{K}_{3}, \mathcal{K}_{5}\},
\end{align}
 as $\mathcal{M}(\Omega),$ the space of all bounded Radon measures on $\Omega,$ is the dual of the space $\left(C_{0}(\Omega), \norm{\cdot}_{\infty}\right)$ ( see \cite{malek1996} for more details). Similarly, dealing with the term $< \tilde{\mathcal{L}}^{\Delta}_{0}, \phi>, $ we obtain a bound
\begin{align}\label{eq:L0bd}
   \norm{\tilde{\mathcal{L}}^{\Delta}_{0}}_{\mathcal{M}(\Omega)} \leq \mathcal{K}_{7}, 
\end{align}
$\mbox{as} \,   \abs{< \tilde{\mathcal{L}}^{\Delta}_{0}, \phi>} \leq \mathcal{K}_{7}\norm{\phi}$ for some constant $\mathcal{K}_{7}>0.$
By Sobolev's imbedding theorem (see Lemma 2.55, page 38 in \cite{malek1996}), we have the inclusion $\mathcal{M}(\Omega) \subset \mathrm{W}^{-1, q}(\Omega)$ with compact imbedding for $q \in [1,2).$ Consequently, from \eqref{eq:L41L43bd} and \eqref{eq:L0bd}, it follows  that the sets 
\begin{nalign}
    \{\tilde{\mathcal{L}}^{\Delta}_{4,1}\}_{\Delta >0}, \{\tilde{\mathcal{L}}^{\Delta}_{4,3}\}_{\Delta >0} \, \,\mbox{and}\, \,  \{\tilde{\mathcal{L}}^{\Delta}_{0}\}_{\Delta > 0} \, \mbox{are compact in} \, \mathrm{W}^{-1,q}(\Omega), \, \forall\, q \in [1,2).
\end{nalign}
\noindent
\textbf{Step 2} (Compactness of $\tilde{\mathcal{L}}^{\Delta}_{1}, \tilde{\mathcal{L}}^{\Delta}_{2}, \tilde{\mathcal{L}}^{\Delta}_{3} \,\mbox{and}\,\tilde{\mathcal{L}}^{\Delta}_{4,2}$):
We first derive some estimates useful in bounding the term $< \tilde{\mathcal{L}}^{\Delta}_{1}, \phi>$ in  \eqref{eq:Lphi}.
Observing that
\begin{nalign}
    \abs[\Big]{u^{n-\frac{1}{2}}_{\jph} - u^{n-\frac{1}{2}}_{\jmh}} &= \abs[\Big]{\Delta u_{j}^{n-1}-\frac{\lambda}{2}f_{u}(k_{\jph},u^{n-1}_{\jph}) \sigma_{\jph}^{n-1} + \frac{\lambda}{2}f_{u}(k_{\jmh},u^{n-1}_{\jmh})\sigma_{\jmh}^{n-1}}\\ &\leq \abs{\Delta u_{j}^{n-1}} \left(1+\kappa\right),
\end{nalign}
we derive the following bound
\begin{align*}
   &\abs{u_{j}^{n}-u_{\jmh}^{n-1}}\\&
   \leq \abs[\Big]{\half( u^{n-1}_{\jph}-u^{n-1}_{\jmh}) - \frac{1}{8}(\sigma_{\jph}^{n-1} -\sigma_{\jmh}^{n-1}) - \lambda\left(f(k_{\jph}, u^{n-1/2}_{\jph})-f(k_{\jmh}, u^{n-1/2}_{\jmh})\right)}\\
&  \leq \half\abs[\Big]{ u^{n-1}_{\jph}-u^{n-1}_{\jmh}} + \frac{1}{8}\left(\abs{\sigma_{\jph}^{n-1}}+\abs{\sigma_{\jmh}^{n-1}}\right)\\ & \spc  + \abs[\Big]{\lambda\left(f(k_{\jph}, u^{n-1/2}_{\jph})-f(k_{\jmh}, u^{n-1/2}_{\jmh})\right)}\\
&\leq \left(\frac{1}{2}+ \frac{1}{4}\right)\abs[\Big]{\Delta u_{j}^{n-1}}  + \lambda \norm{f_{k}}\abs{\Delta k_{j}} +\kappa \abs[\Big]{u^{n-\frac{1}{2}}_{\jph} - u^{n-\frac{1}{2}}_{\jmh}}\\
   &\leq  \abs{\Delta u_{j}^{n-1}}\left(\frac{3}{4}+ \kappa(1+\kappa)\right) + \lambda \norm{f_{k}} \abs{\Delta k_{j}}.
\end{align*}Again, using the above estimate, we obtain
\begin{align*}
     \abs{S_j^n - S_{j-\frac{1}{2}}^{n-1}} &\leq \norm{S_{k}} \abs{k_{j}-k_{\jmh}} + \norm{S_{u}}\abs{u_{j}^{n}-u_{\jmh}^{n-1}}\\
     &\leq \norm{S_{k}} \abs{k_{j}-k_{\jmh}} +\mathcal{K}_{8} \norm{S_{u}}\abs{\Delta u_{j}^{n-1}} + \norm{S_{u}}\lambda \norm{f_{k}}\abs{\Delta k_{j}},
\end{align*}
where $\displaystyle \mathcal{K}_{8}:=\frac{3}{4}+ \kappa(1+\kappa).$ Using this estimate in $< \tilde{\mathcal{L}}^{\Delta}_{1}, \phi>$ and subsequently applying the H\"{o}lder's inequality, we write
\begin{align*}
    &\abs{< \tilde{\mathcal{L}}^{\Delta}_{1}, \phi>} \\\notag &\leq \sum_{n=1}^{N-1}\sum_{\substack{\abs{j} \leq J \\ j+\frac{n}{2} \in \mathbb{Z}}} \mathcal{K}_{8} \norm{S_{u}}\abs{\Delta u_{j}^{n-1}}\int_{x_{j-\frac{1}{2}}}^{x_j}  \abs{\phi_j^n - \phi(x, t_n)} \dif x\\
    \notag& \spc + \sum_{n=1}^{N-1}\sum_{\substack{\abs{j} \leq J \\ j+\frac{n}{2} \in \mathbb{Z}}}\left(\norm{S_{k}}  \abs{k_{j}-k_{\jmh}} + \norm{S_{u}}\lambda \norm{f_{k}} \abs{\Delta k_{j}} \right)\int_{x_{j-\frac{1}{2}}}^{x_j}  \abs{\phi_j^n - \phi(x, t_n)} \dif x\\
    \notag&\leq \frac{1}{2}\mathcal{K}_{8}\sum_{n=1}^{N-1}\sum_{\substack{\abs{j} \leq J \\ j+\frac{n}{2} \in \mathbb{Z}}} \norm{S_{u}}\abs{\Delta u_{j}^{n-1}}\norm{\phi}_{C_{0}^{\alpha}}\Delta x^{\alpha+1}\\
    \notag& \spc + \frac{1}{2}\sum_{n=1}^{N-1}\sum_{\substack{\abs{j} \leq J \\ j+\frac{n}{2} \in \mathbb{Z}}}\left(\norm{S_{k}}  \abs{k_{j}-k_{\jmh}} + \norm{S_{u}}\lambda \norm{f_{k}} \abs{\Delta k_{j}} \right)\norm{\phi}_{C_{0}^{\alpha}}\Delta x^{\alpha+1}\\
    \notag&\leq \frac{1}{2}\mathcal{K}_{8}\norm{S_{u}}\norm{\phi}_{C_{0}^{\alpha}}\Delta x^{\alpha+\frac{2}{3}}\left(\Delta x\sum_{n=1}^{N-1}\sum_{\substack{\abs{j} \leq J \\ j+\frac{n}{2} \in \mathbb{Z}}} \abs{\Delta u_{j}^{n-1}}^{3}\right)^{\frac{1}{3}} \left(\sum_{n=1}^{N-1}\sum_{\substack{\abs{j} \leq J \\ j+\frac{n}{2} \in \mathbb{Z}}}1\right)^{\frac{2}{3}}\\
    \notag&\spc +\frac{1}{2\lambda}\left(\norm{S_{k}}T+\norm{S_{u}}\lambda \norm{f_{k}} T\right)\norm{k}_{BV}\norm{\phi}_{C_{0}^{\alpha}}\Delta x^{\alpha},
\end{align*}
for $\phi \in C_{0}^{\alpha}(\Omega),$ which denotes the space of H\"older continuous functions. Finally, invoking Lemma \ref{lemma:cubicest} in the previous inequality, we obtain
\begin{align}\label{eq:L1tildebd}
    \abs{< \tilde{\mathcal{L}}^{\Delta}_{1}, \phi>} 
    & \leq \frac{1}{2}C(X,T)^{\frac{1}{3}}\mathcal{K}_{8}\norm{S_{u}}\frac{\left(4XT\right)^{\frac{2}{3}}}{\lambda^{\frac{2}{3}}}\norm{\phi}_{C_{0}^{\alpha}}\Delta x^{\alpha-\frac{2}{3}}\\ \notag&\spc+\frac{1}{2\lambda}\left(\norm{S_{k}}T+\norm{S_{u}}\lambda \norm{f_{k}} T\right)\norm{k}_{BV}\norm{\phi}_{C_{0}^{\alpha}}\Delta x^{\alpha} \quad 
\mbox{for} \, \phi \in C_{0}^{\alpha}(\Omega).
\end{align}
Now, by the Sobolev's imbedding theorem  (see \cite{adamsfournier2003}) we have the inclusion  $W_{0}^{1,p}(\Omega) \hookrightarrow C_{0}^{\alpha}(\Omega)$ for  $\alpha \in \left(0, 1 - \frac{2}{p} \right].$ Consequently, for a fixed $\alpha_1 \in  (2/3, 1)$ and  $ p_{1} = \frac{2}{1-\alpha_1},$ there exists a constant $\mathcal{K}_{p_{1}} >0,$ such that \begin{align}\label{eq:phi_sob_embed}
    \norm{\phi}_{C_{0}^{\alpha_{1}}(\Omega)} \leq \mathcal{K}_{p_{1}} \norm{\phi}_{W_{0}^{1,p_{1}}(\Omega)} \,\, \forall \,\,\phi \in W_{0}^{1,p_{1}}(\Omega).
\end{align} 
Using \eqref{eq:phi_sob_embed} in \eqref{eq:L1tildebd} amounts to
\begin{nalign}\label{eq:L1tildebd_b}
     \abs{< \tilde{\mathcal{L}}^{\Delta}_{1}, \phi>} & \leq \frac{1}{2}C(X,T)^{\frac{1}{3}}\mathcal{K}_{p_{1}}  \mathcal{K}_{8}\norm{S_{u}}\frac{\left(4XT\right)^{\frac{2}{3}}}{\lambda^{\frac{2}{3}}}\norm{\phi}_{W_{0}^{1,p_{1}}(\Omega)}\Delta x^{\alpha_{1}-\frac{2}{3}}\\& \spc +\frac{1}{2\lambda}\mathcal{K}_{p_{1}} \left(\norm{S_{k}}T+\norm{S_{u}}\lambda \norm{f_{k}} T\right)\norm{k}_{BV}\norm{\phi}_{W_{0}^{1,p_{1}}(\Omega)}\Delta x^{\alpha_{1}}.
\end{nalign}
This implies that \begin{nalign}\label{eq:L1tildebd_c}
    \lim_{\Delta \rightarrow 0} \norm{\tilde{\mathcal{L}}^{\Delta}_{1}}_{\mathrm{W}^{-1,q_{1}}(\Omega)} \leq \lim_{\Delta \rightarrow 0} ( \mathcal{K}_{9} \Delta x^{\alpha_{1}-\frac{2}{3}}+ \mathcal{K}_{10} \Delta x^{\alpha_{1}}) = 0 \quad \mbox{for} \quad q_{1}=\frac{2}{1+\alpha_{1}}. 
\end{nalign}
where 
\begin{align*}
    \displaystyle \mathcal{K}_{9} &:= \frac{1}{2}C^{\frac{1}{3}}\mathcal{K}_{p_{1}} \mathcal{K}_{8}\norm{S_{u}}\frac{\left(4XT\right)^{\frac{2}{3}}}{\lambda^{\frac{2}{3}}}, \, 
    \mathcal{K}_{10} := \frac{1}{2\lambda}\mathcal{K}_{p_{1}} \left(\norm{S_{k}}T+\norm{S_{u}}\lambda \norm{f_{k}} T\right)\norm{k}_{BV}.
\end{align*}
From this, we conclude that $\{\tilde{\mathcal{L}}^{\Delta}_{1}\}_{\Delta > 0}$ is compact in  $\mathrm{W}^{-1,q_{1}}(\Omega),$ for $\displaystyle q_{1}= \frac{2}{1+\alpha_{1}} \in (1,2).$ In an entirely analogous way, we obtain that the set $\{\tilde{\mathcal{L}}^{\Delta}_{2}\}_{\Delta > 0}$ is compact in $\mathrm{W}^{-1,q_{1}}(\Omega).$
\par
Next, we consider the term $< \tilde{\mathcal{L}}^{\Delta}_{4,2}, \phi>$ from \eqref{eq:L4splitting}, and apply summation by parts to write
\begin{nalign}\label{eq:L42_ibp}
    < \tilde{\mathcal{L}}^{\Delta}_{4,2}, \phi>&= \lambda\Delta x\sum_{n=0}^{N-2}\sum_{\substack{\abs{j} \leq J \\ j+\frac{n+1}{2} \in \mathbb{Z}}}  \left(\phi_{\jmh}^{n+1}-\phi_{\jph}^{n+1}\right) S_{u}(k_{j}, u_{j}^{n})(f(k_{j}, u_{j}^{n}) - g(k_{j}, u_{j}^{n})).
\end{nalign}
From the expression \eqref{eq:gtaylor_j} in the proof of Lemma \ref{lemma:cubicest}, we deduce that
\begin{align}\label{eq:fming_bddelv}
\lambda\abs{f(k_{j},u_j^n)-g(k_j,u_j^{n})} 
& = \lambda\abs[\Big]{\frac{\lambda}{2}(a_{j}^n)^2\sigma_{j}^n- \frac{1}{8}(\lambda a_{j}^{n}\sigma_{j}^n)^{2}f_{uu}(k_j,\zeta_{2})- \frac{1}{8\lambda}\sigma_{j}^n} \\
&\leq \mathcal{K}_{11}\abs{\Delta u_{\jph}^{n}},\notag
\end{align}
where 
$ \displaystyle \mathcal{K}_{11} := \frac{\kappa^{2}}{2}+ \frac{1}{4}C_{u_{0}}\kappa^{2}\gamma_{2}\lambda+\frac{1}{8}.$
Inserting the estimate \eqref{eq:fming_bddelv} into \eqref{eq:L42_ibp} and observing that $\lambda \norm{S_{u}} \leq \lambda \norm{f_{u}} \leq \kappa,$ we  apply H\"{o}lder's inequality on \eqref{eq:L42_ibp} along with the cubic estimate \eqref{lemma:cubicest} of Lemma \ref{lemma:cubicest}, to obtain 
\begin{align}\label{eq:L42bd}
|
\langle \tilde{\mathcal{L}}^{\Delta}_{4,2}, \phi \rangle
| &\leq \kappa \mathcal{K}_{11}
\left(
\Delta x\sum_{n=0}^{N-2}\sum_{\substack{\abs{j} \leq J \\ j+\frac{n+1}{2} \in \mathbb{Z}}}   |\Delta \phi_{j}^{n}|^{3/2}
\right)^{2/3}
\left(
\Delta x \sum_{n=0}^{N-2}\sum_{\substack{\abs{j} \leq J \\ j+\frac{n+1}{2} \in \mathbb{Z}}}  |\Delta u_{\jph}^{n}|^3
\right)^{1/3} \\& \leq  \kappa \mathcal{K}_{11}C(X,T)^{\frac{1}{3}} \left(
\Delta x^{(1+\frac{3}{2}\alpha)}\norm{\phi}_{C^\alpha_0}^{\frac{3}{2}\alpha}\sum_{n=0}^{N-2}\sum_{\substack{\abs{j} \leq J \\ j+\frac{n+1}{2} \in \mathbb{Z}}}1   
\right)^{2/3} \notag\\
&\leq\kappa \mathcal{K}_{11}C(X,T)^{\frac{1}{3}} ||\phi||_{C^\alpha_0} (\Delta x)^{\alpha + \frac{2}{3}} \frac{\left(4XT\right)^{\frac{2}{3}}}{\lambda^{\frac{2}{3}}(\Delta x)^{\frac{4}{3}}}
\leq \mathcal{K}_{12}||\phi||_{C^\alpha_0} (\Delta x)^{\alpha - \frac{2}{3}}. \notag
\end{align}
for $\phi \in C_{0}^{\alpha}(\Omega)$, where $\displaystyle \mathcal{K}_{12}:=\kappa\mathcal{K}_{11}C(X,T)^{\frac{1}{3}} \left(\frac{4XT}{\lambda}\right)^{\frac{2}{3}}$ with $C(X,T)$ as in Lemma \ref{lemma:cubicest}. To get the penultimate inequality in the above estimate, we have used the facts that $X+\Delta x \leq J\Delta x \leq X+ 2\Delta x$ and $(N-1)\Delta x \leq T.$
Now, arguments similar to those in \eqref{eq:L1tildebd}, \eqref{eq:L1tildebd_b} and \eqref{eq:L1tildebd_c},  give compactness of  $\{\tilde{\mathcal{L}}^{\Delta}_{4,2}\}_{\Delta > 0}$ in  $\mathrm{W}^{-1,q_1}(\Omega),$ for the same $q_1 \in (1,2)$ as in the case of $\tilde{\mathcal{L}}^{\Delta}_{1}$ and $\tilde{\mathcal{L}}^{\Delta}_{2}.$
\par 
Next, we consider $\{\tilde{\mathcal{L}}^{\Delta}_{3}\},$  and obtain an estimate on the term ${< \tilde{\mathcal{L}}^{\Delta}_{3}, \phi>}$ using the observation 
\begin{align*}
    \abs{Q_{j+\frac{1}{2}}^{n-1} - Q_{j-\frac{1}{2}}^{n-1}} \leq \norm{Q_{k}}\abs{k_{\jph}-k_{\jmh}} + \norm{Q_{u}} \abs{\Delta u_{j}^{n}}, 
\end{align*}
as follows
\begin{align}
     &\abs{< \tilde{\mathcal{L}}^{\Delta}_{3}, \phi>}\\
     & \leq \lambda  \norm{\phi}_{C_{0}^{\alpha}} (\Delta x)^{\alpha+1}\norm{Q_{k}} \norm{k}_{BV}\left(\frac{2T}{\lambda \Delta x}\right) \notag \\ & \spc + \lambda\norm{\phi}_{C_{0}^{\alpha}} (\Delta x)^{\alpha+\frac{2}{3}} \norm{Q_{u}} \left(\Delta x\sum_{n=1}^{N-1}\sum_{\substack{\abs{j} \leq J \\ j+\frac{n}{2} \in \mathbb{Z}}} \abs{\Delta u_{j}^{n-1}}^{3}\right)^{\frac{1}{3}} \left(\sum_{n=1}^{N-1}\sum_{\substack{\abs{j} \leq J \\ j+\frac{n}{2} \in \mathbb{Z}}}1\right)^{\frac{2}{3}} \notag\\
     & \leq  2\norm{\phi}_{C_{0}^{\alpha}} (\Delta x)^{\alpha}\norm{Q_{k}} \norm{k}_{BV}T  + \lambda C^{\frac{1}{3}}\norm{\phi}_{C_{0}^{\alpha}} (\Delta x)^{\alpha+\frac{2}{3}} \norm{Q_{u}}  \left(\frac{4XT}{\lambda (\Delta x)^{2} }\right)^{\frac{2}{3}}. \notag
\end{align}
Now, with the same arguments leading to  \eqref{eq:L1tildebd_c}, we conclude that the set $\{\tilde{\mathcal{L}}^{\Delta}_{3}\}_{\Delta >0}$ is  compact in  $\mathrm{W}^{-1,q_1}(\Omega)$
for the same $q_1 \in (1,2)$ as in the case of $\tilde{\mathcal{L}}^{\Delta}_{1},$ $\tilde{\mathcal{L}}^{\Delta}_{2}$ and $\tilde{\mathcal{L}}^{\Delta}_{4,2}.$ 
\par
By summarizing Steps 1 and 2, we conclude that  the collection $\{\tilde{\mathcal{L}}^{\Delta}\}_{\Delta >0}$ is compact in $\mathrm{W}^{-1,q_1}(\Omega),$ where $q_1 \in (1,2)$ is as given in Step 2.
Furthermore, combining this result with  \eqref{eq:Lhatcompact},  we deduce that
\begin{nalign}
    \{\mathcal{L}^{\Delta}\}_{\Delta > 0}\,  \mbox{is compact in}\,  \mathrm{W}^{-1,q}(\Omega), \,  \mbox{ for } \, q = q_1\in (1,2).
\end{nalign}
Now, owing to the $\mathrm{L}^{\infty}$ boundedness \eqref{Linf_bd} of the approximate solutions $u_\Delta$ (Theorem \ref{thm:maxpri}),  it is straightforward to see that $\{\mathcal{L}^{\Delta}\}$ is bounded in $\mathrm{W}^{-1,r}(\Omega)$ for any $r>2.$ Finally, applying the the interpolation result Lemma \ref{lemma:interpolation} we conclude that $$\{\mathcal{L}^{\Delta}\}_{\Delta > 0} \,  \mbox{is compact in} \, \mathrm{W}^{-1,2}(\Omega).$$
Since $\Omega$ is an arbitrary bounded open subset, it follows that $\{\mathcal{L}^{\Delta}\}_{\Delta >0}$ is compact in $\mathrm{W}^{-1,2}_{\mathrm{loc}}(\mathbb{R}\times\mathbb{R}_{+}).$ This completes the proof.
\end{proof}
\end{lemma}

We now present the weak solution convergence result for the proposed scheme \eqref{eq:NTscheme}, in the following theorem.
\begin{theorem}\label{thm:weakconv}
 Let the initial datum $u_{0} \in (\mathrm{L}^{\infty}\cap \mathrm{BV})(\mathbb{R})$ and $\{u_{\Delta}\}_{\Delta >0}$ be the approximate solutions  \eqref{eq:pcsoln} obtained from the second-order scheme \eqref{eq:NTscheme} under the CFL condition \eqref{eq:cfl_cubicest}. Then, there exists a subsequence $\{\Delta_{m}\}_{m \in \mathbb{N}}$ with $\displaystyle \lim_{m \to \infty} \Delta_{m}= 0$ such that $\{u_{\Delta_{m}}\}_{m \in \mathbb{N}}$ converges strongly  to a weak solution of the problem \eqref{eq:problem}. i.e., 
\begin{align*}
    u_{\Delta_{m}} \xrightarrow[]{\Delta_m \to 0} u \,\, \mbox{in} \,\, \mathrm{L}^{p}_{\mathrm{loc}}(\mathbb{R}\times\mathbb{R^{+}}) \,\, \, \mbox{for any} \,\, p \in  [1,\infty) \,\, \mbox{and a.e. in} \,\,\mathbb{R}\times\mathbb{R}_{+}. 
\end{align*}
\begin{proof}
Consider a sequence  $\displaystyle\{\Delta_{m}\}_{m \in \mathbb{N}}$ such that $\displaystyle \lim_{m \rightarrow \infty} \Delta_{m} = 0.$ Now, Theorem \ref{thm:maxpri} and Lemma \ref{lemma:W-12compactness} allow us to use the compensated compactness result from Theorem \ref{lemma:compcompactness} for the sequence $\{u_{\Delta_{m}}\}_{m \in \mathbb{N}}.$ Hence, there exists a subsequence, again denoted by $\displaystyle\{\Delta_{m}\}_{m \in \mathbb{N}},$ such that $ u_{\Delta_{m}} \rightarrow u$ pointwise a.e., as $m \rightarrow \infty$. Therefore, for any $p \in [1, \infty),$ it follows that $\abs{u_{\Delta_{m}}-u}^{p} \rightarrow 0$ as $\Delta_{m} \rightarrow 0,$ pointwise a.e..  Since the approximate solutions satisfy the $\mathrm{L}^{\infty}$-estimate $\norm{u_{\Delta}} \leq C_{u_{0}}$ by Theorem \ref{thm:maxpri} and the limit $u \in \mathrm{L}^{\infty}(\mathbb{R}\times \mathbb{R_{+}})$ (by Theorem \ref{lemma:compcompactness}), it follows that $\abs{u_{\Delta}-u}^{p} \leq (C_{u_{0}}+\norm{u})^p.$ Consequently,  applying the dominated convergence theorem, for any compact set $K \subseteq \mathbb{R} \times \mathbb{R}_{+},$ we have
    \begin{nalign}\label{eq:strongconv}
        &\lim_{\Delta \rightarrow 0} \norm{u_{\Delta_{m}}-u}^{p}_{\mathrm{L}^{p}(K)} = 0, \, \mbox{or equivalently,}  \\& u_{\Delta_{m}} \xrightarrow{\Delta_{m} \rightarrow 0} u \, \mbox{in} \, \mathrm{L}^{p}_{\mathrm{loc}}(\mathbb{R} \times \mathbb{R}_{+}) \, \mbox{for any} \, p \in [1, +\infty).
    \end{nalign}
Finally, employing a Lax-Wendroff type argument \cite{laxwendroff1960}  which uses the the $\mathrm{L}^{\infty}$- boundedness of $\{u_\Delta\}_{\Delta >0}$ and  the strong convergence $u_{\Delta_{m}} \xrightarrow{} u$ from \eqref{eq:strongconv}, we can show that the limit $u$ is a weak solution of \eqref{eq:problem}. This completes the proof.
\end{proof}
\end{theorem}

\section{Convergence to the entropy solution}\label{section:entropy}
In order to show the entropy convergence, we follow the approach of \cite{gowda2023, vila1988, vila1989}.  As the first step in this framework, we derive an entropy-convergence result for numerical schemes in the predictor-corrector form:
\begin{align*}
u^{n+1}_{\jph}= \bar{u}^{n+1}_{\jph} - a_{\jpo}^{n+1} + a_{j}^{n+1},
\end{align*} approximating \eqref{eq:problem}. This result builds on the idea of Theorem 3.1 in \cite{vila1988}, incorporating essential modifications to address the discontinuous flux case, as detailed in the following theorem.
\begin{theorem}\label{theorem:vilaconv} Suppose that a numerical scheme approximating \eqref{eq:problem} can be written in the form:
\begin{equation}\label{eq:predcorr_vila}
u^{n+1}_{\jph}= \bar{u}^{n+1}_{\jph} - a_{\jpo}^{n+1} + a_{j}^{n+1},
\end{equation}
 where\\
 (i) $\{\bar{u}_{\jph}^{n+1}\}_{\jinz}$ is computed from $\{u_j^n\}_{\jinz}$ using the Lax-Friedrichs scheme \eqref{eq:lxfscheme} as:
 \begin{align}
     \bar{u}^{n+1}_{\jph} &= \half( u^n_{j}+u^n_{j+1}) - \lambda \left(f(k_{j+1}, u^{n}_{j+1})-f(k_{j}, u^{n}_{j})\right), \label{eq:nt_in_predcorrform}
 \end{align}
 (ii) $ |a^{n+1}_{j}| \leq \mathcal{K}{\Delta x}^\alpha,$ $\jinz$ for some constant $\mathcal{K}>0$ which is independent of $\Delta x$ and for some $\alpha \in (\frac{2}{3},1),$ \\
 (iii) The approximate solutions $u_{\Delta}$ obtained from the predictor-corrector scheme is bounded in the $\mathrm{L}^{\infty}$-norm,  and for any fixed $T>0, X>0,$ with $N:=\lfloor T/\Delta t  \rfloor+1$ and $J:=\lfloor X/\Delta x \rfloor +1,$ satisfies an estimate of the form \begin{align} \label{eq:cubicest_vila} \Delta x\sum_{n=0}^{N-1}\sum_{\substack{\abs{j} \leq J \\ j+\frac{n}{2} \in \mathbb{Z}}}\abs{\Delta u_{\jph}^{n}}^3\leq C, \end{align} for a constant $C$ independent of $\Delta x$ and converges pointwise a.e. to a function $u \in \mathrm{L}^{\infty}(\mathbb{R}\times \mathbb{R}_{+}).$ \\
Then the limit $u$ of the approximate solutions $u_{\Delta}$  generated  by the scheme \eqref{eq:predcorr_vila} is the entropy solution \eqref{eq:entropy_ineq} to the problem \eqref{eq:problem}. 

\begin{proof}
The first-order Lax-Friedrichs scheme \eqref{eq:lxfscheme} used to obtain $\bar{u}^{n+1}_{\jph}$ satisfies a discrete cell entropy inequality (see  \cite{karlsen2004} for more details) given by: 
\begin{align*}
    \abs{\bar{u}^{n+1}_{\jph}-c} - \frac{1}{2}\abs{u^{n}_{j+1}- c}&- \frac{1}{2}\abs{u^{n}_{j}- c} + \lambda \left(F(k_{\jpo}, u_{\jpo}^{n}, c) - F(k_{j}, u_{j}^{n}, c)\right) \\& \spc - \lambda \sign(\bar{u}^{n+1}_{\jph}-c)(f(k_{\jpo},c)- f(k_{j},c)) \leq 0,
\end{align*}
for $F(k,u,c):= \mathrm{sgn}(u-c)(f(k,u)-f(k,c)).$ This in turn implies that
\begin{nalign}\label{eq:lxf_entropy_b}
    \abs{\bar{u}^{n+1}_{\jph}-c} - \frac{1}{2}\abs{u^{n}_{j+1}- c}- \frac{1}{2}\abs{u^{n}_{j}- c} &+ \lambda \left(F(k_{\jpo}, u_{\jpo}^{n}, c) - F(k_{j}, u_{j}^{n}, c)\right) \\& \spc - \lambda \abs{f(k_{\jpo},c)- f(k_{j},c)} \leq 0.
\end{nalign}
Now, consider a non-negative test function $\phi$ with  $\mathrm{supp} (\phi) \subseteq [-X, X]\times[0,T],$ let $N=\lfloor T/\Delta t  \rfloor+1,$  $J=\lfloor X/\Delta x \rfloor +1$ and denote $\phi_{\jph}^n:= \phi(x_{\jph},t^{n}).$ Adding $\abs{u^{n+1}_{\jph}-c}$ to either side of \eqref{eq:lxf_entropy_b}, multiplying it with $\Delta x \phi_{\jph}^{n}$ and summing over $n= 0,1, \dots, N-1$ and $j+\frac{n}{2} \in \mathbb{Z}$ with $\abs{j} \leq J,$  we obtain the following inequality
\begin{align}\label{eq:ntentropy_sum}
    &\Delta x \sum_{n=0}^{N-1}\sum_{\substack{\abs{j} \leq J \\ j+\frac{n}{2} \in \mathbb{Z}}}\phi_{\jph}^{n}\abs{u^{n+1}_{\jph}-c}  - \frac{1}{2} \Delta x \sum_{n=0}^{N-1}\sum_{\substack{\abs{j} \leq J \\ j+\frac{n}{2} \in \mathbb{Z}}}\phi_{\jph}^{n}\abs{u^{n}_{j+1}- c} \\& \spc - \frac{1}{2} \Delta x \sum_{n=0}^{N-1}\sum_{\substack{\abs{j} \leq J \\ j+\frac{n}{2} \in \mathbb{Z}}}\phi_{\jph}^{n}\abs{u^{n}_{j}- c} -  \Delta x \sum_{n=0}^{N-1}\sum_{\substack{\abs{j} \leq J \\ j+\frac{n}{2} \in \mathbb{Z}}}\phi_{\jph}^{n}\lambda \abs{f(k_{\jpo},c)- f(k_{j},c)} \notag\\& \spc +  \Delta x \sum_{n=0}^{N-1}\sum_{\substack{\abs{j} \leq J \\ j+\frac{n}{2} \in \mathbb{Z}}}\phi_{\jph}^{n}\lambda \left(F(k_{\jpo}, u_{\jpo}^{n}, c) - F(k_{j}, u_{j}^{n}, c)\right)     \leq \mathcal{J}, \notag
\end{align}
where 
\begin{align*}
     \mathcal{J} &:= \Delta x \sum_{n=0}^{N-1}\sum_{\substack{\abs{j} \leq J \\ j+\frac{n}{2} \in \mathbb{Z}}} \phi_{\jph}^{n} \left(\abs{u^{n+1}_{\jph}-c}-\abs{\bar{u}^{n+1}_{\jph}-c}\right).
\end{align*}
\par
First, we show that the term $\mathcal{J} \rightarrow 0$ as $\Delta \rightarrow 0.$  For this, we consider the function $\displaystyle S(u) =
\begin{cases}
\frac{u^2}{2} + \frac{1}{2}, & \text{if } |u| < 1, \\
|u|, & \text{if } |u| \geq 1
\end{cases}
$ and for a fixed $c \in \mathbb{R}$ and  $\epsilon> 0, $ we define $\displaystyle S_{\epsilon}(u) = \epsilon S(\frac{u-c}{\epsilon}).$ It is easy to see that $S_{\epsilon} \rightarrow \abs{u-c}$ as $\epsilon \rightarrow 0,$ uniformly in the supremum norm. In particular,
\begin{nalign}\label{eq:s_epsilon_bd}
S_{\epsilon}(u)- \abs{u-c}&=0, \quad \,\,\,\,\mbox{if} \,\, \abs{u-c} \geq \epsilon, \,\, \mbox{and}\\
\abs[\big]{S_{\epsilon}(u)-\abs{u-c}} &\leq \frac{3}{2}\epsilon, \quad \mbox{if} \,\, \abs{u-c} < \epsilon.
\end{nalign}
Upon adding and subtracting suitable terms, $\mathcal{J}$ can be written as 
\begin{align}\label{eq:Jassum}
    \mathcal{J}&= \mathcal{J}_{0}+ \mathcal{J}_{1} + \mathcal{J}_{2},
\end{align}
where
\begin{align*}
\mathcal{J}_{0}&:=
\Delta x \sum_{n=0}^{N-1}\sum_{\substack{\abs{j} \leq J \\ j+\frac{n}{2} \in \mathbb{Z}}} \phi_{\jph}^{n}\left(S_{\epsilon}(u^{n+1}_{\jph})- S_{\epsilon}(\bar{u}^{n+1}_{\jph})\right),\\ 
     \mathcal{J}_{1} &:= \Delta x \sum_{n=0}^{N-1}\sum_{\substack{\abs{j} \leq J \\ j+\frac{n}{2} \in \mathbb{Z}}} \phi_{\jph}^{n} \left(\abs{u^{n+1}_{\jph}-c}- S_{\epsilon}(u^{n+1}_{\jph})\right),\\  \mathcal{J}_{2}&:= - \Delta x\sum_{n=0}^{N-1}\sum_{\substack{\abs{j} \leq J \\ j+\frac{n}{2} \in \mathbb{Z}}} \phi_{\jph}^{n}\left(\abs{\bar{u}^{n+1}_{\jph}-c}-S_{\epsilon}(\bar{u}^{n+1}_{\jph})\right).
\end{align*}
Now, \eqref{eq:s_epsilon_bd} implies that for $\epsilon \leq \Delta t^{2},$ 
\begin{align*}
    \abs{\mathcal{J}_{1}}, \abs{\mathcal{J}_{2}} &\leq \frac{3}{2}\epsilon \Delta x \sum_{n=0}^{N-1}\sum_{\substack{\abs{j} \leq J \\ j+\frac{n}{2} \in \mathbb{Z}}} \norm{\phi}     \leq \frac{3}{2}\Delta t \norm{\phi} (8XT),
\end{align*}
and hence
\begin{align}\label{eq:J1J2goes0}
     \lim_{\Delta \rightarrow 0} \mathcal{J}_{1} =\lim_{\Delta \rightarrow 0} \mathcal{J}_{2} = 0.
\end{align}              
Further, the convexity of $S_{\epsilon}$ gives us
\begin{nalign}\label{eq:S_epsi_convex}
    S_{\epsilon}(u^{n+1}_{\jph})-S_{\epsilon}(\bar{u}^{n+1}_{\jph}) \leq S_{\epsilon}^{\prime}(u^{n+1}_{\jph})\left(u^{n+1}_{\jph} - \bar{u}^{n+1}_{\jph}\right) = S_{\epsilon}^{\prime}(u^{n+1}_{\jph})\left(a^{n+1}_{j} - a^{n+1}_{\jpo}\right).
\end{nalign}
Using \eqref{eq:S_epsi_convex} in $\mathcal{J}_0$ and applying summation by parts, we obtain the following inequality
\begin{align}\label{eq:J0_lb_a}
\mathcal{J}_{0} &\leq \Delta x \sum_{n=0}^{N-1}\sum_{\substack{\abs{j} \leq J \\ j+\frac{n}{2} \in \mathbb{Z}}} a^{n+1}_{\jpo}\left(\phi_{j+\frac{3}{2}}^{n} S_{\epsilon}^{\prime}(u^{n+1}_{j+\frac{3}{2}})-\phi_{\jph}^{n} S_{\epsilon}^{\prime}(u^{n+1}_{\jph})\right).
\end{align}
Now, adding and subtracting the term $\Delta x a_{\jpo}^{n+1}\phi_{j+\frac{3}{2}}^{n} S_{\epsilon}^{\prime}(u^{n+1}_{j+\frac{1}{2}})$ inside the summation in the RHS of \eqref{eq:J0_lb_a}, we have $\mathcal{J}_0 \leq \mathcal{J}_{0}^{a} + \mathcal{J}_{0}^{b},$  where
\begin{align*}
    \mathcal{J}_{0}^{a} &:=  \Delta x \sum_{n=0}^{N-1}\sum_{\substack{\abs{j} \leq J \\ j+\frac{n}{2} \in \mathbb{Z}}} a^{n+1}_{\jpo}\phi_{j+\frac{3}{2}}^{n}\left( S_{\epsilon}^{\prime}(u^{n+1}_{j+\frac{3}{2}})- S_{\epsilon}^{\prime}(u^{n+1}_{\jph})\right),\\
     \mathcal{J}_{0}^{b} &:= \Delta x \sum_{n=0}^{N-1}\sum_{\substack{\abs{j} \leq J \\ j+\frac{n}{2} \in \mathbb{Z}}} a^{n+1}_{\jpo}\left(\phi_{j+\frac{3}{2}}^{n}-\phi_{\jph}^{n} \right)S_{\epsilon}^{\prime}(u^{n+1}_{\jph}).
\end{align*}
Using the assumption that $\abs{a^{n+1}_{\jpo} }\leq \mathcal{K} \Delta x ^{\alpha}$ for $\jinz,
$ along with the condition $\mathrm{supp}(\phi) \subseteq [-X, X]\times[0,T],$ we apply Hölder's inequality, followed by \eqref{eq:cubicest_vila}, to obtain the following estimate for $\mathcal{J}_{0}^{a}$: 
\begin{align*}
    \abs{\mathcal{J}_{0}^{a}} &\leq \mathcal{K} (\Delta x)^{\alpha}\norm{\phi}\norm{S_{\epsilon}^{\prime\prime}} \left(\Delta x \sum_{n=0}^{N-1}\sum_{\substack{\abs{j} \leq J \\ j+\frac{n}{2} \in \mathbb{Z}}}\abs{\Delta u^{n+1}_{j+1}}^{3}\right)^{\frac{1}{3}}\left(\Delta x  \sum_{n=0}^{N-1}\sum_{\substack{\abs{j} \leq J \\ j+\frac{n}{2} \in \mathbb{Z}}} 1\right)^{\frac{2}{3}} \\
    & \leq \mathcal{K} (\Delta x)^{\alpha}\norm{\phi}\norm{S_{\epsilon}^{\prime\prime}} C^{\frac{1}{3}}\left(\frac{8XT}{\Delta t}\right)^{\frac{2}{3}} \leq \frac{1}{\lambda^{\frac{2}{3}}} \mathcal{K} (\Delta x)^{\alpha-\frac{2}{3}}\norm{\phi}\norm{S_{\epsilon}^{\prime\prime}} C^{\frac{1}{3}}\left(8XT\right)^{\frac{2}{3}}.
\end{align*}
Now, since $\alpha > \frac{2}{3}$ by assumption, it follows that $\abs{\mathcal{J}_{0}^{a}} \rightarrow 0$ as $\Delta \rightarrow 0.$ 
Next, it is straightforward to see that 
\begin{align*}
    \abs{\mathcal{J}_{0}^{b}} &= \abs[\Big]{\frac{1}{\lambda} \Delta x \Delta t \sum_{n=0}^{N-1}\sum_{\substack{\abs{j} \leq J \\ j+\frac{n}{2} \in \mathbb{Z}}} a^{n+1}_{\jpo}\left(\frac{\phi_{j+\frac{3}{2}}^{n}-\phi_{\jph}^{n}}{\Delta x} \right)S_{\epsilon}^{\prime}(u^{n+1}_{\jph})}\\ & \leq \frac{1}{\lambda} \mathcal{K} (\Delta x)^{\alpha} \norm{\phi_{x}}\norm{S_{\epsilon}^{\prime}}8XT,
\end{align*}
by which we conclude that $\abs{\mathcal{J}_{0}^{b}} \rightarrow 0$ as $\Delta \rightarrow 0.$ Therefore, $\displaystyle \lim_{\Delta \rightarrow 0}\mathcal{J}_{0}
\leq 0.$ This, together with \eqref{eq:J1J2goes0}  applied on \eqref{eq:Jassum} imply that 
\begin{nalign}\label{eq:J1tendsto0}
    \displaystyle \lim_{\Delta \rightarrow 0}\mathcal{J} \leq 0.
\end{nalign}
Finally, using \eqref{eq:J1tendsto0} and proceeding as  in the proofs of Lemmas 5.3 and 5.4 from \cite{karlsen2004}, it is easy to show from \eqref{eq:ntentropy_sum} that
\begin{align}\label{eq:lhs_entropysum}
&\int\int_{\mathbb{R} \times \mathbb{R}_+} (|u - c|\phi_t + \sign(u - c)(f(k, u) - f(k, c))\phi_x) \, \dif x \, \dif t 
\\ & \spc + \int_\mathbb{R} |u_0 - c|\phi(x, 0) \, \dif x  + \int\int_{(\mathbb{R}\setminus D) \times \mathbb{R}_+ } |f(k(x), c)_x|\phi \, \dif x \, \dif t \notag\\ & \spc 
+ \sum_{m=1}^M \int_0^\infty |f(k^+_m, c) - f(k^-_m, c)|\phi(x_m, t) \, \dif t \geq 0, \notag
\end{align}  
where $D$ is as in \ref{hyp:H6}. This concludes the proof.
\end{proof}
\end{theorem}
\begin{remark}  
We emphasize that the entropy convergence result in Theorem 3.1 of \cite{vila1988}  is based on a $\mathrm{BV}$-estimate for the approximate solutions. However, for conservation laws with discontinuous coefficients, $\mathrm{BV}$-estimates need not be available in general (see \cite{adimurthi2011}). In this context, as a key novelty of our approach, we show in Theorem \ref{theorem:vilaconv} that a significantly weaker cubic estimate of the form \eqref{eq:C_cubiclemma} is sufficient to establish the desired entropy convergence. 
\end{remark}
\par
Now, our strategy essentially is to use Theorem \ref{theorem:vilaconv} to prove the convergence of the proposed second-order scheme \eqref{eq:NTscheme}  to the entropy solution. To this end, we note that the scheme \eqref{eq:NTscheme} can be written in the predictor-corrector form 
\begin{align}\label{eq:predcorr_final}
    \bar{u}^{n+1}_{\jph} &= \half( u^n_{j}+u^n_{j+1}) - \lambda \left(f(k_{j+1}, u^{n}_{j+1})-f(k_{j}, u^{n}_{j})\right),\\
    u^{n+1}_{\jph} &= \bar{u}^{n+1}_{\jph} - a_{\jpo}^{n+1} + a_{j}^{n+1}, \notag
\end{align}
with the correction terms 
\begin{nalign}\label{eq:nt_cor_step}
a_{j}^{n+1} &:= \lambda\left(f(k_{j}, u_{j}^{n+\frac{1}{2}})- f(k_{j}, u_{j}^{n})\right) + \frac{1}{8}\sigma_{j}^n,\\ & = \left(-\frac{\lambda^{2}}{2} f_{u}(k_{j}, \zeta_{j})f_{u}(k_{j}, u_{j}^{n}) + \frac{1}{8}\right)\sigma_{j}^n,
\end{nalign} for some $\zeta_{j} \in \mathcal{I}(u_{j}^{n}, u_{j}^{\nph}).$ Next, we modify the scheme \eqref{eq:predcorr_final} by redefining the slopes \eqref{eq:slopes} as  
\begin{nalign}\label{eq:modifiedslopes}
    &\sigma_{j}^n \\&= \mbox{ minmod}\left(  (u^{n}_{j+1}-u^{n}_{j}), \half( u^{n}_{j+1}-u^{n}_{j-1}), ( u^{n}_{j}- u^{n}_{j-1}), \sign(u^{n}_{\jpo}-u^{n}_j)\tilde{\mathcal{K}}(\Delta x)^{\alpha}\right),
\end{nalign}  for some constants $\tilde{\mathcal{K}}>0$ and $\alpha \in (\frac{2}{3},1).$
 With this modification, it is direct to see from \eqref{eq:nt_cor_step} that the correction terms $\{a_{j}^{n+1}\}_{\jinz}$ in \eqref{eq:predcorr_final} satisfy the estimate 
    \begin{nalign}\label{eq:corr_term_bound}
        \abs{a_j^{n+1}} \leq \mathcal{K} (\Delta x)^{\alpha},
    \end{nalign}
for $\jinz,$ where $\displaystyle \mathcal{K}:= \left(\frac{1}{2}\lambda^{2}\norm{f_{u}}^{2} + \frac{1}{8}\right)\tilde{\mathcal{K}}.$
Finally, we conclude this section with the entropy convergence result, stated below. 
\begin{theorem}\label{theorem:ntconvergence} Let the initial datum $u_{0}$ be such that $u_{0} \in (\mathrm{L}^{\infty}\cap \mathrm{BV})(\mathbb{R}),$ with $\barbelow{u} \leq u_{0}(x) \leq \bar{u},$  for $x \in \mathbb{R}.$  Under the CFL condition \eqref{eq:cfl_cubicest} and hypotheses \ref{hyp:H1}-\ref{hyp:H7}, the approximate solutions $\{u_{\Delta}\}_{\Delta >0}$ \eqref{eq:pcsoln} obtained from the scheme \eqref{eq:NTscheme} with the modified slopes \eqref{eq:modifiedslopes} converge to the entropy solution \eqref{eq:entropy_ineq} of the problem \eqref{eq:problem}.
    \begin{proof}
        Since the second-order scheme \eqref{eq:NTscheme} with the modified slopes \eqref{eq:modifiedslopes} can be reformulated in the predictor-corrector form, it suffices to show that the hypotheses (i)-(iii) of Theorem \ref{theorem:vilaconv} hold true.  From \eqref{eq:predcorr_final}, it is clear that the predictor step employs the Lax-Friedrichs time-stepping, thereby confirming the validity of condition (i) in Theorem \ref{theorem:vilaconv}. From \eqref{eq:corr_term_bound}, it follows that with the slope modification \eqref{eq:modifiedslopes}, the hypothesis (ii) of Theorem \ref{theorem:vilaconv} also holds. Finally, observing that the $\mathrm{L}^{\infty}$- estimate \eqref{Linf_bd} (Theorem \ref{thm:maxpri}), the cubic estimate \eqref{eq:C_cubiclemma} (Lemma \ref{lemma:cubicest}), the $\mathrm{W}^{-1,2}_{\mathrm{loc}}$ compactness result (Lemma \ref{lemma:W-12compactness}), and the convergence theorem (Theorem \ref{thm:weakconv}) remain valid with the modified slopes \eqref{eq:modifiedslopes}, it follows that hypothesis (iii) is also satisfied. Now, an application of Theorem \ref{theorem:vilaconv} yields the desired convergence to the entropy solution.
    \end{proof}
\end{theorem}

\begin{remark} In the implementation of the numerical scheme, the slope modification \eqref{eq:modifiedslopes} is not really required.  This is because, for any given mesh-size $\Delta x,$ we can choose $\tilde{K}>0$ large enough so that the modified slope \eqref{eq:modifiedslopes} reduces to \eqref{eq:slopes}. In particular, for mesh sizes $\Delta x \geq \epsilon$ for some fixed $\epsilon > 0,$ we can choose $\tilde{K}= 2C_{u_0}\epsilon^{-\alpha},$ where $C_{u_0}$ is as given in \eqref{Linf_bd}.  For more details, see \cite{leroux1981} (p. 158, below Eq. (26)), \cite{vila1988} (p. 68, below Fig. 3) and \cite{viallon1991} (p. 577, Remark).

\end{remark}

\section{Numerical experiments}\label{section:numerical}
In this section, we present numerical experiments to illustrate the performance of the second-order scheme \eqref{eq:NTscheme} in comparison to the first-order Lax-Friedrichs scheme \eqref{eq:lxfscheme}. In the following, we denote the second-order scheme \eqref{eq:NTscheme} by SO and the Lax-Friedrichs scheme \eqref{eq:lxfscheme} by LF. Alluding to Remark \ref{remark:concaveflux}, we focus on examples where the flux function is strictly concave, as these are the types of fluxes most commonly considered in the literature. We divide the computational domain into cells of size $\Delta x.$ The CFL condition \eqref{eq:cfl}, required for proving the convergence of the SO scheme, is highly restrictive and need not be optimal. Therefore, to compute the time-step $\Delta t,$ we impose the less restrictive CFL condition \eqref{eq:cfl_maxpri} (which ensures the maximum principle). In each test case, both the LF and SO solutions are computed with the time-step corresponding to the SO scheme. Also, we apply absorbing boundary conditions in both the examples.\\
\begin{example}\label{example:multipl_flux}
We consider an example from \cite{karlsen2004b}, where the flux function in  \eqref{eq:problem} is given by:
\begin{nalign}\label{eq:multi_flux_eg}
 f(k,u)&=ku(1-u), \quad 
 k(x)&=\begin{cases}
    3 \,\,\, \mbox{if} \, x < 0,\\
    1 \,\,\,  \mbox{if} \, x \geq 0.
\end{cases}
\end{nalign}
 The discontinuous flux \eqref{eq:multi_flux_eg} is depicted in Figure \ref{fig:fluxes}(a), where $f_{l}(u)=3u(1-u)$ and $f_{r}(u)=u(1-u).$
\begin{figure}
\centering
     \begin{subfigure}[b]{0.475\textwidth}
         \centering
         \includegraphics[width=\textwidth]{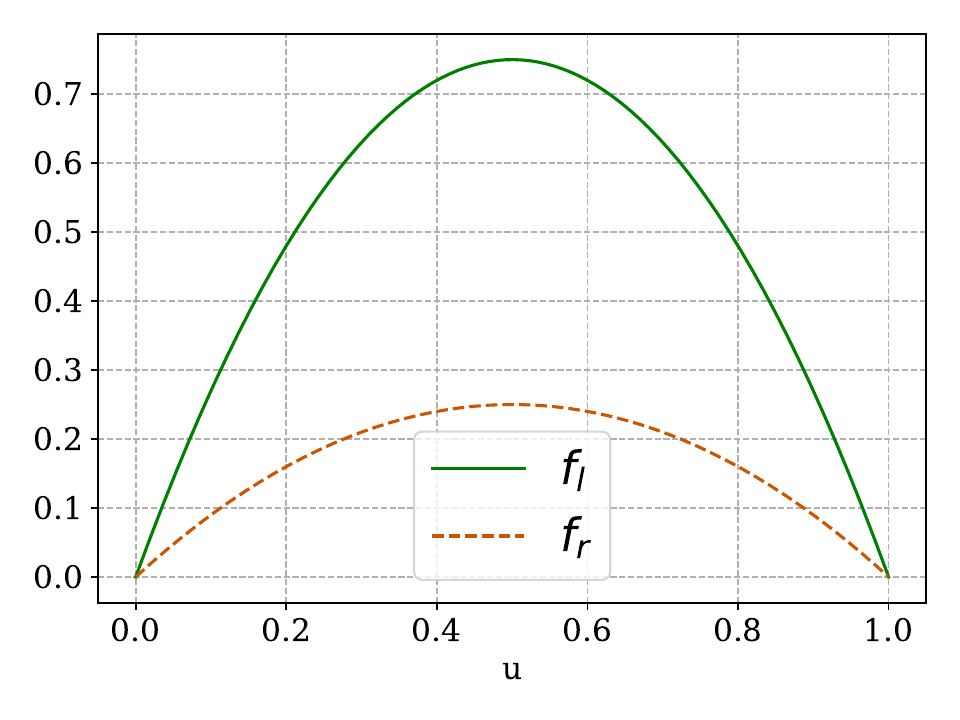}
         \caption{Example \ref{example:multipl_flux}}
     \end{subfigure}\hspace{0.45cm}
     \begin{subfigure}[b]{0.475\textwidth}
         \centering
         \includegraphics[width=\textwidth]{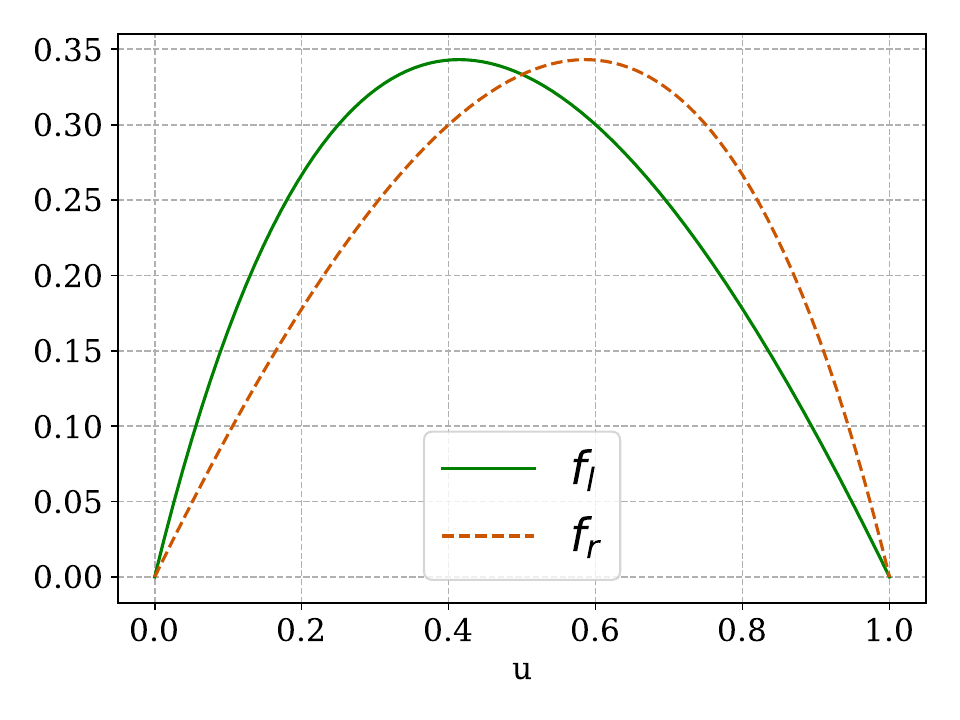}
         \caption{Example \ref{example:mishra2017}}
              \end{subfigure}
     \caption{The fluxes to the left ($f_l$) and right ($f_r$) of $x=0.$}
     \label{fig:fluxes}
\end{figure}
 Note that there is no flux crossing in this case and the crossing condition \ref{hyp:H7} holds trivially for \eqref{eq:multi_flux_eg}.  We set the initial condition as the constant function 
 \begin{align}\label{eq:karlsen_ic}
 u_0(x)= 0.15,
 \end{align}
 and compute the numerical solutions in the domain $[-1,1]$ with a mesh of size $\Delta x = 2/50$ at the time levels $t\in \{0.8, 1.6\}.$ Here, the reference solutions are computed with the LF scheme using a fine mesh of size $\Delta x = 2/1000.$ The results are displayed in Figure \ref{fig:karlsen2004b}. We observe that the SO scheme exhibits lower numerical diffusion and provides a more accurate approximation of the solution, especially near the discontinuities. 
 
\begin{figure}
\centering
     \begin{subfigure}[b]{0.475\textwidth}
         \centering
         \includegraphics[width=\textwidth]{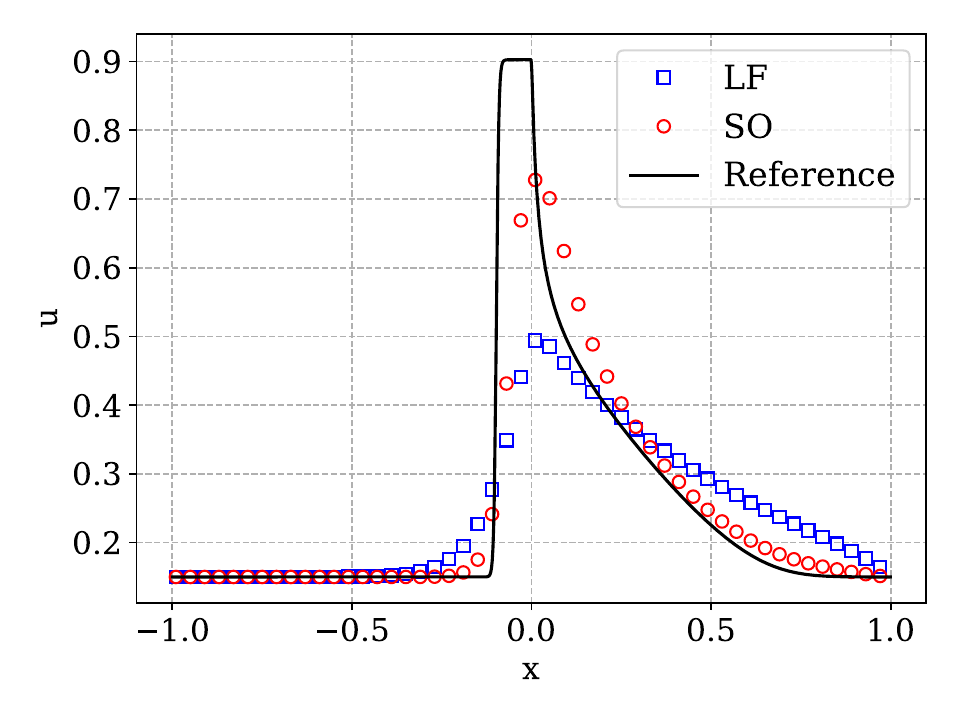}
         \caption{t = 0.8}
     \end{subfigure}\hspace{0.45cm}
     \begin{subfigure}[b]{0.475\textwidth}
         \centering
         \includegraphics[width=\textwidth]{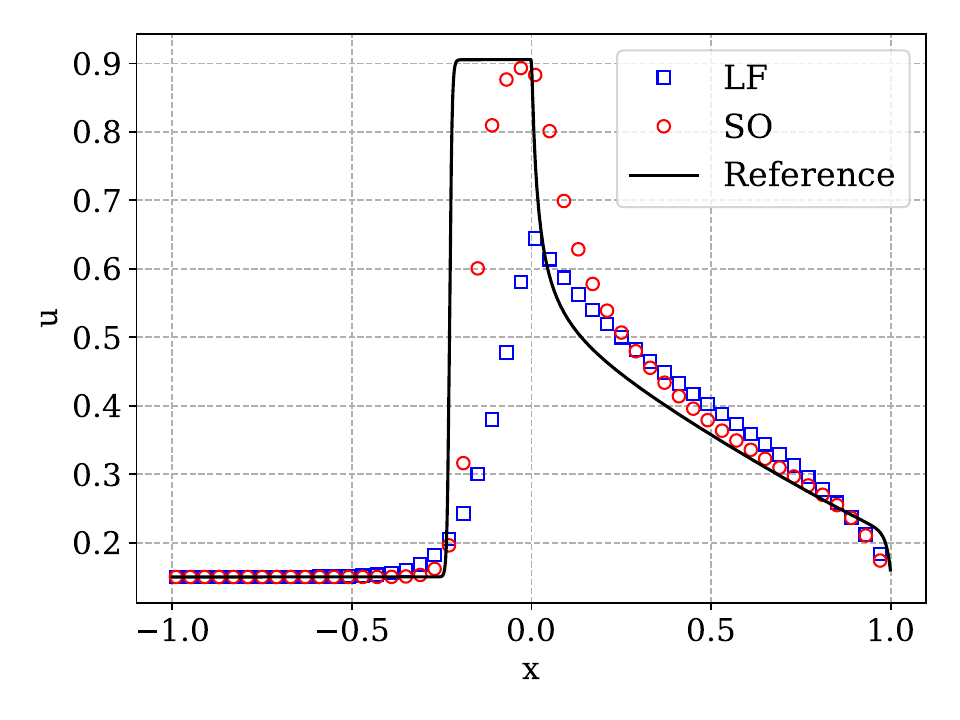}
         \caption{t = 1.6}
              \end{subfigure}
     \caption{Example \ref{example:multipl_flux}. Numerical solutions obtained by evolving \eqref{eq:karlsen_ic} with $\Delta x = 0.04$ and $\Delta t = 1/750.$}
     \label{fig:karlsen2004b}
\end{figure}
\end{example}

\begin{example}\label{example:mishra2017} Next, we consider an example studied in \cite{mishra2017rev}, where the flux function in \eqref{eq:problem} is of the form 
\begin{align}\label{eq:bellflux}
    f(H(x),u(x,t)) &= H(x)f_{r}(u) + (1-H(x))f_{l}(u) = \begin{cases}
        f_{l}(u) \quad \mbox{for} \, x < 0,\\
        f_{r}(u) \quad \mbox{for} \, x \geq 0,
    \end{cases}
\end{align}
with
\begin{align}
    f_{l}(u) := \frac{2u(1-u)}{1+u}, \quad f_{r}(u) := \frac{2u(1-u)}{2-u}.
    \label{eq:dflux}
\end{align}
The discontinuous flux \eqref{eq:bellflux} is depicted in Figure \ref{fig:fluxes}(b), where we note that the crossing condition \ref{hyp:H7} is satisfied by \eqref{eq:bellflux}. We set the initial datum to be
\begin{align}
u(x,0) = \begin{cases}
    0.9 & \mbox{ if } x \le 0,\\
    0.2 & \mbox { otherwise.} 
\end{cases}
\label{eq:icdflu}
\end{align}
and compute the solutions using $\Delta x =8/50$ at the time levels $t=1.0$ and $t=2.0,$ which are compared in Figures \ref{fig:mishra2017}(a) and \ref{fig:mishra2017}(b), respectively. The reference solutions are obtained with the LF scheme on a fine mesh of size $\Delta x =8/2000.$ As shown in Figure \ref{fig:mishra2017}, the proposed SO scheme yields more accurate solutions, in comparison to the LF scheme.
\begin{figure}
\centering
     \begin{subfigure}[b]{0.47\textwidth}
         \centering
         \includegraphics[width=\textwidth]{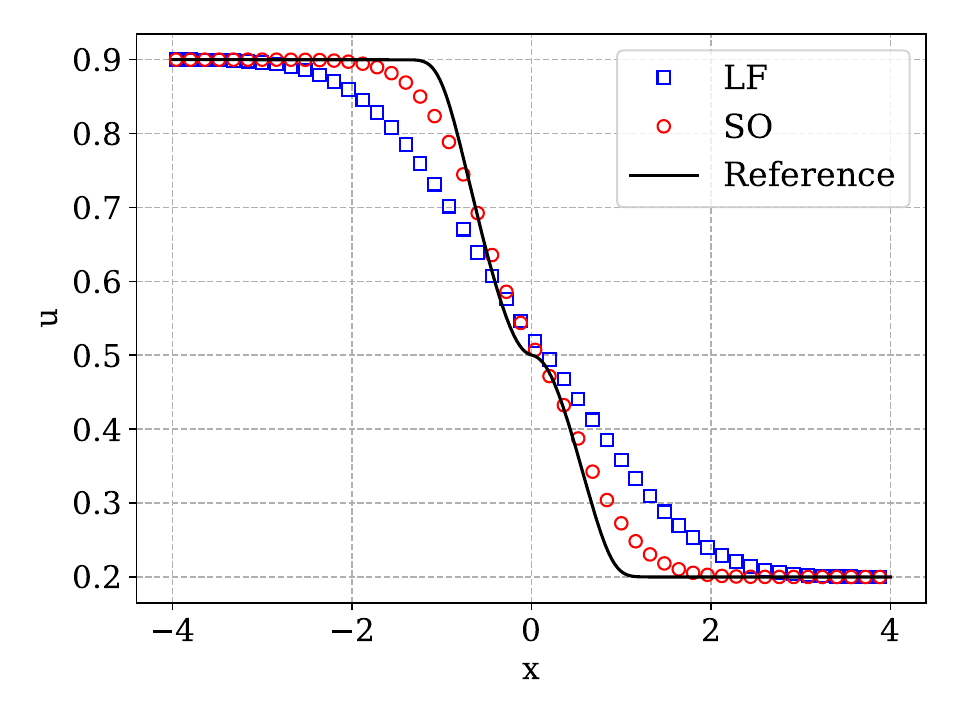}
         \caption{t = 1.0}
     \end{subfigure}\hspace{0.45cm}
     \begin{subfigure}[b]{0.47\textwidth}
         \centering
         \includegraphics[width=\textwidth]{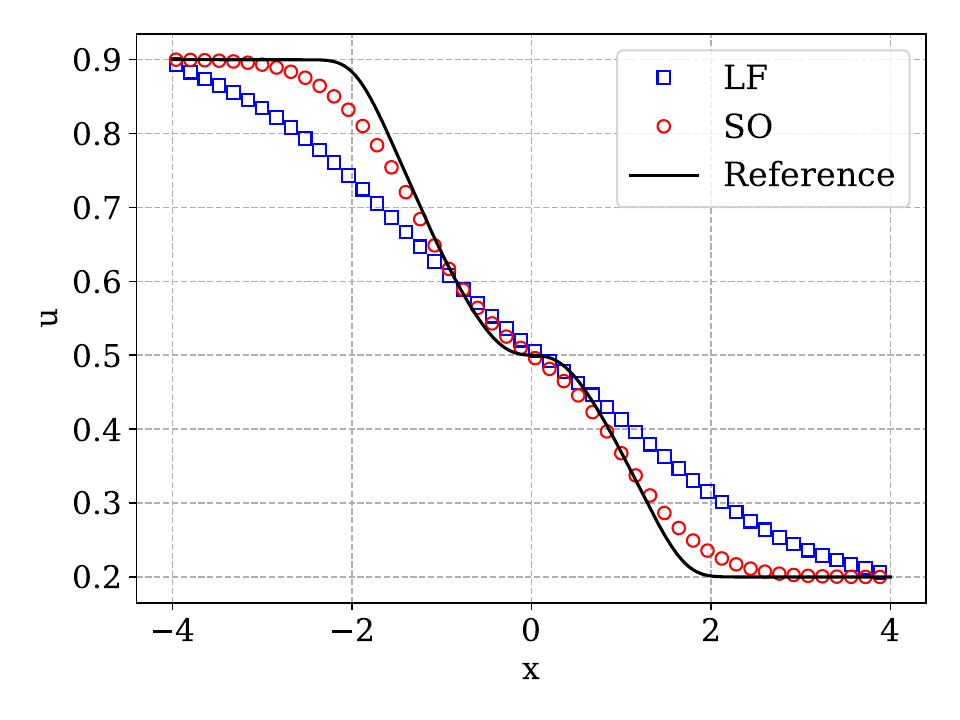}
         \caption{t = 2.0}
              \end{subfigure}
     \caption{Example \ref{example:mishra2017}. Numerical solutions obtained by evolving \eqref{eq:icdflu} with $\Delta x = 0.16$ and $\Delta t = 0.008.$} 
     \label{fig:mishra2017}
\end{figure}
\end{example}

\section{Conclusion} \label{sec:conclusion}
In this work, we have proposed and analyzed a second-order numerical scheme based on MUSCL-type reconstruction for conservation laws with discontinuous flux functions, where the flux may exhibit multiple discontinuities. The scheme features a simple formulation inspired by the Nessyahu-Tadmor central scheme. We employ the framework of compensated compactness to establish the convergence of the proposed scheme. This analysis requires several key estimates, which we carefully derive by exploiting the structural properties of the numerical method.

The uniqueness of solutions to conservation laws with discontinuous flux is a topic of considerable interest, and various approaches exist in the literature. In this work, we focus on a Kruzkov-type entropy condition equivalent to the vanishing viscosity solution approach. We establish the convergence of entropy solutions by incorporating a mesh-dependent parameter into the scheme. To the best of our knowledge, this is the first comprehensive convergence result, including entropy solutions, for a MUSCL-type reconstruction scheme applied to conservation laws with discontinuous flux.
The numerical solutions presented in Section~\ref{section:numerical}, where we compare the proposed second-order scheme and the first-order Lax-Friedrichs (LF) scheme, clearly indicate the superiority of the SO scheme over its FO counterpart.

This work opens up several avenues for further research. For instance, investigating the convergence of the MUSCL-Hancock \cite{vanleer1984}     scheme and Runge-Kutta MUSCL \cite{gottlieb1998otalVD} schemes, as well as the convergence of schemes using Godunov-type numerical interface fluxes. Additionally, it would be interesting to explore convergence to the (A, B)-entropy solution for arbitrary admissible (A, B)-connections. Extending the convergence analysis of MUSCL-type schemes to multidimensional cases is also a compelling area for future research. Furthermore, we note that 
the convergence analysis presented in this work relies on certain a priori estimates (Lemmas \ref{lemma:osle}, \ref{lemma:cubicest}, \ref{lemma:quadest}), which require the flux function to be strictly convex (or strictly concave). Extending this analysis to encompass more general flux functions remains an open question for future research. We aim to address these topics in forthcoming papers.

\section*{Acknowledgments} Nikhil Manoj acknowledges financial support in the form of a doctoral fellowship from the Council of Scientific and Industrial Research (CSIR), Government of India.

\appendix
\section{Proof of Lemma \ref{lemma:osle}}\label{app:A}The corresponding  estimate for the case when $k\equiv 1,$ was derived in  \cite{popov2006b}. We follow a similar approach for the problem \eqref{eq:problem}, where the coefficient $k$ can be discontinuous.
\subsection{A bound on non-decreasing sequences} In this section, 
we establish a result  regarding the jumps generated by applying the scheme \eqref{eq:NTscheme} to non-decreasing sequences with finitely many non-zero jumps.  The main result is stated in the following lemma, and its proof is provided at the end of this section.

\begin{lemma}\label{lemma:mainlemma_inseq}
 Consider a non-decreasing sequence $\{u_{j}\}_{\jinz}$ such that for some $l,r \in \mathbb{Z},$ 
  the jumps $\delta_{\jph}:= u_{j+1}- u_{j} = 0$ for all $j \notin \{l-1,l, \dots, r-1\}.$ Further, let $\abs{u_j} \leq \norm{C_{u_{0}}}, \, \forall \jinz,$ where $C_{u_{0}}$ is as in \eqref{Linf_bd}.
Let $\{u^{\prime}_{\jph}\}_{\jinz}$ be obtained from $\{u_{j}\}_{\jinz}$ by applying the time-update formula \eqref{eq:NTscheme} and
     denote the corresponding jumps $\delta_{j}^{\prime}:= u_{\jph}^{\prime}- u_{\jmh}^{\prime}$ for $\jinz$. Under the CFL condition \eqref{eq:cfl} there exists a constant $\Theta \geq 0$ independent of $\Delta x$ such that
    the jump sequence $\{\delta_{j}^{\prime}\}_{\jinz}$ satisfies the following estimate
    \begin{nalign}\label{eq:incr_seq_result}
    \sum_{j=l-1}^{r+1}\left((\delta_{\jmh})^2- (\delta_{j-1}^{\prime})^2\right)& \geq \frac{1}{500} \lambda \gamma_{1}\sum_{j=l-1}^{r+1} (\delta_{\jmh})^{3} +  \frac{1}{6400}\sum_{j=l-1}^{r+1}\left(\Delta^{2}\delta_{\jmtbt}\right)^{2}\\ & \spc - \Theta\sum_{j=l-2}^{r+2}\abs{\Delta k_{\jmtbt}},
\end{nalign}
where $\Delta^{2}\delta_{\jmtbt}:= \Delta\delta_{j-1}- \Delta\delta_{j-2}.$ 
\end{lemma}
To prove Lemma \ref{lemma:mainlemma_inseq}, we need to establish several auxiliary lemmas, starting with the following one.
\begin{lemma}\label{lemma:inseqjumpbd}
 Let $\{u_{j}\}_{\jinz}$ and  $\{\delta^\prime_j\}_{\jinz}$  be as in Lemma \ref{lemma:mainlemma_inseq}. Under the CFL condition \eqref{eq:cfl}
    the jumps $\{\delta_{j}^{\prime}\}_{\jinz}$ satisfy the following estimate
\begin{nalign}
    \abs{\delta_{j}^\prime} &\leq  (\delta_{\jph}+ \delta_{\jmh}) + \lambda \norm{f_{k}}\left(\abs{\Delta k_{\jph}}+\abs{\Delta k_{\jmh}}\right).
\end{nalign}
\end{lemma}
\begin{proof}
From the scheme \eqref{eq:NTscheme}, the jump $\delta_{j}^{\prime}$ can be written as 
\begin{nalign}\label{eq:deltaprime}
    \delta_{j}^{\prime}
    & = \frac{1}{2}(\delta_{\jph}+ \delta_{\jmh}) -\lambda\left( f(k_{j+1}, u_{j+1}^{\half})- 2 f(k_{j}, u_{j}^{\half})+ f(k_{j-1}, u_{j-1}^{\half}) \right)\\ & \spc - \frac{1}{8}(s_{j+1}- 2s_{j}+ s_{j-1}),
\end{nalign}
where we define $u_{j}^{\half}:=u_j- \frac{\lambda}{2} f_{u}(k_{j}, u_j)s_j,$ $s_{j}:=\min\{\delta_{\jmh},\delta_{\jph}\}, \jinz.$ We can also write $$u_{j+1}^{\half}- u_{j}^{\half}= \delta_{\jph} - \frac{\lambda}{2}(a_{j+1} s_{j+1} - a_{j} s_{j}),$$
where, \begin{align}\label{eq:ajn_defn}
    a_{j} := f_{u}(k_{j}, u_{j}), \, \mbox{for} \, j \in \mathbb{Z}.
\end{align} Now, by adding and subtracting the term $f(k_{j}, u_{j+1}^{\half})$ and using the mean value theorem, we can write
\begin{align*}
    f(k_{j+1}, u_{j+1}^{\half})-  f(k_{j}, u_{j}^{\half}) &= f_{k}(c_{12}, u_{j+1}^{\half})\Delta k_{\jph} + f_{u}(k_{j},\zeta_{12})(u_{j+1}^{\half}-u_{j}^{\half}) \\
    &= \bar{b}_{\jph}\Delta k_{\jph} + \bar{a}_{\jph}\left(\delta_{\jph} - \frac{\lambda}{2}(a_{j+1} s_{j+1} - a_{j} s_{j})\right),
\end{align*}
where we define 
\begin{nalign}\label{eq:abar_bbar}
    \bar{a}_{\jph} &:=f_{u}(k_{j},\zeta), \quad    \bar{b}_{\jph}:= f_{k}(c, u_{j+1}^{\half}),
\end{nalign} for some $\zeta \in \mathcal{I}(u_{j}^{\half}, u_{j+1}^{\half})$ and $c \in \mathcal{I}(k_{j}, k_{j+1}).$
Now, \eqref{eq:deltaprime} reduces to
\begin{nalign}\label{eq:delta_pr_b}
    \delta_{j}^{\prime}
    & = \frac{1}{2}(\delta_{\jph}+ \delta_{\jmh}) - \frac{1}{8}(\Delta s_{\jph}- \Delta s_{\jmh} ) -\lambda \left[\bar{b}_{\jph}\Delta k_{\jph}- \bar{b}_{\jmh}\Delta k_{\jmh}\right] \\ & \hspace{0.5cm} -\lambda\left[\bar{a}_{\jph}\left(\delta_{\jph} - \frac{\lambda}{2}(a_{j+1} s_{j+1} - a_{j} s_{j})\right)\right.\\& \spc \spc \left. -\bar{a}_{\jmh}\left(\delta_{\jmh} - \frac{\lambda}{2}(a_{j} s_{j} - a_{j-1} s_{j-1})\right) \right],
\end{nalign}
The CFL condition \eqref{eq:cfl} implies that $\kappa \leq \frac{1}{4},$ which yields the following upper bound 
\begin{nalign}\label{eq:deltpr_ub}
\delta_{j}^\prime &\leq \left(\frac{1}{2}+\frac{1}{8}+\kappa+\kappa^2\right)(\delta_{\jph}+ \delta_{\jmh}) + \lambda \norm{f_{k}}\left(\abs{\Delta k_{\jph}}+\abs{\Delta k_{\jmh}}\right)\\
&\leq (\delta_{\jph}+ \delta_{\jmh}) + \lambda \norm{f_{k}}\left(\abs{\Delta k_{\jph}}+\abs{\Delta k_{\jmh}}\right).
\end{nalign}
\par
Analogously, we obtain a lower bound 
\begin{nalign}\label{eq:deltapr_lb}
\delta_{j}^\prime &\geq \left(\frac{1}{2}-\frac{1}{8}-\kappa-\kappa^2\right)(\delta_{\jph}+ \delta_{\jmh}) - \lambda \norm{f_{k}}\left(\abs{\Delta k_{\jph}}+\abs{\Delta k_{\jmh}}\right)\\
& \geq - \lambda \norm{f_{k}}\left(\abs{\Delta k_{\jph}}+\abs{\Delta k_{\jmh}}\right).
\end{nalign}
\par
The estimates \eqref{eq:deltpr_ub} and \eqref{eq:deltapr_lb} together imply that
\begin{nalign}\label{eq:deltaprime_bound}
    \abs{\delta_{j}^\prime} &\leq  (\delta_{\jph}+ \delta_{\jmh}) + \lambda \norm{f_{k}}\left(\abs{\Delta k_{\jph}}+\abs{\Delta k_{\jmh}}\right).
\end{nalign} This completes the proof.
\end{proof}
\begin{lemma}\label{lemma:deltamin_deltapr_est}
    Let $\{u_{j}\}_{\jinz}$ and  $\{\delta^\prime_j\}_{\jinz}$  be as in Lemma \ref{lemma:mainlemma_inseq}, with $\abs{u_{j}}\leq C_{u_{0}},$ for all $\jinz,$ (where $C_{u_{0}}$ is as in \eqref{Linf_bd}). Define $\delta^{\prime\prime}_j$ by replacing the  term $(a_{j+1}^n s_{j+1}^{n} - a_{j}^n s_{j}^{n})$ with $\bar{a}_{\jph}^{n}\Delta s_{\jph}$  in $\delta_{j}^{\prime}$ \eqref{eq:delta_pr_b}, where $\bar{a}_{\jph}^{n}$ is as defined in \eqref{eq:abar_bbar}.
Then \begin{nalign}\label{eq:djpr_minus_djprpr}
       \abs[\Big]{\sum_{j=l-2}^{r}(\delta_{j}^{\prime})^2- \sum_{j=l-2}^{r}(\delta_{j}^{\prime\prime})^2} &\leq  55\lambda^2  \gamma_{2}\norm{f_{u}}\sum_{j=l-2}^{r}(\delta_{\jph}^{n})^{3}+ P_{1},
\end{nalign}
where $P_{1} := (16\lambda^{2} C_{u_{0}}^{2}\norm{f_{u}}\norm{f_{uk}}+ 44\lambda^{3}\norm{f_k}\gamma_{2}\norm{f_{u}}C_{u_{0}}^{2} + 16\lambda^{3}\norm{f_k}C_{u_{0}}\norm{k})\\\sum_{j=l-2}^{r+1}\abs{\Delta k_{\jmh}}.$  
\end{lemma}
\begin{proof}
The jump term $\delta_{j}^{\prime \prime}$ is given by 
\begin{align}\label{eq:deltadobpr}
    \delta_{j}^{\prime \prime} &=\frac{1}{2}(\delta_{\jph}+ \delta_{\jmh}) - \frac{1}{8}(\Delta s_{\jph}- \Delta s_{\jmh} ) -\lambda \left(\bar{b}_{\jph}\Delta k_{\jph}- \bar{b}_{\jmh}\Delta k_{\jmh}\right)\\ & \hspace{0.5cm} -\lambda\left(\bar{a}_{\jph}\left(\delta_{\jph} - \frac{\lambda}{2}a_{\jph} \Delta s_{\jph} \right) -\bar{a}_{\jmh}\left(\delta_{\jmh} - \frac{\lambda}{2}\bar{a}_{\jmh} \Delta s_{\jmh}\right) \right).\notag
\end{align} With  arguments analogous to those used in \eqref{eq:deltpr_ub}, it follows that
\begin{nalign}\label{eq:deltdoubprime_bound}
    \abs{\delta_{j}^{\prime\prime}} &\leq  (\delta_{\jph}+ \delta_{\jmh}) + \lambda \norm{f_{k}}\left(\abs{\Delta k_{\jph}}+\abs{\Delta k_{\jmh}}\right).
\end{nalign}
Now, we write the difference $\delta_{j}^{\prime}- \delta_{j}^{\prime\prime}$ as
\begin{nalign}\label{eq:deltapr_dif_deldobpr}
    \delta_{j}^{\prime}- \delta_{j}^{\prime\prime} 
    & = \frac{\lambda^2}{2}(\tau_{\jph}-\tau_{
    \jmh}),
\end{nalign}
where $\tau_{\jph} := \bar{a}_{\jph}\left(a_{j+1} s_{j+1} - a_{j} s_{j}- \bar{a}_{\jph}\Delta s_{\jph}\right),$ which can be bounded as \begin{nalign}\label{eq:tau_bd}
    \abs{\tau_{\jph}} &\leq \norm{f_{u}}\abs{ s_{j+1}(a_{j+1}-\bar{a}_{\jph}) -  s_{j}(a_{j}-\bar{a}_{\jph})}.
\end{nalign}
Further, adding and subtracting the term $f_{u}(k_{j},u_{j+1})$ in $a_{j+1}-\bar{a}_{\jph}$ yields
\begin{nalign}\label{eq:a_min_abar}
    a_{j+1}-\bar{a}_{\jph} &= f_{uk}(c,u_{j+1}) \Delta k_{\jph} + f_{uu}(k_{j},\zeta^{\prime}) \left(u_{j+1}-\zeta\right),
\end{nalign}
where $\zeta$ is as in \eqref{eq:abar_bbar}, $c \in \mathcal{I}(k_{j},k_{j+1})$ and  $\zeta^{\prime} \in \mathcal{I}(\zeta,u_{\jpo}).$
Now, due to the CFL condition \eqref{eq:cfl}, it follows that
\begin{nalign}\label{eq:halfvaluebds}
u_{j+1}^{\half} &= u_{j+1} - \frac{1}{2}\lambda a_{\jpo}s_{j+1} \leq u_{j+1} + \frac{\delta_{\jph}}{4} =: u_{j+1}^{+} \quad \mbox{and}\\
u_{j+1}^{\half} &= u_{j+1} - \frac{1}{2}\lambda a_{j+1}s_{j+1} \geq u_{j+1} - \frac{\delta_{\jph}}{4} =: u_{j+1}^{-}, \quad \mbox{for} \, \jinz. 
\end{nalign}
 The CFL condition \eqref{eq:cfl} and \eqref{eq:halfvaluebds} ensure that $u_{j}^{\half} \in [u_{j}^{-}, u_{j+1}]$ and $u_{\jpo}^{\half} \in [u_{j}, u_{j+1}^{+}].$ Combining this with the fact that $\zeta \in \mathcal{I}(u_{j}^{\half}, u_{j+1}^{\half}),$ we derive the following estimate:
\begin{nalign}\label{eq:ujpo-zeta_bd}\abs{u_{j+1}-\zeta} \leq \max\{u_{j+1}- u_{j}^{-} ,  u_{j+1}^{+} - u_{j+1}\} 
    & \leq \max\left\{\delta_{\jph}+ \frac{\delta_{\jmh}}{4} ,  \frac{\delta_{\jph}}{4}\right\} \\ & \leq \delta_{\jph}+ \frac{\delta_{\jmh}}{4},
\end{nalign}
which, when combined with \eqref{eq:a_min_abar} and hypothesis \ref{hyp:H2}, provides the following upper bound
\begin{nalign}\label{eq:ajpominajph}
     \abs{a_{j+1}-\bar{a}_{\jph}} \leq \norm{f_{uk}} \abs{\Delta k_{\jph}} + \gamma_{2} \left(\delta_{\jph} + \frac{\delta_{\jmh}}{4}\right).
\end{nalign}
In a similar way, we have 
\begin{nalign}\label{eq:ajminajph}
    \abs{a_{j}-\bar{a}_{\jph}} &=  \abs{f_{u}(k_{j},u_{j}) - f_{u}(k_{j},\zeta)}= \abs{f_{uu}(k_{j},\bar{\zeta}) \left(u_{j}-\zeta\right)} \leq  \gamma_{2}\abs{u_{j}-\zeta},
\end{nalign}
where $\bar{\zeta} \in \mathcal{I}(\zeta, u_{j}).$
Here, analogous to \eqref{eq:ujpo-zeta_bd}, $\abs{u_{j}-\zeta}$ can be bounded as
\begin{nalign}\label{eq:ujminzeta}
    \abs{u_{j}-\zeta} \leq \max\{u_{j}- u_{j}^{-},  u_{j+1}^{+} - u_{j} \}& \leq \max\left\{\frac{\delta_{\jmh}}{4} ,  \delta_{\jph}+ \frac{\delta_{\jph}}{4} \right\} \\ & \leq \frac{5}{4} \delta_{\jph} + \frac{1}{4}  \delta_{\jmh}.
\end{nalign}
Using \eqref{eq:ajpominajph}, \eqref{eq:ajminajph} and \eqref{eq:ujminzeta} in \eqref{eq:tau_bd}, we obtain an estimate on the term $\tau_{\jph}$ in \eqref{eq:deltapr_dif_deldobpr}:
\begin{nalign}\label{eq:tau_bd_final}
    \abs{\tau_{\jph}} 
    &\leq \gamma_{2}\norm{f_{u}}\delta_{\jph} \left(\frac{9}{4}\delta_{\jph} + \frac{1}{2}\delta_{\jmh}\right) + \norm{f_{u}}\norm{f_{uk}}\delta_{\jph}\abs{\Delta k_{\jph}}.
\end{nalign}
The estimate \eqref{eq:tau_bd_final} applied on \eqref{eq:deltapr_dif_deldobpr} yields
\begin{nalign}\label{eq:delta_pr_min_dou_pr_bd}
\abs{\delta_{j}^{\prime}- \delta_{j}^{\prime\prime}}  
& \leq \frac{\lambda^2}{2}  \gamma_{2}\norm{f_{u}}\left(\frac{9}{4}+\frac{1}{2}+\frac{9}{4}+\frac{1}{2}\right)
M_{\delta_{\jmh}}\\
& \hspace{0.5cm} + \frac{\lambda^2}{2}  \norm{f_{u}}\norm{f_{uk}}\left(\delta_{\jph}\abs{\Delta k_{\jph}}+ \delta_{\jmh}\abs{\Delta k_{\jmh}}\right),
\end{nalign}
where $M_{\delta_{\jmh}}:=\max\left\{ (\delta_{\jph})^{2}, (\delta_{\jmh})^{2}, (\delta_{j-\frac{3}{2}})^{2}\right\}.$
Making use of the estimate \eqref{eq:delta_pr_min_dou_pr_bd} together with  \eqref{eq:deltaprime_bound} and \eqref{eq:deltdoubprime_bound} and using the inequality $ab^{2} \leq a^{3}+b^{3},$ for $a,b \geq 0,$  it follows that
\begin{nalign}
       \abs[\Big]{\sum_{j=l-2}^{r}(\delta_{j}^{\prime})^2- \sum_{j=l-2}^{r}(\delta_{j}^{\prime\prime})^2} &\leq \sum_{j=l-2}^{r}\abs{\delta_{j}^{\prime}- \delta_{j}^{\prime\prime}}\abs{\delta_{j}^{\prime}+ \delta_{j}^{\prime\prime}}\\&\leq  55\lambda^2  \gamma_{2}\norm{f_{u}}\sum_{j=l-2}^{r}(\delta_{\jph})^{3}+ P_{1}.
\end{nalign}
To obtain $P_1$, we have used the bound $\abs{\delta_{\jph}} \leq 2C_{u_0}.$ This completes the proof. \end{proof}

Now we define the term $\displaystyle\mathcal{D} := \sum_{j=l-1}^{r+1}\left((\delta_{\jmh})^2 - (\delta_{j-1}^{\prime\prime})^2\right)$ and write this in a suitable form to obtain a lower bound. We begin with the notations 
\begin{nalign}\label{eq:alpha_phi_def}
    \alpha_{\jmh} :=  \frac{1}{2}+\lambda\bar{a}_{\jmh} \,\, \mbox{and} \,\, \varphi_{\jmh}:= \alpha_{\jmh}(1-\alpha_{\jmh}),
\end{nalign}
and reformulate  the modified jumps \eqref{eq:deltadobpr}  as
\begin{nalign}\label{eq: deltadobpr_mod}
    \delta_{j-1}^{\prime\prime} &= (1-\alpha_{\jmh})\delta_{\jmh} + \alpha_{\jmtbt}\delta_{j-\frac{3}{2}} -\frac{1}{2}\left(\varphi_{\jmh}\Delta s_{\jmh}-\varphi_{\jmtbt}\Delta s_{\jmtbt}\right) \\
    &\spc-\lambda \left(\bar{b}_{\jmh}\Delta k_{\jmh}- \bar{b}_{\jmtbt}\Delta k_{\jmtbt}\right).
\end{nalign}
Plugging in the expression \eqref{eq: deltadobpr_mod}, we can write  $\mathcal{D}$ as 
\begin{nalign}\label{eq:D_a}
\mathcal{D} = \mathcal{I}_{1} + \mathcal{I}_{2} + \mathcal{I}_{3} + P_{2},  
\end{nalign}
where 
\begin{align}\label{eq:I1I2I3P2}
    \mathcal{I}_{1} & := \sum_{j=l-1}^{r+1}\delta_{\jmh}^2 - \sum_{j=l-1}^{r+1}\left(\alpha_{\jmtbt}\delta_{\jmtbt}+ (1-\alpha_{\jmh})\delta_{\jmh}\right)^{2},\\
    \mathcal{I}_{2} &:= \sum_{j=l-1}^{r+1}\left((1-\alpha_{\jmh})\delta_{\jmh} + \alpha_{\jmtbt}\delta_{j-\frac{3}{2}}\right)\left(\varphi_{\jmh}\Delta s_{\jmh}-\varphi_{\jmtbt}\Delta s_{\jmtbt}\right),\notag\\
    \mathcal{I}_{3} &:= -\frac{1}{4}\sum_{j=l-1}^{r+1}\left(\varphi_{\jmh}\Delta s_{\jmh}-\varphi_{\jmtbt}\Delta s_{\jmtbt}\right)^{2},\notag\\
    P_2 &:= -\lambda^{2}\sum_{j=l-1}^{r+1} \left(\bar{b}_{\jmh}\Delta k_{\jmh}- \bar{b}_{\jmtbt}\Delta k_{\jmtbt}\right)^{2} \notag\\& \spc +2\lambda\sum_{j=l-1}^{r+1}\left(\bar{b}_{\jmh}\Delta k_{\jmh}- \bar{b}_{\jmtbt}\Delta k_{\jmtbt}\right)\left(\alpha_{\jmtbt}\delta_{\jmtbt}+ (1-\alpha_{\jmh})\delta_{\jmh}\right)\notag\\
    & \spc -\lambda\sum_{j=l-1}^{r+1}\left(\bar{b}_{\jmh}\Delta k_{\jmh}- \bar{b}_{\jmtbt}\Delta k_{\jmtbt}\right)\left(\varphi_{\jmh}\Delta s_{\jmh}-\varphi_{\jmtbt}\Delta s_{\jmtbt}\right). \notag
\end{align}
Noting that $\abs{\alpha_{\jmh}}\leq 1$ and $\abs{\varphi_{\jmh}}\leq \frac{1}{2},$ for $\jinz$ (by the CFL condition \eqref{eq:cfl}), we obtain
\begin{align}\label{eq:P2bound}
    \abs{P_{2}} \leq \lambda\norm{f_{k}}(8\lambda\norm{f_{k}}\norm{k}+24C_{u_{0}}). 
\end{align}

Upon a change of index in the summation and grouping the terms appropriately, the term $\mathcal{I}_{1}$ can be rewritten as follows:
\begin{align}\label{eq:I1transform}
\mathcal{I}_{1}
& = \sum_{j=l-1}^{r+1}\left(\delta_{\jmh}^{2}\left(1-\alpha_{\jmh}^{2}-(1-\alpha_{\jmh})^{2}\right)-2\alpha_{\jmtbt} (1-\alpha_{\jmh})\delta_{\jmtbt}\delta_{\jmh}\right) \\
& = \sum_{j=l-1}^{r+1} \left(\alpha_{\jmh}(1-\alpha_{\jmh})\delta_{\jmh}^{2}+\alpha_{\jmtbt}(1-\alpha_{\jmtbt})\delta_{\jmtbt}^{2}\right. \notag\\&\left. \spc -2\alpha_{\jmtbt}(1-\alpha_{\jmh})\delta_{\jmtbt}\delta_{\jmh}\right) \notag\\
& = \sum_{j=l-1}^{r+1} \alpha_{\jmtbt}(1-\alpha_{\jmh})(\delta_{\jmh}-\delta_{\jmtbt})^{2} \notag\\ & \spc + \sum_{j=l-1}^{r+1}(\alpha_{\jmh}-\alpha_{\jmtbt})\left(\alpha_{\jmtbt}\delta_{\jmtbt}^{2}+(1-\alpha_{\jmh})\delta_{\jmh}^{2}\right).\notag
\end{align}
Now, we write
\begin{align}\label{eq:I3mod}
    \mathcal{I}_{3} &= -\frac{1}{4}\sum_{j=l-1}^{r+1}\left(\varphi_{\jmtbt}\Delta^{2}s_{j-1}\right)^{2}+ \mathcal{E}_{1}, 
\end{align}
where $\mathcal{E}_{1} := \mathcal{I}_{3} +\frac{1}{4}\sum_{j=l-1}^{r+1}\left(\varphi_{\jmtbt}\Delta^{2}s_{j-1}\right)^{2}.$
Incorporating the expressions \eqref{eq:I1transform} and \eqref{eq:I3mod} into \eqref{eq:D_a} and by adding and subtracting $\frac{1}{2}\sum_{j=l-1}^{r+1}(\varphi_{\jmtbt}\Delta^{2}\delta_{\jmtbt})$  the term $\mathcal{D}$ is reformulated as
\begin{nalign}\label{eq:D}
\mathcal{D} = \mathcal{Q}_{1} + \mathcal{Q}_{2} + \mathcal{Q}_{3}+ P_{2} + \mathcal{E}_{1}, 
\end{nalign}
where 
\begin{align}\label{eq:Q1Q3Q3}
    \mathcal{Q}_{1} &:= \sum_{j=l-1}^{r+1}\alpha_{\jmtbt}(1-\alpha_{\jmh})(\Delta \delta_{\jmh})^{2} + \mathcal{I}_{2} - \frac{1}{2}\sum_{j=l-1}^{r+1}(\varphi_{\jmtbt}\Delta^{2}\delta_{\jmtbt}),\\
    \mathcal{Q}_{2} &:= \frac{1}{4}\left(2\sum_{j=l-1}^{r+1}\left(\varphi_{\jmtbt}\Delta^{2}\delta_{\jmtbt}\right)^{2}-\sum_{j=l-1}^{r+1}\left(\varphi_{\jmtbt}\Delta^{2}s_{j-1}\right)^{2} \right) \quad \mbox{and} \notag\\
    \mathcal{Q}_{3} &:= \sum_{j=l-1}^{r+1} \Delta \alpha_{j-1}\left(\alpha_{\jmtbt}(\delta_{\jmtbt})^{2} + (1-\alpha_{\jmh})(\delta_{\jmh})^{2}\right). \notag
 \end{align}
 \begin{remark}\label{remark: E1bound} 
     It is possible to derive a bound for the term $\mathcal{E}_{1}$ in  \eqref{eq:D}. To show this,  we begin by rewriting it as
\begin{nalign}\label{eq:E1_exp}
    \mathcal{E}_{1} &= -\frac{1}{4}\sum_{j=l-1}^{r+1} \left((\Delta s_{\jmh})^{2}\left(\varphi_{\jmh}-\varphi_{\jmtbt}\right)\left(\varphi_{\jmh}+\varphi_{\jmtbt}\right) \right)\\ & \spc -\frac{1}{4}\sum_{j=l-1}^{r+1} \left(2 \left(\varphi_{\jmtbt}-\varphi_{\jmh}\right)\varphi_{\jmtbt} \Delta s_{\jmh} \Delta s_{\jmtbt} \right).
\end{nalign}
Immediately, we have the following bound:
\begin{nalign}\label{eq:phidiff}
\abs{\varphi_{\jmh}-\varphi_{\jmtbt}} & = \lambda^{2}\abs{(\bar{a}_{\jmh})^{2}- (\bar{a}_{\jmtbt})^{2} } \leq \lambda^{2}\abs{\bar{a}_{\jmh}- \bar{a}_{\jmtbt}}\abs{\bar{a}_{\jmh}+ \bar{a}_{\jmtbt}}.
\end{nalign}
Next, using the definition of $\bar{a}_{\jmh}$ from \eqref{eq:abar_bbar} and adding and subtracting the term $f_{u}(k_{j-2},\bar{u}_{\jmh}),$ we get
 \begin{nalign}\label{eq:abar_dif}
    \bar{a}_{\jmh}- \bar{a}_{\jmtbt} 
    & = f_{uk}(\bar{k}_{\jmtbt}, \bar{u}_{\jmh})\Delta k_{\jmtbt} + f_{uu}(k_{j-2}, \bar{\bar{u}}_{j-1})(\bar{u}_{\jmh}-\bar{u}_{\jmtbt}),
\end{nalign}
for some $\bar{k}_{\jmtbt} \in \mathcal{I}(k_{j-2}, k_{j-1}),
$  $\bar{u}_{\jmh} \in \mathcal{I}(u_{j-1}^{\half}, u_{j}^{\half}),
$ $\bar{u}_{\jmtbt} \in \mathcal{I}(u_{j-2}^{\half}, u_{j-1}^{\half})
$ and $\bar{\bar{u}}_{j-1} \in \mathcal{I}(\bar{u}_{\jmtbt},\bar{u}_{\jmh}).$
Further, by the CFL restriction \eqref{eq:cfl} and using arguments similar to those in \eqref{eq:halfvaluebds}, we have $\displaystyle \bar{u}_{\jmh} \in [u_{j-1}-\frac{\delta_{\jmtbt}}{4},u_{j}+\frac{\delta_{\jmh}}{4}]$ and $\displaystyle \bar{u}_{\jmtbt} \in [u_{j-2}-\frac{\delta_{\jmtbt}}{4},u_{j-1}+\frac{\delta_{\jmtbt}}{4}].$ This in turn, implies that 
\begin{align*}
    \bar{u}_{\jmh}-\bar{u}_{\jmtbt} &\leq u_{j} + \frac{\delta_{\jmh}}{4} - u_{j-2} + \frac{\delta_{\jmtbt}}{4}\\&= (u_{j}-u_{j-1}) +(u_{j-1}- u_{j-2}) + \frac{1}{4}(\delta_{\jmh}+\delta_{\jmtbt}) \leq \frac{5}{4}(
    \delta_{\jmh}+\delta_{\jmtbt}), \\
    \bar{u}_{\jmh}-\bar{u}_{\jmtbt}&\geq -\frac{\delta_{\jmtbt}}{2},
\end{align*}
and hence \begin{nalign}\label{eq:ubardiff_bd}
    \abs{\bar{u}_{\jmh}-\bar{u}_{\jmtbt}} \leq \frac{5}{4}(
    \delta_{\jmh}+\delta_{\jmtbt}),
\end{nalign} due to the fact that $\delta_{\jmh} \geq 0, \jinz.$
The estimates \eqref{eq:phidiff} and \eqref{eq:ubardiff_bd} together with the hypothesis \ref{hyp:H2} imply that
\begin{nalign}\label{eq:phi_difbd}
\abs{\varphi_{\jmh}-\varphi_{\jmtbt}} 
& \leq  2 \lambda^{2}\norm{f_{u}} \gamma_{2} \frac{5}{4}(
    \delta_{\jmh}+\delta_{\jmtbt}) +  2 \lambda^{2}\norm{f_{u}}\norm{f_{uk}}\abs{\Delta k_{\jmtbt}}.
\end{nalign}
By the CFL condition \eqref{eq:cfl}, we have $\abs{\varphi_{\jmh}}\leq \frac{1}{2}, \jinz.$ Moreover, $\abs{\Delta s_{\jmh}}\leq 2\delta_{\jmh}, \jinz$ and $ab^{2} \leq a^{3}+b^{3},$ for any $a,b \geq 0.$ Invoking these observations and the estimate \eqref{eq:phi_difbd} into the expression \eqref{eq:E1_exp}, we obtain the desired bound
\begin{nalign}\label{eq:E1bd}
    \abs{\mathcal{E}_{1}} &\leq \frac{35}{2} \lambda^{2}\norm{f_{u}}\gamma_{2}\sum_{j=l-1}^{r+1}(\delta_{\jmh})^{3} + 24\lambda^{2}\norm{f_{uk}}C_{u_{0}}^{2}\sum_{j=l-1}^{r+1}\abs{\Delta k_{\jmtbt}}.
\end{nalign}
 \end{remark}
\begin{lemma}\label{lemma:Q3lb}
 Let $\{u_{j}\}_{\jinz}$ be as in Lemma \ref{lemma:mainlemma_inseq}. Then for $\alpha_{\jmh}$ in \eqref{eq:alpha_phi_def}, 
\begin{nalign}
    \Delta \alpha_{j-1}
    & \geq \frac{3}{8}\lambda \gamma_{1} \left(\delta_{\jmh}^{n}+\delta_{\jmtbt}^{n}\right) + \lambda \Delta k_{\jmtbt}f_{ku}(\bar{\bar{k}}_{\jmtbt}, \bar{u}_{\jmh}),
\end{nalign} for $\bar{\bar{k}}_{\jmtbt} \in \mathcal{I}(k_{j-1}, k_{j-2}).$
\end{lemma}

\begin{proof}
Using the fact that $\displaystyle \bar{a}_{\jmh}=\frac{f(k_{j-1},u_{j}^{\half})- f(k_{j-1},u_{j-1}^{\half})}{u_{j}^{\half}-u_{j-1}^{\half}} \, \mbox{for} \, \jinz,$ the term $\Delta \alpha_{j-1}=\lambda(\bar{a}_{\jmh}-\bar{a}_{\jmtbt})$ can be written as follows
\begin{align}\label{eq:alphadif}
&\Delta \alpha_{j-1} \\
&= \lambda \left(\frac{f(k_{j-2},u_{j}^{\half})- f(k_{j-2},u_{j-1}^{\half})}{u_{j}^{\half}-u_{j-1}^{\half}}- \frac{f(k_{j-2},u_{j-1}^{\half})- f(k_{j-2},u_{j-2}^{\half})}{u_{j-1}^{\half}-u_{j-2}^{\half}}\right) \notag\\ 
& \spc + \lambda \left( \frac{\left(f(k_{j-1},u_{j}^{\half})- f(k_{j-2},u_{j}^{\half})\right)-\left(f(k_{j-1},u_{j-1}^{\half})- f(k_{j-2},u_{j-1}^{\half})\right)}{u_{j}^{\half}-u_{j-1}^{\half}}\right) \notag\\
&= \lambda \bar{f}_{j-2}\left[u_{j-2}^{\half},u_{j-1}^{\half},u_{j}^{\half}\right]\left(u_{j}^{\half}-u_{j-2}^{\half}\right) \notag\\ & \spc + \lambda \Delta k_{\jmtbt}\left(\frac{f_{k}(\bar{k}_{\jmtbt},u_{j}^{\half})-f_{k}(\bar{\bar{k}}_{\jmtbt},u_{j-1}^{\half})}{u_{j}^{\half}-u_{j-1}^{\half}}\right) \notag\\
&= \frac{\lambda}{2} f_{uu}(k_{j-2},\bar{\zeta})\left(u_{j}^{\half}-u_{j-2}^{\half}\right) + \lambda \Delta k_{\jmtbt}f_{ku}(\bar{\bar{k}}_{\jmtbt}, \bar{u}_{j-\half}),\notag
\end{align}
where $f_{j-2}(\cdot)$ denotes the denote the divided difference of the function $f(k_{j-2}, \cdot),$  $\bar{\zeta} \in \left(\min\{u_{j-2}^{\half},u_{j-1}^{\half},u_{j}^{\half}\}, \max\{u_{j-2}^{\half},u_{j-1}^{\half},u_{j}^{\half}\}\right),$   $\bar{u}_{\jmh} \in \mathcal{I}(u_{j-1}^{\half}, u_{j}^{\half}),$  and $\bar{k}_{\jmtbt}, \bar{\bar{k}}_{\jmtbt} \in \mathcal{I}(k_{j-1}, k_{j-2}).$ The second term in the last step of \eqref{eq:alphadif} is derived by applying assumption \ref{hyp:H3}, $f_{kk}=0,$ and adding and subtracting $f_{k}(\bar{\bar{k}}_{\jmtbt},u_{j}^{\half}).$ 
\par
Now, the CFL condition \eqref{eq:cfl} allows us to write \begin{nalign}\label{eq:ujminujmtwo}
    u_{j}^{\half} &= u_{j}^{n} - \frac{\lambda}{2}f_{u}(k_{j}, u_{j}^{n})s_{j}^{n} \geq u_{j}^{n} - \frac{\delta_{\jmh}^{n}}{4} \quad \mbox{and}\\
    u_{j-2}^{\half} &= u_{j-2}^{n} - \frac{\lambda}{2}f_{u}(k_{j-2}, u_{j-2}^{n})s_{j-2}^{n} \leq u_{j-2}^{n} + \frac{\delta_{\jmtbt}^{n}}{4}, 
\end{nalign}
and hence $u_{j}^{\half}-u_{j-2}^{\half} \geq \frac{3}{4}\left(\delta_{\jmh}^{n}+\delta_{\jmtbt}^{n}\right).$ Using this bound in \eqref{eq:alphadif} together with \ref{hyp:H2}, yields
\begin{nalign}\label{eq:alphadif_lb}
    \Delta \alpha_{j-1}
    & \geq \frac{3}{8}\lambda \gamma_{1} \left(\delta_{\jmh}^{n}+\delta_{\jmtbt}^{n}\right) + \lambda \Delta k_{\jmtbt}f_{ku}(\bar{\bar{k}}_{\jmtbt}, \bar{u}_{\jmh}).
\end{nalign}
This completes the proof.
\end{proof}

\begin{lemma}\label{lemma:Q1rewrite}
    Let $\{u_{j}\}_{\jinz}$ be as in Lemma \ref{lemma:mainlemma_inseq}. The term $\mathcal{Q}_{1}$ in \eqref{eq:Q1Q3Q3} can be represented as follows
\begin{nalign}
    \mathcal{Q}_{1} = \mathcal{R}_{1} + \mathcal{Q}_{1}^{*} -\mathcal{E}_{4} -\mathcal{E}_{5} -\mathcal{E}_{6},
\end{nalign}
where \begin{align}\label{eq:R1}
    \mathcal{R}_{1} &:= \sum_{j=l-1}^{r+1}\alpha_{\jmtbt}(1-\alpha_{\jmh})(\Delta \delta_{j})^{2} + \sum_{j=l-1}^{r+1} \varphi_{\jmh}\Delta s_{\jmh} \left( \delta_{\jmh}(\Delta \alpha_{j}- \Delta \alpha_{j-1})\right) \\
    &\spc - \sum_{\Delta \delta_{j-1}\leq 0} \varphi_{\jmh}(1-\alpha_{\jmh})(\Delta \delta_{j-1})^{2} \notag\\ & \spc - \frac{1}{2}\sum_{\Delta \delta_{j-1}\geq 0}\varphi_{\jmh}\left((1-\alpha_{\jph})+(1-\alpha_{\jmh})\right)(\Delta \delta_{j-1})^{2}\notag\\
    &\spc -\sum_{\Delta \delta_{j-1}\geq 0} \varphi_{\jmh}\alpha_{\jmtbt}(\Delta \delta_{j-1})^{2} - \frac{1}{2}\sum_{\Delta \delta_{j-1}\leq 0}\varphi_{\jmh}(\alpha_{\jmtbt}+\alpha_{j-\frac{5}{2}})(\Delta \delta_{j-1})^{2},\notag\end{align}
\begin{align}
     \mathcal{Q}_1^{*} &:= \sum_{\Delta \delta_{j-1}\geq 0, \Delta \delta_{j} < 0 }\frac{1}{8}(\Delta \delta_{j-1})(\Delta \delta_{j})- \frac{1}{2}\sum_{\Delta \delta_{j-1}\geq 0, 
\Delta \delta_{j}\leq 0}\frac{1}{8}(\Delta \delta_{j-1})^{2}\\
& \spc -  \frac{1}{2} \sum_{\Delta \delta_{j-1}\geq 0, \Delta 
\delta_{j-2}\leq 0}\frac{1}{8}(\Delta \delta_{j-1})^{2}- \frac{1}{2}\sum_{\Delta \delta_{j-2} \geq 0, \Delta \delta_{j-1} \geq 0} \frac{1}{8}(\Delta^{2}\delta_{\jmtbt})^{2}\notag\\
& \spc +\sum_{\Delta \delta_{j-1}\geq 0, \Delta \delta_{j} < 0 }\frac{1}{8}(\Delta \delta_{j-1})(\Delta \delta_{j}) - \frac{1}{2}\sum_{\Delta \delta_{j-1}\leq 0, 
\Delta \delta_{j}\geq 0}\frac{1}{8}(\Delta \delta_{j-1})^{2}\notag\\
& \spc -  \frac{1}{2} \sum_{\Delta \delta_{j-1}\leq 0, \Delta 
\delta_{j-2}\geq 0}\frac{1}{8}(\Delta \delta_{j-1})^{2} -\frac{1}{2}\sum_{\Delta \delta_{j-2} \leq 0, \Delta \delta_{j-1} \leq 0} \frac{1}{8}(\Delta^{2}\delta_{\jmtbt})^{2}\notag
\\ & \spc  -\frac{1}{32}\sum_{j=l-1}^{r+1}\left(\Delta^{2}\delta_{\jmtbt}\right)^{2}, \notag
\end{align}
and $\displaystyle \mathcal{E}_{4}, \mathcal{E}_{5}, \, \mathcal{E}_{6}$ satisfy the following bounds
\begin{align}\label{eq:E4bd}
    \abs{\mathcal{E}_{4}} &\leq \norm{f_{u}}\lambda^{2} \gamma_{2} \frac{205}{4}\sum_{j=l-1}^{r+1}(
\delta_{\jmh})^{3}
 + 96 C_{u_{0}}^{2} \norm{f_{u}}\lambda^{2}\norm{f_{uk}}\sum_{j=l-1}^{r+1}\abs{\Delta k_{\jmtbt}} \\ & \spc  + \frac{1}{2}\lambda\norm{f_{u}}\sum_{j=l-1}^{r+1}\left(\Delta^{2}\delta_{\jmtbt}\right)^{2}, \notag
\end{align}
\begin{align}\label{eq:E5bd}    \abs{\mathcal{E}_{5}} 
    & \leq  \norm{f_{u}}\lambda^{2} \gamma_{2} \frac{205}{4}\sum_{j=l-1}^{r+1}(
\delta_{\jmh})^{3}
+  96 C_{u_{0}}^{2} \norm{f_{u}}\lambda^{2}\norm{f_{uk}}\sum_{j=l-1}^{r+1}\abs{\Delta k_{\jmtbt}}  \\ & \spc + \frac{5}{8}\lambda\norm{f_{u}}\sum_{j=l-1}^{r+1}\left(\Delta^{2}\delta_{\jmtbt}\right)^{2}, \notag
\end{align}
and
\begin{nalign}
\abs{\mathcal{E}_{6}} \leq \frac{1}{16}  \lambda \norm{f_{u}}\sum_{j=l-1}^{r+1}\left(\Delta^{2}\delta_{\jmtbt}\right)^{2}.
\end{nalign}

\end{lemma}

\begin{proof}
 Applying summation by parts, the term $\mathcal{I}_{2}$ in \eqref{eq:D_a} can be expressed as
\begin{align*}
    \mathcal{I}_{2} 
    &= -\sum_{j=l-1}^{r+1} \varphi_{\jmh}\Delta s_{\jmh} \left(\alpha_{\jmtbt}\Delta \delta_{j-1}+ (1-\alpha_{\jph})\Delta\delta_{j} + \delta_{\jmh}(\Delta \alpha_{j-1}- \Delta \alpha_{j})\right),
\end{align*}
using which the term $\mathcal{Q}_{1}$ in \eqref{eq:Q1Q3Q3} can be reformulated as
\begin{nalign}\label{eq:Q1_mod}
    \mathcal{Q}_1 &= \sum_{j=l-1}^{r+1}\alpha_{\jmtbt}(1-\alpha_{\jmh})(\Delta \delta_{\jmh})^{2}  -\tilde{\mathcal{I}}_{2}-\tilde{\tilde{\mathcal{I}}}_{2} \\ & \spc + \sum_{j=l-1}^{r+1} \varphi_{\jmh}\Delta s_{\jmh} \delta_{\jmh}(\Delta \alpha_{j}- \Delta \alpha_{j-1}) -\frac{1}{2}\sum_{j=l-1}^{r+1}\left(\varphi_{\jmtbt}\Delta^{2}\delta_{\jmtbt}\right)^{2},
\end{nalign}
where 
\begin{align*}
    \tilde{\mathcal{I}}_{2} &:= \sum_{j=l-1}^{r+1}\varphi_{\jmh}(1-\alpha_{\jph})\Delta\delta_{j}\Delta s_{\jmh} \quad \mbox{and} \quad
    \tilde{\tilde{\mathcal{I}}}_{2} := \sum_{j=l-1}^{r+1}\varphi_{\jmh}\alpha_{\jmtbt}\Delta \delta_{j-1}\Delta s_{\jmh}.
\end{align*}
Clearly, $\Delta s_{\jmh}:= (\Delta \delta_{j-1})_{+}+ (\Delta \delta_{j})_{-}.$ Utilizing this relation and following arguments similar to those in equations (50)-(53) of \cite{popov2006a}, we expand the term $\tilde{\mathcal{I}}_{2}$ as 
\begin{align}\label{eq:A}
\tilde{\mathcal{I}}_{2} &=  \mathcal{E}_{2} + \sum_{\Delta \delta_{j-1}\leq 0} \varphi_{\jmh}(1-\alpha_{\jmh})(\Delta \delta_{j-1})^{2}  \notag \\ & \spc +\frac{1}{2}\sum_{\Delta \delta_{j-1}\geq 0}\varphi_{\jmh}\left((1-\alpha_{\jph})+ (1-\alpha_{\jmh})\right)(\Delta \delta_{j-1})^{2}\notag\\
&\spc +\sum_{\Delta \delta_{j-1}\geq 0, \Delta \delta_{j} < 0 }\varphi_{\jmh}(1-\alpha_{\jph})(\Delta \delta_{j-1})(\Delta \delta_{j})\notag\\ & \spc  - \frac{1}{2}\sum_{\Delta \delta_{j-1}\geq 0, 
\Delta \delta_{j}\leq 0}\varphi_{\jmh}(1-\alpha_{\jph})(\Delta \delta_{j-1})^{2}\notag\\
& \spc -  \frac{1}{2} \sum_{\Delta \delta_{j-1}\geq 0, \Delta 
\delta_{j-2}\leq 0}\varphi_{\jmh}(1-\alpha_{\jmh})(\Delta \delta_{j-1})^{2}\notag\\ & \spc -\frac{1}{2}\sum_{\Delta \delta_{j-2} \geq 0, \Delta \delta_{j-1} \geq 0} \varphi_{\jmh}(1-\alpha_{\jmh})(\Delta^{2}\delta_{\jmtbt})^{2},\notag
\end{align}
and 
\begin{align*}
    \mathcal{E}_{2}&:= \sum_{\Delta\delta_{j-1} \leq 0}(\varphi_{\jmtbt}-\varphi_{\jmh})(1-\alpha_{\jmh})(\Delta \delta_{j-1})^{2} \\&\spc+ \frac{1}{2}\sum_{\Delta \delta_{j-2}\geq 0, \Delta \delta_{j-1}\geq 0} (\varphi_{\jmtbt}-\varphi_{\jmh})(1-\alpha_{\jmh})(\Delta \delta_{j-1})^{2}\\
    & \spc -\frac{1}{2}\sum_{\Delta\delta_{j-1}\geq 0, \Delta \delta_{j-2} \geq 0} (\varphi_{\jmtbt} - \varphi_{\jmh})(1-\alpha_{\jmh})(\Delta^{2}\delta_{\jmtbt})^{2}.
\end{align*}
Recollecting the estimate \eqref{eq:phi_difbd}, a bound can be obtained on the term $\mathcal{E}_{2}$ as follows:
\begin{align}\label{eq:E2bound}
    \abs{\mathcal{E}_{2}} &\leq \sum_{j=l-1}^{r+1}\abs{\varphi_{\jmtbt} - \varphi_{\jmh}} \left((\Delta\delta_{j-1})^{2}+\frac{1}{2}(\Delta^{2}\delta_{\jmtbt})^{2}\right) \\
& \leq   \norm{f_{u}}\lambda^{2} \gamma_{2} \frac{5}{2}\sum_{j=l-1}^{r+1}(
    \delta_{\jmh}+\delta_{\jmtbt})\left((\Delta\delta_{j-1})^{2}+\frac{1}{2}(\Delta^{2}\delta_{\jmtbt})^{2}\right) \notag\\
    &\spc +  2 \norm{f_{u}}\lambda^{2}\norm{f_{uk}}\sum_{j=l-1}^{r+1}\abs{\Delta k_{\jmtbt}} \left((\Delta\delta_{j-1})^{2}+\frac{1}{2}(\Delta^{2}\delta_{\jmtbt})^{2}\right) \notag\\
& \leq   \norm{f_{u}}\lambda^{2} \gamma_{2} \frac{205}{4}\sum_{j=l-1}^{r+1}(
\delta_{\jmh})^{3}
+  96 C_{u_{0}}^{2}\norm{f_{u}}\lambda^{2}\norm{f_{uk}}\sum_{j=l-1}^{r+1}\abs{\Delta k_{\jmtbt}}. \notag  
\end{align}
Similarly, the term $\tilde{\tilde{\mathcal{I}}}_{2}$ can be expanded as
\begin{align}\label{eq:B}
\tilde{\tilde{\mathcal{I}}}_{2} &= \mathcal{E}_{3} + \sum_{\Delta \delta_{j-1}\geq 0} \varphi_{\jmh}\alpha_{\jmtbt}(\Delta \delta_{j-1})^{2} +\sum_{\Delta \delta_{j-1}\geq 0, \Delta \delta_{j} < 0 }\varphi_{\jmh}\alpha_{\jmtbt}(\Delta \delta_{j-1})(\Delta \delta_{j})\\
&\spc +\frac{1}{2}\sum_{\Delta \delta_{j-1}\leq 0}\varphi_{\jmh}(\alpha_{\jmtbt}+\alpha_{j-\frac{5}{2}})(\Delta \delta_{j-1})^{2} \notag\\ & \spc - \frac{1}{2}\sum_{\Delta \delta_{j-1}\leq 0, 
\Delta \delta_{j}\geq 0}\varphi_{\jmh}\alpha_{\jmtbt}(\Delta \delta_{j-1})^{2} -  \frac{1}{2} \sum_{\Delta \delta_{j-1}\leq 0, \Delta 
\delta_{j-2}\geq 0}\varphi_{\jmh}\alpha_{j-\frac{5}{2}}(\Delta \delta_{j-1})^{2} \notag\\ & \spc -\frac{1}{2}\sum_{\Delta \delta_{j-2} \leq 0, \Delta \delta_{j-1} \leq 0} \varphi_{\jmh}\alpha_{j-\frac{5}{2}}(\Delta^{2}\delta_{\jmtbt})^{2}, \notag
\end{align}
where the term $\mathcal{E}_{3}$ admits an estimate of the form: 
\begin{align}\label{eq:E3bound}
    \abs{\mathcal{E}_{3}} &\leq  \norm{f_{u}}\lambda^{2} \gamma_{2} \frac{205}{4}\sum_{j=l-1}^{r+1}(
\delta_{\jmh})^{3}
+  96 C_{u_{0}}^{2} \norm{f_{u}}\lambda^{2}\norm{f_{uk}}\sum_{j=l-1}^{r+1}\abs{\Delta k_{\jmtbt}}.  
\end{align}
It is immediate to see that
\begin{nalign}\label{eq:alpha_phi_bds}
    \abs[\Big]{\alpha_{\jmh}-\frac{1}{2}} \leq \lambda \norm{f_{u}} \quad \mbox{and} \quad \abs[\Big]{\varphi_{\jmh}-\frac{1}{4}} \leq \lambda^{2} \norm{f_{u}}^{2}. 
\end{nalign}
In view of this estimate, we add and subtract $\frac{1}{2}$ to $\alpha_{\jmh}$ and  $\frac{1}{4}$ to $\varphi_{\jmh},$  thereby simplifying $\tilde{\mathcal{I}}_{2}$ and $\tilde{\tilde{\mathcal{I}}}_{2}$, as given below
\begin{align}\label{eq:A_simple}
\tilde{\mathcal{I}}_{2} &=   \sum_{\Delta \delta_{j-1}\leq 0} \varphi_{\jmh}(1-\alpha_{\jmh})(\Delta \delta_{j-1})^{2} \\ & \spc + \frac{1}{2}\sum_{\Delta \delta_{j-1}\geq 0}\varphi_{\jmh}\left((1-\alpha_{\jph})+(1-\alpha_{\jmh})\right)(\Delta \delta_{j-1})^{2} \notag\\
&\spc +\sum_{\Delta \delta_{j-1}\geq 0, \Delta \delta_{j} < 0 }\frac{1}{8}(\Delta \delta_{j-1})(\Delta \delta_{j})- \frac{1}{2}\sum_{\Delta \delta_{j-1}\geq 0, 
\Delta \delta_{j}\leq 0}\frac{1}{8}(\Delta \delta_{j-1})^{2} \notag\\
& \spc -  \frac{1}{2} \sum_{\Delta \delta_{j-1}\geq 0, \Delta 
\delta_{j-2}\leq 0}\frac{1}{8}(\Delta \delta_{j-1})^{2}- \frac{1}{2}\sum_{\Delta \delta_{j-2} \geq 0, \Delta \delta_{j-1} \geq 0} \frac{1}{8}(\Delta^{2}\delta_{\jmtbt})^{2} + \mathcal{E}_{4}, \notag
\end{align}
where the term $\mathcal{E}_{4}$ can be bounded as follows
\begin{nalign}
    \abs{\mathcal{E}_{4}} &\leq \abs{\mathcal{E}_{2}} + \frac{1}{2}\lambda\norm{f_{u}}\sum_{j=l-1}^{r+1}\left(\Delta^{2}\delta_{\jmtbt}\right)^{2},
\end{nalign} on which an application of \eqref{eq:E2bound} yields the estimate \eqref{eq:E4bd}.
Analogous substitutions for the term $\tilde{\tilde{\mathcal{I}}}_{2}$ in \eqref{eq:B}
yields 
\begin{nalign}\label{eq:Bsimple} \tilde{\tilde{\mathcal{I}}}_{2} &=  \sum_{\Delta \delta_{j-1}\geq 0} \varphi_{\jmh}\alpha_{\jmtbt}(\Delta \delta_{j-1})^{2} +  \frac{1}{2}\sum_{\Delta \delta_{j-1}\leq 0}\varphi_{\jmh}(\alpha_{\jmtbt}+\alpha_{j-\frac{5}{2}})(\Delta \delta_{j-1})^{2}\\
& \spc +\sum_{\Delta \delta_{j-1}\geq 0, \Delta \delta_{j} < 0 }\frac{1}{8}(\Delta \delta_{j-1})(\Delta \delta_{j}) - \frac{1}{2}\sum_{\Delta \delta_{j-1}\leq 0, 
\Delta \delta_{j}\geq 0}\frac{1}{8}(\Delta \delta_{j-1})^{2}\\
& \spc -  \frac{1}{2} \sum_{\Delta \delta_{j-1}\leq 0, \Delta 
\delta_{j-2}\geq 0}\frac{1}{8}(\Delta \delta_{j-1})^{2} -\frac{1}{2}\sum_{\Delta \delta_{j-2} \leq 0, \Delta \delta_{j-1} \leq 0} \frac{1}{8}(\Delta^{2}\delta_{\jmtbt})^{2} + \mathcal{E}_{5}, 
\end{nalign}
where the term $\mathcal{E}_{5}$ can be bounded as follows
\begin{nalign}\abs{\mathcal{E}_{5}} & \leq \abs{\mathcal{E}_{3}} + \frac{5}{8}\lambda\norm{f_{u}}\sum_{j=l-1}^{r+1}\left(\Delta^{2}\delta_{\jmtbt}\right),
\end{nalign}
which combined with \eqref{eq:E3bound} gives the estimate \eqref{eq:E5bd}.
Further, adding and subtracting $\displaystyle \frac{1}{4}$ to $\varphi_{\jmtbt}$  in the last term of \eqref{eq:Q1_mod} and using the expressions \eqref{eq:A_simple} and \eqref{eq:Bsimple}, we represent the term $\mathcal{Q}_1$ as follows
\begin{nalign}\label{eq:Q1_rearranged}
    \mathcal{Q}_{1} = \mathcal{R}_{1} + \mathcal{Q}_{1}^{*} -\mathcal{E}_{4} -\mathcal{E}_{5} -\mathcal{E}_{6},
\end{nalign}
where $\displaystyle \mathcal{E}_{6} := -\frac{1}{2}\sum_{j=l-1}^{r+1}\left(\varphi_{\jmtbt}-\frac{1}{4}\right)^{2}\left(\Delta^{2}\delta_{\jmtbt}\right)^{2}.$
Using the CFL condition \eqref{eq:cfl}, we can bound $\mathcal{E}_{6}$ as follows:\begin{nalign}\label{eq:E6bd}
\abs{\mathcal{E}_{6}} &\leq \frac{1}{2}(\lambda \norm{f_{u}})^{4}\sum_{j=l-1}^{r+1}\left(\Delta^{2}\delta_{\jmtbt}\right)^{2} \leq \frac{1}{16}  \lambda \norm{f_{u}}\sum_{j=l-1}^{r+1}\left(\Delta^{2}\delta_{\jmtbt}\right)^{2}.
\end{nalign}
\end{proof}
\begin{lemma}\label{lemma:Q2}
    Let $\{u_{j}\}_{\jinz}$ be as in Lemma \ref{lemma:mainlemma_inseq}. The term $\mathcal{Q}_{2}$ in \eqref{eq:Q1Q3Q3} can be expressed as \begin{nalign}\label{eq:Q2rearranged}
    \mathcal{Q}_{2} = \mathcal{Q}_{2}^{*} + \mathcal{E}_{7},
\end{nalign}
where $\displaystyle
    \mathcal{Q}_{2}^{*} := \frac{1}{64}\left(2\sum_{j=l-1}^{r+1}\left(\Delta^{2}\delta_{\jmtbt}\right)^{2}- \sum_{j=l-1}^{r+1}\left(\Delta^{2}s_{j-1}\right)^{2} \right)$ and \begin{nalign}\label{eq:E7bd}
    \abs{\mathcal{E}_{7}} \leq \frac{1}{8} \lambda\norm{f_{u}}\sum_{j=l-1}^{r+1}\left(\Delta^{2}\delta_{\jmtbt}\right)^{2}. 
\end{nalign} Further, for $\mathcal{Q}_{1}^{*}$ as in \eqref{lemma:Q1rewrite}, we have the lower bound \begin{nalign}\label{eq:Q1starpQ2star}
    \mathcal{Q}_{1}^{*}+ \mathcal{Q}_{2}^{*} \geq \frac{1}{2048}\sum_{j=l-1}^{r+1}\left(\Delta^{2}\delta_{\jmtbt}\right)^{2}.
\end{nalign}
\end{lemma}
\begin{proof} 
Adding and subtracting $\frac{1}{4}$ to $\varphi_{\jmtbt},$ the identity \eqref{eq:Q2rearranged} is immediate, 
where 
\begin{align}\label{eq:E7defn}
    \mathcal{E}_{7} := \frac{1}{4}\left(2\sum_{j=l-1}^{r+1}\left(\varphi_{\jmtbt}-\frac{1}{4}\right)^{2}\left(\Delta^{2}\delta_{\jmtbt}\right)^{2}-\sum_{j=l-1}^{r+1}\left(\varphi_{\jmtbt}-\frac{1}{4}\right)^{2}\left(\Delta^{2}s_{j-1}\right)^{2} \right).
\end{align}
Further, invoking Lemma 4 from \cite{popov2006a}, we have
$$\sum_{j=l-1}^{r+1}\left(\Delta^{2}s_{j-1}\right)^{2}  \leq 2\sum_{j=l-1}^{r+1}\left(\Delta^{2}\delta_{\jmtbt}\right)^{2},$$
which, when applied to  \eqref{eq:E7defn}, implies that
\begin{nalign}
    \abs{\mathcal{E}_{7}} &\leq \lambda^{4}\norm{f_{u}}^{4}\sum_{j=l-1}^{r+1}\left(\Delta^{2}\delta_{\jmtbt}\right)^{2} \leq \frac{1}{8} \lambda\norm{f_{u}}\sum_{j=l-1}^{r+1}\left(\Delta^{2}\delta_{\jmtbt}\right)^{2}, 
\end{nalign}
where the last inequality follows from the CFL condition \eqref{eq:cfl}. Furthermore, following the proof of Lemma 4 in \cite{popov2006a}, we obtain the desired lower bound \eqref{eq:Q1starpQ2star}, thus completing the proof. 
\end{proof}
Towards our objective of getting a lower bound on $\mathcal{D},$ we now reformulate the term $R_{1}$ from \eqref{eq:R1}. By applying summation by parts and subsequently adding and subtracting $\varphi_{\jmh}$ to $\varphi_{\jmtbt},$  we rewrite the second term of \eqref{eq:R1} as
\begin{align}\label{eq:R1secterm}
&\sum_{j=l-1}^{r+1} \varphi_{\jmh}\Delta s_{\jmh} \left( \delta_{\jmh}(\Delta \alpha_{j}- \Delta \alpha_{j-1})\right) \\&= -\sum_{j=l-1}^{r+1}\Delta \alpha_{j-1}\left(\varphi_{\jmh}\Delta s_{\jmh}\delta_{\jmh}- \varphi_{\jmtbt}\Delta s_{\jmtbt}\delta_{\jmtbt}\right) \notag\\ 
& = - \sum_{j=l-1}^{r+1}\varphi_{\jmh}\Delta \alpha_{j-1}\left(\Delta s_{\jmh}\delta_{\jmh}- \Delta s_{\jmtbt}\delta_{\jmtbt}\right) + \mathcal{E}_{8}, \notag
\end{align}
where $
     \mathcal{E}_{8}:=\sum_{j=l-1}^{r+1} \Delta\alpha_{j-1}\left(\varphi_{\jmtbt}-\varphi_{\jmh}\right)\Delta s_{\jmtbt}\delta_{\jmtbt},$ and is bounded as follows
\begin{nalign}\label{eq:E8bound}    \abs{\mathcal{E}_{8}} \leq 15 \lambda^{3}\norm{f_{u}}^{2}\gamma_{2}\sum_{j=l-1}^{r+1}\left(\delta_{\jmtbt}\right)^{3} + 16 \lambda^{3}\norm{f_{u}}^{2}\norm{f_{uk}}C_{u_{0}}^{2}\sum_{j=l-1}^{r+1}\abs{\Delta k_{\jmtbt}}.
\end{nalign}
Now, in view of  \eqref{eq:R1secterm},  \eqref{eq:R1} can be written as 
\begin{align}\label{eq:R1minusE8}
    \mathcal{R}_{1}
    &=  \mathcal{E}_{8} +\sum_{j=l-1}^{r+1}\alpha_{\jmtbt}(1-\alpha_{\jmh})(\Delta \delta_{j-1})^{2} \\ & \spc + \sum_{j=l-1}^{r+1}\varphi_{\jmh}\Delta \alpha_{j-1}\left(\Delta s_{\jmtbt}\delta_{\jmtbt}- \Delta s_{\jmh}\delta_{\jmh}\right)  \notag\\
    &\spc - \sum_{\Delta \delta_{j-1}\leq 0} \varphi_{\jmh}\left(1-\alpha_{\jmh} + \frac{\alpha_{\jmtbt}+\alpha_{j-\frac{5}{2}}}{2}\right)(\Delta \delta_{j-1})^{2} \notag \\
    &\spc -\sum_{\Delta \delta_{j-1}\geq 0} \varphi_{\jmh}\left(1+\alpha_{\jmtbt}-\frac{\alpha_{\jph}+ \alpha_{\jmh}}{2}\right)(\Delta \delta_{j-1})^{2}. \notag
\end{align}
In \eqref{eq:R1minusE8}, adding and subtracting $\alpha_{\jmtbt}$ to $\alpha_{j-\frac{5}{2}}$ in the fourth term and $\alpha_{\jmh}$ to $\alpha_{\jph}$ in the fifth term, we write
\begin{align}\label{eq:R1rewrite}
     \mathcal{R}_{1} = \mathcal{E}_{8} + \mathcal{R}_{2} + \mathcal{R}_{3} + \mathcal{R}_{4} + \mathcal{R}_{5},
\end{align}
where 
\begin{nalign}\label{eq:R2R3R4R5}
    \mathcal{R}_{2} &:=  \sum_{j=l-1}^{r+1}\varphi_{\jmh}\Delta \alpha_{j-1}\left(\Delta s_{\jmtbt}\delta_{\jmtbt}- \Delta s_{\jmh}\delta_{\jmh}\right),  \\
    \mathcal{R}_{3} &:= \sum_{j=l-1}^{r+1}\left(\alpha_{\jmtbt}(1-\alpha_{\jmh}) - \varphi_{\jmh}\left(1-\Delta \alpha_{j-1}\right) \right)(\Delta \delta_{j-1})^{2},\\
    \mathcal{R}_{4} &:=  \sum_{\Delta \delta_{j-1}\leq 0} \varphi_{\jmh}\left( \frac{\alpha_{\jmtbt}-\alpha_{j-\frac{5}{2}}}{2}\right)(\Delta \delta_{j-1})^{2} \quad \mbox{and}\\
    \mathcal{R}_{5} &:= \sum_{\Delta \delta_{j-1}\geq 0} \varphi_{\jmh}\left(\frac{\alpha_{\jph}- \alpha_{\jmh}}{2}\right)(\Delta \delta_{j-1})^{2}. 
\end{nalign}
Next, $\mathcal{R}_{4}$ admits a lower bound \begin{nalign}\label{eq:R4_lb}
    \mathcal{R}_{4} \geq -16(C_{u_{0}})^{2}\lambda\norm{f_{ku}}\sum_{j=l-2}^{r} \abs{\Delta k_{j-\frac{3}{2}}}=:\tilde{\mathcal{R}}_{4}.
\end{nalign}
This holds true because, by \eqref{eq:alphadif_lb}, we have
 \begin{align*}
    \mathcal{R}_{4} & \geq  \sum_{\Delta \delta_{j-1}\leq 0} \frac{\varphi_{\jmh}}{2}\left(\frac{3}{8}\lambda \gamma_{1} \left(\delta_{\jmtbt}^{n}+\delta_{j-\frac{5}{2}}^{n}\right) + \lambda \Delta k_{j-\frac{5}{2}}f_{ku}(\bar{k}_{j-\frac{5}{2}}, \bar{u}_{\jmtbt})\right)(\Delta \delta_{j-1})^{2}\\
    &\geq \frac{3}{16}\lambda \gamma_{1}\sum_{\Delta \delta_{j-1}\leq 0} \varphi_{\jmh} \left(\delta_{\jmtbt}^{n}+\delta_{j-\frac{5}{2}}^{n}\right)(\Delta \delta_{j-1})^{2} \\ & \spc + \sum_{\Delta \delta_{j-1}\leq 0}\lambda \Delta k_{j-\frac{5}{2}}f_{ku}(\bar{k}_{j-\frac{5}{2}}, \bar{u}_{\jmtbt})(\Delta \delta_{j-1})^{2}\\
    &\geq  \sum_{\Delta \delta_{j-1}\leq 0}\lambda \Delta k_{j-\frac{5}{2}}f_{ku}(\bar{k}_{j-\frac{5}{2}}, \bar{u}_{\jmtbt})(\Delta \delta_{j-1})^{2},
\end{align*}
and also by noting the fact that $(\Delta \delta_{j-1})^{2} \leq 16 (C_{u_{0}})^{2}.$ Similar arguments on $\mathcal{R}_5$ yield
\begin{nalign}\label{eq:R5_lb}
    \mathcal{R}_5 \geq -16(C_{u_{0}})^{2}\lambda\norm{f_{ku}}\sum_{j=l}^{r+2} \abs{\Delta k_{\jmtbt}}=:\tilde{\mathcal{R}}_{5}
\end{nalign}
Now, using Lemmas \ref{lemma:Q1rewrite} and \ref{lemma:Q2}, along with with the expression \eqref{eq:R1rewrite} in \eqref{eq:D}, we  write
\begin{nalign}\label{eq:D_lb}
    \mathcal{D} =   \mathcal{Q}_{1}^{*} + \mathcal{Q}_{2}^{*}+\mathcal{Q}_{3}^{*}+ \mathcal{R}_{4}+ \mathcal{R}_{5}+ P_{2} +  \mathcal{E}_{9}, 
\end{nalign}
where $Q_{3}^{*} := \mathcal{Q}_{3}+ \mathcal{R}_{2}+ \mathcal{R}_{3}$ and $\mathcal{E}_{9} :=  \mathcal{E}_{1}-\mathcal{E}_{4} -\mathcal{E}_{5} -\mathcal{E}_{6}+  \mathcal{E}_{7} +\mathcal{E}_{8}.$    Further, combining the expressions \eqref{eq:E1bd},\eqref{eq:E4bd}, \eqref{eq:E5bd}, \eqref{eq:E6bd}, \eqref{eq:E7bd} and \eqref{eq:E8bound}, using the CFL condition \eqref{eq:cfl} and after suitable algebraic manipulations, we obtain the following bound 
\begin{align}\label{eq:E9bound}
    \mathcal{E}_{9} 
&\geq -255\lambda^{2}\norm{f_{u}}\gamma_{2}\sum_{j=l-1}^{r+1}(\delta_{\jmh})^{3} - \frac{21}{16}\lambda\norm{f_{u}}\sum_{j=l-1}^{r+1}\left(\Delta^{2}\delta_{\jmtbt}\right)^{2} \\
&\spc - \left(24+192\norm{f_{u}}+16\lambda\norm{f_{u}}^{2}\right)C_{u_{0}}^{2}\lambda^{2}\norm{f_{uk}}\sum_{j=l-1}^{r+1}\abs{\Delta k_{\jmtbt}}. \notag
\end{align}  Next, we establish a lower bound on the term $\mathcal{Q}_{3}^{*}.$
\begin{lemma}\label{lemma:Q3pr_bd} The term $\mathcal{Q}_{3}^{*}$ in \eqref{eq:D_lb} admits the lower bound
\begin{nalign}\label{eq:q3primebd}
    \mathcal{Q}_{3}^{*} \geq \frac{3}{64}\sum_{j=l-1}^{r+1}\lambda \gamma_{1} (\delta_{\jmh})^{3} - 6\lambda \norm{f_{ku}}C_{u_{0}}^{2}\sum_{j=l-2}^{r+1}\abs{\Delta k_{\jmh}}+ \mathcal{E}_{10},
\end{nalign}
where $C_{u_0}$ is as in \eqref{Linf_bd} and $\mathcal{E}_{10}$ is bounded as follows
\begin{align*}
    \abs{\mathcal{E}_{10}} \leq \left(28\lambda\norm{f_{u}}\right) \lambda\norm{f_{uk}} C_{u_{0}}^{2}\sum_{j=l-1}^{r+1} \abs{\Delta k_{\jmtbt}}+  \left(\frac{35}{2}\lambda\norm{f_{u}}\right)\lambda \gamma_{2} \sum_{j=l-1}^{r+1}(\delta_{\jmh})^{3}.
\end{align*}
\end{lemma}

\begin{proof}Recalling the notation $\varphi_{\jmh}:= \alpha_{\jmh}(1-\alpha_{\jmh})$ and using it in $\mathcal{R}_{3}$ of \eqref{eq:R2R3R4R5}, we write $\mathcal{Q}_{3}^{*}$ as 
\begin{align*}  \mathcal{Q}_{3}^{*} &= \mathcal{Q}_{3}+ \mathcal{R}_{2}- \sum_{j=l-1}^{r+1}(1-\alpha_{\jmh})^{2}(\alpha_{\jmh}-\alpha_{\jmtbt})(\Delta \delta_{j-1})^{2}.
\end{align*}
Writing $\varphi_{\jmh}= \frac{1}{4}+( \varphi_{\jmh} -\frac{1}{4}),$ $\alpha_{\jmtbt}=\frac{1}{2}+(\alpha_{\jmtbt}-\frac{1}{2}) \, \mbox{and} \, \alpha_{\jmh}=\frac{1}{2}+(\alpha_{\jmh}-\frac{1}{2}),$ the term $\mathcal{Q}_{3}^{*}$ can be expressed as
\begin{nalign}\label{eq:Q3prime}
    \mathcal{Q}_{3}^{*} &:= \sum_{j=l-1}^{r+1} \Delta \alpha_{j-1}\left(\frac{1}{2}(\delta_{\jmtbt})^{2} + \frac{1}{2}(\delta_{\jmh})^{2}\right) \\ & \spc + \sum_{j=l-1}^{r+1}\frac{1}{4}\Delta \alpha_{j-1}\left(\Delta s_{\jmtbt}\delta_{\jmtbt}- \Delta s_{\jmh}\delta_{\jmh}\right) - \sum_{j=l-1}^{r+1}\frac{1}{4}\Delta \alpha_{j-1}(\Delta \delta_{j-1})^{2}+ \mathcal{E}_{10}\\
    &= \frac{1}{4}\sum_{j=l-1}^{r+1} \Delta \alpha_{j-1} z_{j-1} + \mathcal{E}_{10},
\end{nalign}
where
\begin{align*}
z_{j-1} &:= 2(\delta_{\jmtbt})^{2} + 2(\delta_{\jmh})^{2}+ \Delta s_{\jmtbt}\delta_{\jmtbt}- \Delta s_{\jmh}\delta_{\jmh}-(\Delta \delta_{j-1})^{2} \quad \mbox{and}\\
    \mathcal{E}_{10} &:= \sum_{j=l-1}^{r+1} \Delta \alpha_{j-1} \left[
    \left(\alpha_{\jmtbt}-\frac{1}{2}\right)(\delta_{\jmtbt}^{n})^{2} + \left(\frac{1}{2}-\alpha_{\jmh}\right)(\delta_{\jmh}^{n})^{2} \right. \\
    &\spc \left. + \left(\varphi_{\jmh}-\frac{1}{4}\right)\left(\Delta s_{\jmtbt}\delta_{\jmtbt}
    - \Delta s_{\jmh}\delta_{\jmh}\right) 
    \right. \\ & \left. \spc -\left(\frac{3}{4} +(\alpha_{\jmh})^{2}-2\alpha_{\jmh}\right)(\Delta \delta_{j-1})^{2}\right].
\end{align*}
Recalling the notation $\alpha_{\jmh} :=  \frac{1}{2}+\lambda\bar{a}_{\jmh},$ it is clear from \eqref{eq:abar_dif} and \eqref{eq:ubardiff_bd} that $\Delta \alpha_{j-1}$ admits a bound $ \abs{\Delta \alpha_{j-1}} \leq  \lambda\norm{f_{uk}}\abs{\Delta k_{\jmtbt}} + \frac{5}{4}\lambda \gamma_{2}(\delta_{\jmh}+\delta_{\jmtbt}).$ Consequently, we obtain the following bound on $\mathcal{E}_{10}:$ 
\begin{align}\label{eq:E10bd}
\abs{\mathcal{E}_{10}} &\leq \left(\frac{7}{2}\lambda\norm{f_{u}}\right)  \sum_{j=l-1}^{r+1} \lambda\norm{f_{uk}}\abs{\Delta k_{\jmtbt}}\left((\delta_{\jmh})^{2}+(\delta_{\jmtbt})^{2}\right)  \\&\spc+  \left(\frac{7}{2}\lambda\norm{f_{u}}\right)\sum_{j=l-1}^{r+1}\frac{5}{4}\lambda \gamma_{2}(\delta_{\jmh}+\delta_{\jmtbt})\left((\delta_{\jmh})^{2}+(\delta_{\jmtbt})^{2}\right) \notag \\
& \leq \left(\frac{7}{2}\lambda\norm{f_{u}}\right)8\lambda\norm{f_{uk}} C_{u_{0}}^{2}\sum_{j=l-1}^{r+1} \abs{\Delta k_{\jmtbt}}+  \left(\frac{7}{2}\lambda\norm{f_{u}}\right)\frac{5}{4}\lambda \gamma_{2}6\sum_{j=l-1}^{r+1}(\delta_{\jmh})^{3}\notag\\
& \leq \left(28\lambda\norm{f_{u}}\right) \lambda\norm{f_{uk}} C_{u_{0}}^{2}\sum_{j=l-1}^{r+1} \abs{\Delta k_{\jmtbt}}+  \left(\frac{105}{4}\lambda\norm{f_{u}}\right)\lambda \gamma_{2} \sum_{j=l-1}^{r+1}(\delta_{\jmh})^{3}.\notag
\end{align}
Next, to obtain a lower bound on $\mathcal{Q}_{3}^{*}$, we need a bound on the first term in the RHS of \eqref{eq:Q3prime}. To this end, we first make the following observation on $\displaystyle z_{j-1}$ in \eqref{eq:Q3prime}:\begin{nalign}\label{eq:z_lb}
z_{j-1} & =  \left(\delta_{\jmtbt}+\delta_{\jmh}\right)^{2} + s_{j-1}(\delta_{\jmtbt}+\delta_{\jmh}) - \delta_{\jmtbt}  s_{j-2}-\delta_{\jmh}s_{j}\\
&\geq 2\delta_{\jmh}\delta_{\jmtbt} + s_{j-1}(\delta_{\jmtbt}+\delta_{\jmh}).
\end{nalign}
\cblue{This establishes that $z_{j-1} \geq 0.$}
Furthermore, combining the expression \eqref{eq:z_lb} with the $\mathrm{L}^{\infty}$-stability \eqref{Linf_bd} of the scheme \eqref{eq:NTscheme}, we deduce that $\abs{z_{j-1}} \leq 24C_{u_0}^2,$ for $\jinz.$ Next, we obtain a bound
\begin{nalign}\label{eq:lbd_deltaalph}
\Delta \alpha_{j-1} z_{j-1}+  \Delta \alpha_{j-2} z_{j-2} &\geq  \frac{3}{8}\lambda \gamma_{1}  (\delta_{\jmtbt})^{3}  \\ & \spc - 24 C_{u_0}^2\lambda \norm{f_{ku}}(\abs{\Delta k_{\jmtbt}}+ \abs{\Delta k_{j-\frac{5}{2}}}),
\end{nalign}
which is verified by considering the two possible cases.\\
\textbf{Case 1} ($\delta_{\jmh} \leq \delta_{\jmtbt}$ and  $\delta_{j-\frac{5}{2}} \leq \delta_{\jmtbt}$): In this case, as $s_{j-2} = \delta_{j-\frac{5}{2}},$ from the definition of $z_{j-1},$ it follows that
\begin{align*}
z_{j-1} &= (\delta_{\jmtbt})^{2}+(\delta_{\jmh})^{2}+2\delta_{\jmh}\delta_{\jmtbt} + s_{j-1}(\delta_{\jmh}+\delta_{\jmtbt})-\delta_{\jmtbt}s_{j-2} -\delta_{\jmh}s_{j}\\
&\geq (\delta_{\jmtbt})^{2} - \delta_{\jmtbt}\delta_{j-\frac{5}{2}}.
\end{align*}
This bound, combined with \eqref{eq:alphadif_lb} implies that
 \begin{align*}
      \Delta \alpha_{j-1} z_{j-1} &\geq \frac{3}{8}\lambda \gamma_{1} \left(\delta_{\jmh}+\delta_{\jmtbt}\right) \left((\delta_{\jmtbt})^{2} - \delta_{\jmtbt}\delta_{j-\frac{5}{2}}\right) -24 C_{u_{0}}^2\lambda \norm{f_{ku}}\abs{\Delta k_{\jmtbt}} \\
     &\geq \frac{3}{8}\lambda \gamma_{1}  \left((\delta_{\jmtbt})^{3} - (\delta_{\jmtbt})^{2}\delta_{j-\frac{5}{2}}\right) - 24 C_{u_{0}}^2\lambda \norm{f_{ku}}\abs{\Delta k_{\jmtbt}}.
 \end{align*}
Furthermore, in this case, we also have
$z_{j-2} \geq 2\delta_{j-\frac{5}{2}} \delta_{\jmtbt},$
which yields
 \begin{align*}
       \Delta \alpha_{j-2} z_{j-2} &\geq \frac{3}{8}\lambda \gamma_{1} 2\delta_{j-\frac{5}{2}} \delta_{\jmtbt}\left(\delta_{\jmtbt}^{n}+\delta_{j-\frac{5}{2}}^{n}\right) -24C_{u_{0}}^{2} \lambda \norm{f_{ku}} \abs{\Delta k_{j-\frac{5}{2}}}\\
        &\geq \frac{3}{4}\lambda \gamma_{1} \delta_{j-\frac{5}{2}} (\delta_{\jmtbt})^{2} -24C_{u_{0}}^{2} \lambda \norm{f_{ku}} \abs{\Delta k_{j-\frac{5}{2}}}.
\end{align*}
This leads to the conclusion
\begin{nalign}\label{eq:case3}
    \Delta \alpha_{j-1} z_{j-1}+  \Delta \alpha_{j-2} z_{j-2} \geq  \frac{3}{8}\lambda \gamma_{1}  (\delta_{\jmtbt})^{3}  - 24 C_{u_{0}}^2\lambda \norm{f_{ku}}(\abs{\Delta k_{\jmtbt}}+ \abs{\Delta k_{j-\frac{5}{2}}}),
\end{nalign} for Case 1. \\
\textbf{Case 2} ($\delta_{\jmh} \geq \delta_{\jmtbt}$ or $\delta_{j-\frac{5}{2}} \geq \delta_{\jmtbt}$): In the case when $\delta_{\jmh} \geq \delta_{\jmtbt}$, we have $s_{j-1} = \delta_{\jmtbt}.$ Consequently, applying \eqref{eq:z_lb} yields $z_{j-1} \geq  4(\delta_{\jmtbt})^{2}.$ Further, using \eqref{eq:alphadif_lb}, we can write 
\begin{align*}
 \Delta \alpha_{j-1} z_{j-1} &\geq  \frac{3}{8}\lambda \gamma_{1} 4(\delta_{\jmtbt})^{2}\left(\delta_{\jmh}^{n}+\delta_{\jmtbt}^{n}\right) - 24 C_{u_0}^2 \lambda \norm{f_{ku}} \abs{\Delta k_{\jmtbt}}\\
&\geq  \frac{3}{2}\lambda \gamma_{1} (\delta_{\jmtbt})^{3} - 24 C_{u_{0}}^2 \lambda \norm{f_{ku}} \abs{\Delta k_{\jmtbt}}.
\end{align*} Moreover, \eqref{eq:alphadif_lb} also implies that $\Delta \alpha_{j-2} z_{j-2} \geq -24 C_{u_{0}}^2 \lambda \norm{f_{ku}} \abs{\Delta k_{j-\frac{5}{2}}}.$ Thus,
\begin{nalign}\label{eq:case1}
    \Delta \alpha_{j-1} z_{j-1}+\Delta \alpha_{j-2} z_{j-2} &\geq \frac{3}{4}\lambda \gamma_{1} (\delta_{\jmtbt})^{3}- 24 C_{u_{0}}^2 \lambda \norm{f_{ku}} (\abs{\Delta k_{\jmtbt}} + \abs{\Delta k_{j-\frac{5}{2}}}),
\end{nalign} for the case when $\delta_{\jmh} \geq \delta_{\jmtbt}.$
Next, for the case when
$\delta_{j-\frac{5}{2}} \geq \delta_{\jmtbt},$ we have $s_{j-2} = \delta_{\jmtbt}.$ Therefore, using \eqref{eq:z_lb}, we have $z_{j-2} \geq  4(\delta_{\jmtbt})^{2}.$
Applying \eqref{eq:alphadif_lb} again, we obtain
\begin{align*}
 \Delta \alpha_{j-2} z_{j-2} &\geq  \frac{3}{8}\lambda \gamma_{1} 4(\delta_{\jmtbt})^{2}\left(\delta_{\jmtbt}^{n}+\delta_{j-\frac{5}{2}}^{n}\right)  -24 C_{u_{0}}^{2}\lambda \norm{f_{ku}}\abs{\Delta k_{j-\frac{5}{2}}}\\
&\geq  \frac{3}{2}\lambda \gamma_{1} (\delta_{\jmtbt})^{3} - 24 C_{u_{0}}^{2}\lambda \norm{f_{ku}}\abs{\Delta k_{j-\frac{5}{2}}}.
\end{align*}
Again, using \eqref{eq:z_lb} it follows that $ \Delta \alpha_{j-1} z_{j-1} \geq -24 C_{u_{0}}^2 \lambda \norm{f_{ku}} \abs{\Delta k_{j-\frac{3}{2}}},$ leading to
\begin{nalign}\label{eq:case2}
    \Delta \alpha_{j-1} z_{j-1}+\Delta \alpha_{j-2} z_{j-2} &\geq \frac{3}{4}\lambda \gamma_{1} (\delta_{\jmtbt})^{3}- 24 C_{u_{0}}^2 \lambda \norm{f_{ku}} (\abs{\Delta k_{\jmtbt}} + \abs{\Delta k_{j-\frac{5}{2}}}),
\end{nalign}
thus concluding Case 2.
\par
Combining both the cases, the estimates \eqref{eq:case3},\eqref{eq:case1}  and \eqref{eq:case2} establish the bound \eqref{eq:lbd_deltaalph}. Furthermore, using \eqref{eq:lbd_deltaalph}, we obtain
\begin{nalign}\label{eq:alphazsum_lb}\sum_{j=l-1}^{r+1}\Delta \alpha_{j-1} z_{j-1} &=\frac{1}{2}\left(\sum_{j=l}^{r+1}\Delta \alpha_{j-1} z_{j-1}+\sum_{j=l+1}^{r+2}\Delta \alpha_{j-2} z_{j-2}\right)\\&= \frac{1}{2}\left(\sum_{j=l}^{r+2}\left(\Delta \alpha_{j-1} z_{j-1}+\Delta \alpha_{j-2} z_{j-2}\right)\right)\\ &\geq \frac{3}{16}\sum_{j=l-1}^{r+1}\lambda \gamma_{1} (\delta_{\jmh})^{3} - 24\lambda \norm{f_{ku}}C_{u_{0}}^{2}\sum_{j=l-2}^{r+1}\abs{\Delta k_{\jmh}}.\end{nalign}
The bound \eqref{eq:alphazsum_lb} together with \eqref{eq:Q3prime} yields the desired estimate \eqref{eq:q3primebd} on $Q_{3}^{*}$.
\end{proof}
Finally, to conclude this section on non-decreasing sequences, we provide the proof of Lemma \ref{lemma:mainlemma_inseq}.
\begin{proof}[Proof of Lemma \ref{lemma:mainlemma_inseq}]
By invoking the bound \eqref{eq:Q1starpQ2star} for $\mathcal{Q}_{1}^{*}+\mathcal{Q}_{2}^{*}$ from Lemma \ref{lemma:Q2} and the bound \eqref{eq:q3primebd} for $\mathcal{Q}_{3}^{*}$ from Lemma \ref{lemma:Q3pr_bd}, in the expression \eqref{eq:D_lb}, we obtain
\begin{nalign}\label{eq:D_final}
\mathcal{D} &\geq  \frac{3}{64}\sum_{j=l-1}^{r+1}\lambda \gamma_{1} (\delta_{\jmh})^{3} +  \frac{1}{2048}\sum_{j=l-1}^{r+1}\left(\Delta^{2}\delta_{\jmtbt}\right)^{2} + \mathcal{E}_{9}+ \mathcal{E}_{10} \\&\spc {+\tilde{\mathcal{R}}_{4}+\tilde{\mathcal{R}}_{5} + P_{2} - 6\lambda \norm{f_{ku}}C_{u_{0}}^{2}\sum_{j=l-2}^{r+1}\abs{\Delta k_{\jmh}}.}
\end{nalign}
Now, using the bounds \eqref{eq:E9bound} and   \eqref{eq:E10bd} for $\mathcal{E}_{9}$ and $\mathcal{E}_{10},$ respectively, and the CFL condition \eqref{eq:cfl}, we further simplify \eqref{eq:D_final} as follows 

\begin{align}\label{eq:Dfinal_2}
     \mathcal{D} &\geq \frac{3}{64}\lambda \gamma_{1}\sum_{j=l-1}^{r+1}(\delta_{\jmh})^{3}  +  \frac{1}{2048}\sum_{j=l-1}^{r+1}\left(\Delta^{2}\delta_{\jmtbt}\right)^{2}\\
    &\spc-\frac{1125}{4}\lambda^{2}\norm{f_{u}}\gamma_{2}\sum_{j=l-1}^{r+1}(\delta_{\jmh})^{3}- \frac{21}{16}\lambda\norm{f_{u}}\sum_{j=l-1}^{r+1}\left(\Delta^{2}\delta_{\jmtbt}\right)^{2}\notag\\
    &\spc - \left(24+220\norm{f_{u}} +16\lambda\norm{f_{u}}^{2}\right)C_{u_{0}}^{2}\lambda^{2}\norm{f_{uk}}\sum_{j=l-1}^{r+1}\abs{\Delta k_{\jmtbt}} \notag\\
    & \spc {+\tilde{\mathcal{R}}_{4}+\tilde{\mathcal{R}}_{5} + P_{2} - 6\lambda \norm{f_{ku}}C_{u_{0}}^{2}\sum_{j=l-2}^{r+1}\abs{\Delta k_{\jmh}}}\notag\\
    & \geq \frac{15}{1600} \lambda \gamma_{1}\sum_{j=l-1}^{r+1}(\delta_{\jmh})^{3} +  \frac{1}{6400}\sum_{j=l-1}^{r+1}\left(\Delta^{2}\delta_{\jmtbt}\right)^{2} \notag\\
    & \spc - \left(24+220\norm{f_{u}} +16\lambda\norm{f_{u}}^{2}\right)C_{u_{0}}^{2}\lambda^{2}\norm{f_{uk}}\sum_{j=l-1}^{r+1}\abs{\Delta k_{\jmtbt}}  \notag\\
    & \spc {+\tilde{\mathcal{R}}_{4}+\tilde{\mathcal{R}}_{5} + P_{2} - 6\lambda \norm{f_{ku}}C_{u_{0}}^{2}\sum_{j=l-2}^{r+1}\abs{\Delta k_{\jmh}}}, \notag
\end{align}
where the last inequality follows from the CFL condition \eqref{eq:cfl}.
Adding and subtracting $\sum_{j=l-1}^{r+1}(\delta_{j-1}^{\prime})^2$ to $\mathcal{D}$ in \eqref{eq:Dfinal_2} and subsequently  using \eqref{eq:djpr_minus_djprpr} and the CFL condition \eqref{eq:cfl}, we have 
\begin{align}\label{eq:posjumpsfinal}
    \sum_{j=l-1}^{r+1}\left((\delta_{\jmh})^2- (\delta_{j-1}^{\prime})^2\right) &\geq \sum_{j=l-1}^{r+1}\left((\delta_{\jmh})^2- (\delta_{j-1}^{\prime\prime})^2\right) - \abs[\Big]{\sum_{j=l-1}^{r+1}(\delta_{j-1}^{\prime})^2- (\delta_{j-1}^{\prime\prime})^2}\\
    &\geq \frac{1}{500} \lambda \gamma_{1}\sum_{j=l-1}^{r+1} (\delta_{\jmh})^{3} +  \frac{1}{6400}\sum_{j=l-1}^{r+1}\left(\Delta^{2}\delta_{\jmtbt}\right)^{2} - P_{1}\notag\\
    &\spc - \left(24+220\norm{f_{u}} \right.\notag \\ & \spc \spc \left. +16\lambda\norm{f_{u}}^{2}\right)C_{u_{0}}^{2}\lambda^{2}\norm{f_{uk}}\sum_{j=l-1}^{r+1}\abs{\Delta k_{\jmtbt}}  \notag\\
    & \spc {+\tilde{\mathcal{R}}_{4}+\tilde{\mathcal{R}}_{5} + P_{2} - 6\lambda \norm{f_{ku}}C_{u_{0}}^{2}\sum_{j=l-2}^{r+1}\abs{\Delta k_{\jmh}}}.\notag
\end{align}
Further, substituting the values of $P_{1}, \tilde{\mathcal{R}}_{4}$ and $\tilde{\mathcal{R}}_{5}$ from \eqref{eq:djpr_minus_djprpr}, \eqref{eq:R4_lb} and \eqref{eq:R5_lb}, respectively, and using the estimate  \eqref{eq:P2bound} for $P_{2},$ we write
\begin{align}\label{eq:incr_prefinal}
    \sum_{j=l-1}^{r+1}\left((\delta_{\jmh})^2- (\delta_{j-1}^{\prime})^2\right)& \geq \frac{1}{500} \lambda \gamma_{1}\sum_{j=l-1}^{r+1} (\delta_{\jmh})^{3} +  \frac{1}{6400}\sum_{j=l-1}^{r+1}\left(\Delta^{2}\delta_{\jmtbt}\right)^{2}\\ & \spc - \Theta\sum_{j=l-2}^{r+2}\abs{\Delta k_{\jmtbt}},\notag
\end{align}
where
\begin{align*}
    \Theta &:= 24\lambda^{2}(C_{u_{0}})^{2}\norm{f_{uk}}+38\lambda(C_{u_{0}})^{2}\norm{f_{ku}} + \left(236(C_{u_{0}})^{2}+16\lambda\norm{f_{u}}\right)\lambda^{2}\norm{f_{u}}\norm{f_{uk}}\\
    & \spc + \left(16\lambda^{2} C_{u_0}\norm{k}+44\lambda^{2} (C_{u_0})^{2}\gamma_{2}\norm{f_{u}}+8\lambda\norm{f_k}\norm{k}+24 C_{u_0}\right)\lambda\norm{f_{k}}.
\end{align*}
 Finally, we obtain the desired estimate \eqref{eq:incr_seq_result} since $ \displaystyle  \sum_{j=l-1}^{r+1}\left((\delta_{\jmh})^2- (\delta_{j-1}^{\prime})^2_{+}\right)  \geq  \sum_{j=l-1}^{r+1}\left((\delta_{\jmh})^2- (\delta_{j-1}^{\prime})^2\right).$  
\end{proof}

\textit{\textbf {A summary of the proof of Lemma \ref{lemma:mainlemma_inseq}}:} To derive the desired estimate for the term $\sum_{j=l-1}^{r+1}\left((\delta_{\jmh})^2- (\delta_{j-1}^{\prime})^2\right)$, we begin by introducing the modified jumps $\delta_{j}^{\prime \prime}$ (see \eqref{eq:deltadobpr})    and show in Lemma \ref{lemma:deltamin_deltapr_est} that the term $\abs[\Big]{\sum_{j=l-1}^{r+1}(\delta_{j-1}^{\prime})^2- \sum_{j=l-1}^{r+1}(\delta_{j-1}^{\prime\prime})^2}$ is appropriately bounded. With this estimate in hand, we focus on the term $\mathcal{D} := \sum_{j=l-1}^{r+1}\left((\delta_{\jmh})^2 - (\delta_{j-1}^{\prime\prime})^2\right),$ which we decompose as $\mathcal{D} = \mathcal{Q}_{1} + \mathcal{Q}_{2} + \mathcal{Q}_{3}+ P_{2} + \mathcal{E}_{1}.$  Next, in Lemmas \ref{lemma:Q1rewrite} and \ref{lemma:Q2}, we reformulate the terms $\mathcal{Q}_{1}$ and $\mathcal{Q}_{2}$ as $\mathcal{Q}_{1} = \mathcal{R}_{1} + \mathcal{Q}_{1}^{*} -\mathcal{E}_{4} -\mathcal{E}_{5} -\mathcal{E}_{6}$  and $\mathcal{Q}_{2} = \mathcal{Q}_{2}^{*} + \mathcal{E}_{7}.$ We then rewrite the term $\mathcal{R}_{1}$ in $\mathcal{Q}_{1}$ as $ \mathcal{R}_{1} = \mathcal{E}_{8} + \mathcal{R}_{2} + \mathcal{R}_{3} + \mathcal{R}_{4} + \mathcal{R}_{5}$ (see \eqref{eq:R1rewrite}). The above reformulations allow us to write $\mathcal{D} =   \mathcal{Q}_{1}^{*} + \mathcal{Q}_{2}^{*}+\mathcal{Q}_{3}^{*}+ \mathcal{R}_{4}+ \mathcal{R}_{5}+ P_{2} +  \mathcal{E}_{9}$ (see \eqref{eq:D_lb}), where $\mathcal{E}_{9} :=  \mathcal{E}_{1}-\mathcal{E}_{4} -\mathcal{E}_{5} -\mathcal{E}_{6}+  \mathcal{E}_{7} +\mathcal{E}_{8}$ and $\mathcal{Q}_{3}^{*} := \mathcal{Q}_{3}+ \mathcal{R}_{2}+ \mathcal{R}_{3}.$
Now, we obtain suitable bounds for $P_{2}$ (see \eqref{eq:P2bound}) and $\mathcal{E}_{1}$ (see \ref{eq:E1bd}). Further, through Lemmas  \ref{lemma:Q1rewrite} and \ref{lemma:Q2}, we show that the terms $\mathcal{E}_{4},\mathcal{E}_{5},\mathcal{E}_{6}$ and $\mathcal{E}_{7}$ are bounded. Lemma \ref{lemma:Q2} also provides a lower bound for $\mathcal{Q}_1^{*}+ \mathcal{Q}_2^{*}.$ In  \eqref{eq:R4_lb} and \eqref{eq:R5_lb}, we derive lower bounds for $\mathcal{R}_{4}\, \mbox{and}\, \mathcal{R}_{5},$ respectively.  Furthermore, in \eqref{eq:E9bound} we estimate the term $\mathcal{E}_{9}$, while in Lemma \ref{lemma:Q3pr_bd}, we establish a lower bound for $\mathcal{Q}_{3}^{*}$. Combining these results, we derive a lower bound for $\mathcal{D}.$ Using this lower bound alongside Lemma \ref{lemma:deltamin_deltapr_est}, we finally obtain a lower bound for the term $\sum_{j=l-1}^{r+1}\left((\delta_{\jmh})^2- (\delta_{j-1}^{\prime})^2\right),$ thereby completing the proof.

\subsection{Proof of Lemma \ref{lemma:osle}}
We present an auxiliary lemma for non-increasing sequences, which will be used in the proof of Lemma \ref{lemma:osle}.
\begin{lemma}\label{lemma:dec_seq_lemma}
Consider a non-increasing sequence $\{u_{j}\}_{\jinz},$ with jumps $\delta_{\jph}:= u_{j+1}- u_{j}.$  
Let $\{u^{\prime}_{\jph}\}_{\jinz}$ be obtained from $\{u_{j}\}_{\jinz}$ by applying the time-update formula \eqref{eq:NTscheme} and
     denote the corresponding jumps $\delta_{j}^{\prime}:= u_{\jph}^{\prime}- u_{\jmh}^{\prime}$ for $\jinz$. Under the CFL condition \eqref{eq:cfl},
    the jump sequence $\{\delta_{j}^{\prime}\}_{\jinz}$ satisfies the following estimate    \begin{nalign}\label{eq:decseq_bd}
    \left(\delta_{j}^{\prime}\right)_{+} \leq \lambda \norm{f_{k}}\left(\abs{\Delta k_{\jph}}+\abs{\Delta k_{\jmh}}\right).
\end{nalign}
\end{lemma}
\begin{proof} As  $\delta_{\jph} \leq 0$ for $\jinz,$ applying the CFL condition \eqref{eq:cfl} to \eqref{eq:deltaprime}, it follows that 
\begin{nalign}\label{eq:deltajnegbd}
    \delta_{j}^{\prime} &\leq \left(\frac{1}{2}-\frac{1}{8}-\kappa-\kappa^2\right)(\delta_{\jph}+ \delta_{\jmh}) + \lambda \norm{f_{k}}\left(\abs{\Delta k_{\jph}}+\abs{\Delta k_{\jmh}}\right)\\
&\leq  \lambda \norm{f_{k}}\left(\abs{\Delta k_{\jph}}+\abs{\Delta k_{\jmh}}\right),
\end{nalign}
from which the estimate \eqref{eq:decseq_bd} follows directly.
\end{proof}

\begin{proof}[Proof of Lemma \ref{lemma:osle}] Denote the jumps $\delta_{j}^{n+1}:= u_{\jph}^{n+1}- u_{\jmh}^{n+1}$ and $\delta_{\jmh}^{n}:= u_{j}^{n}- u_{j-1}^{n}$ for $\jinz$. The key step in deriving the estimate \eqref{eq:oslsfinal} involves  decomposing $\{u_{j}^{n}\}_{\jinz}$ into monotone sequences and invoking Lemmas \ref{lemma:mainlemma_inseq} and \ref{lemma:dec_seq_lemma} to obtain useful estimates.
To this end, we split the index set $\mathbb{Z}$ of the sequence $\{u_{j}^{n}\}_{\jinz}$ into maximal subsets of the form $\Gamma_m:=\{l-1,l,l+1,\cdots, r\}, m \in S\subseteq \mathbb{N}$ such that $\{u_j^n\}_{j \in \Gamma_m}$ is either non-decreasing or non-increasing. Clearly, $\bigcup_{m \in S }\Gamma_{m} =\mathbb{Z}.$ 
Now, for each $m\in S,$ we define a new sequence $\{u^{m}_{j}\}_{\jinz},$ which is either non-decreasing or non-increasing, as follows:
\begin{nalign}\label{eq:corr_increa}
    u^{m}_{j} &:= \begin{cases}
        u_{j}^{n}, \quad \,\,\,\, \,\mbox{if} \,l-1 \leq  j \leq r,\\
        u_{l-1}^{n},  \quad  \mbox{if} \,\, j < l-1,\\
        u_{r}^{n},  \,\,\,\, \, \quad \mbox{if} \,\, j >  r.
    \end{cases}
\end{nalign}
With this definition, we observe that for each $m \in S,$ the jump sequence  $\{\delta_{\jph}^{m}\}_{\jinz}$ associated with $u^{m}_{j}$ is either non-negative or non-positive. Additionally, 
\begin{nalign}\label{eq:delta_n_square_sum}
\sum_{\jinz}(\delta^{n}_{\jph})_{+}^{2} = 
\sum_{m \in S}\sum_{\jinz}(\delta^{m}_{\jph})_{+}^{2} = \sum_{m \in S}\sum_{j \in \Gamma_{m}}(\delta^{m}_{\jph})_{+}^{2} . 
\end{nalign}

Now, our aim is to compare the sums $\sum_{\jinz}(\delta^{n+1}_{j})_{+}^{2}$ and $\sum_{m \in S}\sum_{\jinz}(\delta^{m^\prime}_{j})_{+}^{2},$ where for $m \in S,$ $\delta^{m^\prime}_{j}:=u^{m^{\prime}}_{\jph}-u^{m^{\prime}}_{\jmh},$ and $\{u^{m^{\prime}}_{\jph}\}_{\jinz}$ is obtained by applying the scheme \eqref{eq:NTscheme} to $\{u^{m}_{j}\}_{\jinz}.$
For this, first we note that the jumps $\{\delta_{j}^{n+1}\}_{\jinz}$ can be expressed as
\begin{nalign}\label{eq:jumpeqn}
    \delta_{j}^{n+1}
    & = \frac{1}{2}(\delta_{\jph}^n+ \delta_{\jmh}^{n}) -\lambda\left( f(k_{j+1}, u_{j+1}^{\nph})- 2 f(k_{j}, u_{j}^{\nph})+ f(k_{j-1}, u_{j-1}^{\nph}) \right)\\ & \spc - \frac{1}{8}(\sigma_{j+1}^{n}- 2\sigma_{j}^{n}+ \sigma_{j-1}^{n}).
\end{nalign}
  Next, consider a region where the sequence $\{u_{j}^{n}\}_{\jinz}$ is monotone, say for $j\in \{l-1, l, \dots, r\}.$ In this case it is evident that $\delta_{\jph}^{n}=\delta_{\jph}^{m}$ for $j \in \{l-1,l, \dots,r-1\},$ for some $m \in S.$  Consequently, from \eqref{eq:jumpeqn} we have $\delta_{j}^{n+1}=\delta_{j}^{m^{\prime}},$ for $j \in \{l,l+1, \dots, r-1\}.$ 
In other words, away from the local extremum points of $\{u_{j}^{n}\}_{\jinz},$ the relation $\delta_{j}^{n+1} = \delta_{j}^{m^{\prime}} $ holds  for some $m.$ 
On the other hand, at the local extremum points,       say $j\in\{l-1, r\},$ of $\{u_{j}^{n}\}_{\jinz}$  with $\{u_{j}^{n}
\}$ non-decreasing on $j\in\{l-1, l,  \dots, r\},$ we have
\begin{align}\label{eq:jumpatextrema}
\delta_{r}^{n+1} &= \delta_{r}^{m^{\prime}}+ \delta_{r}^{(m+1)^{\prime}} +\lambda f_{k}(\bar{k}_{r+\half}, u_{r}^{n})\Delta k_{r+\frac{1}{2}} - \lambda f_{k}(\bar{k}_{r-\half}, u_{r}^{n})\Delta k_{r-\frac{1}{2}}, \\\delta_{l-1}^{n+1} &= \delta_{l-1}^{(m-1)^{\prime}}+ \delta_{l-1}^{m^{\prime}} +\lambda f_{k}(\bar{k}_{l-\half}, u_{l-1}^{n})\Delta k_{l-\frac{1}{2}} - \lambda f_{k}(\bar{k}_{l-\frac{3}{2}}, u_{l-1}^{n})\Delta k_{l-\frac{3}{2}}, \notag
\end{align}
for $\bar{k}_{r-\half} \in \mathcal{I}(k_{r-1}, k_{r}),$ $\bar{k}_{r+\half} \in \mathcal{I}(k_{r}, k_{r+1}),$  $\bar{k}_{l-\half} \in \mathcal{I}(k_{l-1}, k_{l})$ and   $\bar{k}_{l-\frac{3}{2}} \in \mathcal{I}(k_{l-2},\\ k_{l-1}).$

Further, using Lemma \ref{lemma:dec_seq_lemma} for the term $(\delta_{r}^{(m+1)^{\prime}})_{+}$(which is generated from the non-increasing sequence $\{u^{m+1}\}_{\jinz}$) in \eqref{eq:jumpatextrema}, along with the property that  $a, b \in \mathbb{R},$ $(a+b)_{+} \leq a_+ + b_+,$ we obtain
\begin{align}\label{eq:deltar_plus_bd}
\left(\delta_{r}^{n+1}\right)_{+} &\leq  (\delta_{r}^{m^{\prime}})_{+}+ (\delta_{r}^{(m+1)^{\prime}})_{+} + \lambda\norm{f_{k}}(\abs{\Delta k_{r+\frac{1}{2}}}+ \abs{\Delta k_{r-\frac{1}{2}}})\\
& \leq (\delta_{r}^{m^{\prime}})_{+} + 2\lambda\norm{f_{k}}(\abs{\Delta k_{r+\frac{1}{2}}}+ \abs{\Delta k_{r-\frac{1}{2}}}). \notag
\end{align} 
Upon squaring both sides of \eqref{eq:deltar_plus_bd} and using the bound $\abs{\delta_{r}^{m^{\prime}}} \leq 2C_{u_{0}}$  (a consequence of Theorem \ref{thm:maxpri}), we deduce that
\begin{nalign}\label{eq:jumpnearext_final}
\left(\delta_{r}^{n+1}\right)_{+}^{2} 
& \leq  (\delta_{r}^{m^{ \prime}})_{+}^{2} + \tilde{\Psi}\left(\abs{\Delta k_{r+\frac{1}{2}}} + \abs{\Delta k_{r-\frac{1}{2}}}\right),
\end{nalign} where $\tilde{\Psi}:=16\lambda^{2}\norm{f_k}^{2}\norm{k}+ 8C_{u_{0}}\lambda\norm{f_k}.$ 
Analogous arguments for $j=l-1$ yield \begin{align}\label{eq:jumpnearmin}
    \left(\delta_{l-1}^{n+1}\right)_{+}^{2} 
& \leq  (\delta_{l-1}^{m^ \prime})_{+}^{2} + \tilde{\Psi}\left(\abs{\Delta k_{l-\frac{1}{2}}} + \abs{\Delta k_{l-\frac{3}{2}}}\right).
\end{align}
Thus, we observe that for any $\jinz,$ either $(\delta_{j}^{n+1})_{+}^{2}=(\delta_{j}^{m^\prime})_{+}^{2}$ or  $(\delta_{j}^{n+1})_{+}^{2} \leq (\delta_{j}^{m^{\prime}})_{+}^{2} + \tilde{\Psi}(\abs{\Delta k_{j+\frac{1}{2}}}+ \abs{\Delta k_{j-\frac{1}{2}}}),$  for some $m \in S.$ Therefore, \begin{align}\label{eq:osls_peniltimate}
\sum_{\jinz}\left(\delta_{j}^{n+1}\right)_{+}^{2}& \leq \sum_{m \in S}\sum_{j \in \Gamma_{m}}\left(\delta_{j}^{m^{\prime}}\right)_{+}^{2} + 2\tilde{\Psi}\norm{k}_{BV} \\ &= \sum_{m \in S_{\uparrow}}\sum_{j \in \Gamma_{m}}\left(\delta_{j}^{m^{\prime}}\right)_{+}^{2} + \sum_{m \in S_{\downarrow}}\sum_{j \in \Gamma_{m}}\left(\delta_{j}^{m^{\prime}}\right)_{+}^{2}+ 2\tilde{\Psi}\norm{k}_{BV}    \notag
\end{align} where $S_{\uparrow}:=\{m \in S: \{u^{m}_{j}\}_{\jinz} \, \mbox{is a non-decreasing sequence}\} $ and $S_{\downarrow}:=\{m \in S : \{u^{m}_{j}\}_{\jinz}\mbox{is a non-increasing sequence}\}.$ Finally, invoking Lemma \ref{lemma:mainlemma_inseq} for $\{u^{m}_{j}\}_{\jinz},\, m \in S_{\uparrow}$ and Lemma \ref{lemma:dec_seq_lemma} for  $\{u^{m}_{j}\}_{\jinz},\, m \in S_{\downarrow}$ in \eqref{eq:osls_peniltimate} and recollecting \eqref{eq:delta_n_square_sum}, we obtain 
\begin{nalign}\label{eq:oslsfinal_a}
\sum_{\jinz}\left(\delta_{j}^{n+1}\right)_{+}^{2}
    & \leq \sum_{m \in S_{\uparrow}}\sum_{j \in \Gamma_{m}}(\delta_{\jmh}^{m})_{+}^2- \frac{1}{500} \lambda \gamma_{1}\sum_{m \in S_{\uparrow}}\sum_{j \in \Gamma_{m}} (\delta_{\jmh}^{m})_{+}^{3}\\& \spc -\frac{1}{6400}\sum_{m \in S_{\uparrow}}\sum_{j \in \Gamma_{m}}\left(\Delta^{2}\delta_{\jmtbt}^{m}\right)^{2}  +3\Theta\norm{k}_{BV} \\ & \spc + \sum_{m \in S_{\downarrow}}\sum_{j \in \Gamma_{m}} \lambda^{2} \norm{f_{k}}^{2}\left(\abs{\Delta k_{\jph}}+\abs{\Delta k_{\jmh}}\right)^{2} + 2\tilde{\Psi}\norm{k}_{BV}\\& \leq \sum_{\jinz}(\delta_{\jmh}^{n})_{+}^2- \frac{1}{500} \lambda \gamma_{1}\sum_{\jinz} (\delta_{\jmh}^{n})_{+}^{3} + 3\Theta\norm{k}_{BV} + 8\norm{k}\lambda^{2}\norm{f_k}^{2}\norm{k}_{BV}\\ & \spc + 2\tilde{\Psi}\norm{k}_{BV}\\ & \leq \sum_{\jinz}(\delta_{\jmh}^{n})_{+}^2- \frac{1}{500} \lambda \gamma_{1}\sum_{\jinz} (\delta_{\jmh}^{n})_{+}^{3} + \Psi\norm{k}_{BV}, 
\end{nalign}
where $\Psi := 3\Theta +8\norm{k}\lambda^{2}\norm{f_k}^{2}+2\tilde{\Psi}.$
Here, we have used the fact that $$\sum_{\jinz}\left(\abs{\Delta k_{\jph}}+\abs{\Delta k_{\jmh}}\right)^{2} \leq 4\norm{k}\sum_{\jinz}\left(\abs{\Delta k_{\jph}}+\abs{\Delta k_{\jmh}}\right)\leq 8\norm{k}\norm{k}_{BV}.$$ This completes the proof.
\end{proof}

	\bibliographystyle{siamplain}
	\bibliography{references}
\end{document}

%% file: ex_shared.tex
\usepackage{lipsum}
\usepackage{amsfonts}
\usepackage{graphicx}
\usepackage{epstopdf}
\usepackage{algorithmic}
\ifpdf
  \DeclareGraphicsExtensions{.eps,.pdf,.png,.jpg}
\else
  \DeclareGraphicsExtensions{.eps}
\fi

\newsiamremark{remark}{Remark}
\newsiamremark{hypothesis}{Hypothesis}
\crefname{hypothesis}{Hypothesis}{Hypotheses}
\newsiamthm{claim}{Claim}

\headers{Conservation Laws with Discontinuous Flux}{NIKHIL MANOJ  AND SUDARSHAN KUMAR K.}

\title{Analysis of a central MUSCL-type scheme for conservation laws with discontinuous flux}

\author{Nikhil Manoj \orcidlink{0009-0003-6812-8633}\thanks{School of Mathematics\\
Indian Institute of Science Education and Research\\
Thiruvananthapuram --695551\\}
\and Sudarshan~Kumar~K. \orcidlink{0009-0005-1186-4167}
\thanks{School of Mathematics\\
Indian Institute of Science Education and Research\\
Thiruvananthapuram --695551\\
(\email{nikhilmanoj2020@iisertvm.ac.in}, \email{sudarshan@iisertvm.ac.in}).}
}

\usepackage{amsopn}
